\documentclass[preprint,12pt,3p,times,numbers,sort&compress]{elsarticle}

\usepackage{amssymb}
\usepackage{booktabs}
\usepackage{amsmath}
\usepackage{commath}
\usepackage{import}
\usepackage{tikz}
\usetikzlibrary{arrows.meta,decorations.markings}
\usepackage{epstopdf}
\usepackage{hyperref}
\usepackage{cleveref}
\usepackage{algorithm}
\usepackage{algpseudocode}
\usepackage{graphicx} 
\usepackage{array}
\usepackage{caption}
\usepackage{float}
\usepackage{placeins}
\usepackage{url}

\journal{Computers and Structures}

\definecolor{blue2}{RGB}{25, 132, 197}
\definecolor{red2}{RGB}{194, 55, 40}
\definecolor{gray2}{RGB}{204, 204, 204}

\newcommand{\genvec}{\raise0.5pt\hbox{\mathsurround=0pt$\mbox{\rule{0.1pt}{0pt}\tiny{}$\boldsymbol{\square}$\rule{0.2pt}{0pt}}$}}

\renewcommand{\url}{\mathbf{u}_\mathrm{r}^{j}}

\newcommand{\zs}{z_\mathrm{s}}

\begin{document}

\begin{frontmatter}

\title{Reduced order computation of 2D elastodynamic Green's functions in layered soil using a low-rank tensor approximation}

\author{Z. Farooq\corref{cor1}}
\ead{zainab.farooq@kuleuven.be}
\author{A. Pashov}
\author{P. Reumers}
\author{S. Fran\c{c}ois}
\author{G. Degrande}

\address{KU Leuven, Department of Civil Engineering, Structural Mechanics Section, \\Kasteelpark Arenberg 40, B-3001 Leuven, Belgium}

\cortext[cor1]{Corresponding author}

\begin{abstract}
The evaluation of elastodynamic Green's functions across numerous source–receiver locations, frequencies, and material properties, particularly in the context of parametric studies or boundary element computations, is computationally demanding and memory intensive. This paper presents a reduced order modeling strategy based on the Greedy Tucker Approximation (GTA), which incrementally constructs a low-rank representation of the Green's tensor through rank-one enrichments obtained via a Proper Generalized Decomposition (PGD)-type alternating least squares procedure. A Petrov-Galerkin formulation is employed to improve convergence and approximation accuracy. The resulting multi-dimensional tensor, expressed in terms of one-dimensional basis functions and a compact core, achieves substantial reductions in memory requirements. The methodology is demonstrated for two cases: a soil layer on rigid bedrock and a layered halfspace. Different separable dimensions are considered to capture various combinations of source and receiver configurations, frequencies, and material parameters. Results are validated against those obtained with the direct stiffness method and computation times and memory requirements are compared.
\end{abstract}

\end{frontmatter}

\vspace{-1em}
%
%
\section{Introduction}
\label{sec:intro}
Elastodynamic Green's functions, fundamental solutions to the equations governing wave propagation in elastic media \cite{bedf23a}, are central to vibration analysis \cite{chop12a}, seismic modeling \cite{kaus06a}, and form the basis of the boundary element (BE) formulation \cite{bonn95a}. In applications such as railway-induced ground vibration, the dynamic response spans a wide frequency range and requires high spatial resolution to accurately capture wave propagation. Solving such problems using a BE approach requires Green's functions to be evaluated over dense spatial and frequency grids, covering many source–receiver pairs and multiple frequencies, which is both computationally intensive and memory-demanding. This results in a high-dimensional Green's function tensor, making its evaluation and storage a major bottleneck in large-scale simulations. While translational invariance in the horizontal direction can be exploited through a shifting property to reduce the number of required evaluations, more complex problems involving embedded foundations in a horizontally layered halfspace necessitate dense resolution along the receiver depth. This leads to a higher dimensionality of the Green's tensor and increased computational cost.
\par
In BE applications, closed-form Green's functions for a homogeneous full space are often used due to their analytical simplicity. However, this assumption becomes inadequate for more complex scenarios, such as a layered halfspace, as it requires meshing the free surface and layer interfaces. In such cases, numerical methods such as the direct stiffness method \citep{kaus81a} are employed to compute Green’s functions with greater accuracy. Despite efficient implementations \cite{sche09a}, evaluating these functions over large spatial grids and a wide frequency range remains computationally expensive and memory-intensive, since Green's functions must be repeatedly evaluated for numerous source–receiver pairs and frequencies.
\par
Model order reduction (MOR) techniques \cite{benn17a} have emerged as powerful tools for reducing the computational demands of high-fidelity simulations in engineering and applied sciences. These methods rely on the idea that the response of a complex, computationally intensive model can be efficiently approximated by a surrogate model. Such a reduced model is constructed by projecting the original system onto a low-dimensional subspace, thereby significantly lowering computational cost while preserving essential system behavior. MOR approaches have been successfully applied for fluid mechanics \cite{lass14a}, aerodynamics \cite{mend19a}, vibroacoustics \cite{deck20a, phd-wall-18a} and more diverse fields such as real-time thermal process monitoring \cite{agua15a}. These approaches have enabled efficient analyses of wave propagation \cite{pere08a, reko25a}, structural vibration \cite{bess13a, lu21a, luo24a}, and parametric studies \cite{benn15a, jens20a, kim23a}. For instance, reduced order models of harmonic Green’s functions in layered soils have been developed to mitigate the prohibitive cost of large-scale evaluations \cite{othm20a}, while projection-based approaches have been employed to accelerate seismic response analyses and vibration problems with parametric variability \cite{daby21a, lenz23a, hawk25a, bame17a}.

In this paper, MOR is employed to construct low-rank approximations of elastodynamic Green's functions using Proper Generalized Decomposition (PGD) \cite{chin13a, nouy10a} and the Greedy Tucker Approximation (GTA) algorithm  \cite{geor19a}. PGD has emerged as a versatile approach for solving high-dimensional problems by exploiting separability in space, time, and parameters. Its foundations and scope have been reviewed in \cite{chin11a} and it has since been successfully applied to a broad spectrum of static and dynamic problems. Applications include harbor agitation \cite{mode15a}, seismic response analysis \cite{germ16a}, structural dynamics \cite{chev11a}, and transient wave propagation \cite{gout21a}. Beyond elastodynamics, PGD has also been used in multi-scale multi-physics simulations \cite{nero10a}, topology optimization \cite{pauw24a}, non-linear material modeling \cite{vits19a}, and stochastic analysis \cite{louf13a}. More recently, GTA has been introduced by Georgieva and Hofreither  \cite{geor19a} as a reduction strategy tailored to high-dimensional tensor problems. By combining successive rank-one enrichments with a Tucker projection, it constructs low-rank representations with flexible rank adaptation across modes. Although still relatively new, GTA offers promising opportunities for reduced order modeling of elastodynamic wave propagation. This paper explores that potential through the construction of low-rank Green's function approximations.

Here, a Petrov-Galerkin formulation is employed to enhance numerical stability and convergence. These techniques construct the solution a priori as a finite sum of separable functions, each depending on a single physical or parametric variable. This separable structure enables an efficient tensor decomposition \cite{kold09a} of the Green’s function. PGD yields a Canonical Polyadic Decomposition (CPD), where the solution is represented as a sum of rank-one tensors, while GTA produces a Tucker format, with a small core tensor and low-dimensional basis functions in each direction. In this way, the so-called curse of dimensionality is circumvented, making the approach particularly advantageous for problems involving multiple parameters and high-resolution grids. Compared to traditional dense representations, these formats offer substantial efficiency gains, particularly when repeated evaluations over multiple parameters are required.  
\par
The paper is organized as follows. Section \ref{sec:methodology} presents the problem formulation, including the governing equations, boundary conditions, and the derivation of the weak form.  Section \ref{sec:mor} discusses in detail the MOR strategy used to approximate the Green’s functions in low-rank tensor formats. Section \ref{sec:examples} demonstrates the proposed methodology through numerical examples. Green’s functions are first computed for a soil layer on rigid bedrock, with varying separable dimensions to illustrate the method’s ability to capture accurate approximations across different combinations of source and receiver configurations, frequencies, and material parameters. This is followed by a realistic application involving the computation of Green’s functions for a layered halfspace representative of the Groene Hart site in the Netherlands. The computational performance of the proposed algorithm is compared with results obtained using the direct stiffness method implemented in the Elastodynamic Toolbox (EDT) \cite{sche09a}, and evaluated in terms of accuracy, memory usage, and computation time, followed by a discussion on the merits and limitations of the formulation. Finally, section \ref{sec:conclusions} summarises the main conclusions.

\section{Methodology}
\label{sec:methodology}
In this paper, two-dimensional wave propagation problems in the $(x,z)$-plane are considered. For such problems, the wave motion decouples into in-plane motion (P-SV) and out-of-plane (SH) components. The governing equations are presented and weak formulations for both the P–SV and SH problems are derived based on these equations and the corresponding boundary conditions.

\subsection{Governing equations}
Consider the two-dimensional, horizontally layered soil domain $\Omega_\mathrm{S}$ of infinite extent, as shown in figure~\ref{fig:halfspace}. The layered soil may rest on either a rigid bedrock (figure~\ref{fig:halfspace}a) or a homogeneous halfspace (figure~\ref{fig:halfspace}b), representing typical boundary conditions in site-response analysis. 
The soil domain $\Omega_\mathrm{S}$ is excited by a concentrated pulse load $\rho\hat{\mathbf{b}}(\mathbf{x},\omega)=\delta(\mathbf{x}-\mathbf{x}_\mathrm{s})\mathbf{e}_i$ applied at location~$\mathbf{x}_\mathrm{s}$ in direction~$\mathbf{e}_i$, where $\hat{\mathbf{b}}$ represents body forces, $\rho$ is the material density, $\omega$ is the excitation frequency and $\delta(\cdot)$ denotes the Dirac delta function.  A hat above a variable denotes its representation in the frequency domain. 
The corresponding displacement field $\hat{\mathbf{u}}(\mathbf{x},\omega)$ evaluated at location $\mathbf{x}$ is called a Green's function or fundamental solution of the layered soil domain. Since the soil domain is invariant in the $x$-direction, the concentrated pulse load is applied at $\mathbf{x}_\mathrm{s}=\{0,\zs\}^\mathrm{T}$ resulting in body forces $\rho\hat{\mathbf{b}}(\mathbf{x},\omega)=\delta(x)\delta(z-\zs)\mathbf{e}_i$,  The elastodynamic equilibrium equations in the frequency domain are:

\begin{figure}[!htb]
  \centering
  (a) \includegraphics[width=0.35\linewidth]{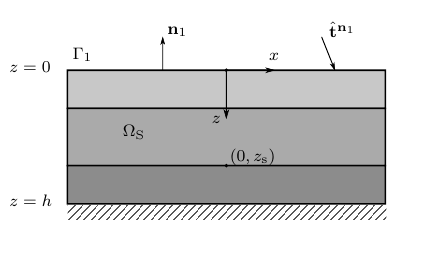}
  (b) \includegraphics[width=0.35\linewidth]{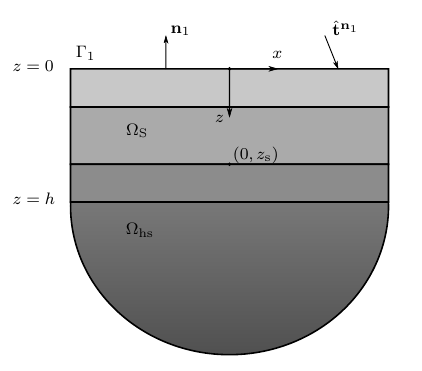}
  \caption{Horizontally layered soil medium overlying (a) a rigid bedrock and (b) a homogeneous halfspace. A concentrated pulse load $\rho\hat{\mathbf{b}}(\mathbf{x},\omega)=\delta(x)\delta(z-z_\mathrm{s})\mathbf{e}_i$ is applied at $\mathbf{x}_\mathrm{s}=\{0,\,z_\mathrm{s}\}^\mathrm{T}$ in direction $\mathbf{e}_i$.}
  \label{fig:halfspace}
\end{figure}

\begin{equation}
	\mathbf{L}^\mathrm{T}\hat{\boldsymbol{\sigma}}+\rho\hat{\mathbf{b}}=-\rho\omega^2\hat{\mathbf{u}}\quad\text{in}\quad\Omega_\mathrm{S}
	\label{eq:equi}
\end{equation}
where $\hat{\boldsymbol{\sigma}}$ denotes the stress vector. The matrix $\mathbf{L}$ contains differential operators and varies depending on whether P-SV or SH wave motion is considered. 
For linear elastic materials, the relation between the stress vector~$\hat{\boldsymbol{\sigma}}$ and strain vector~$\hat{\boldsymbol{\epsilon}}$ is:
\begin{equation}
	\hat{\boldsymbol{\sigma}}=\mathbf{C}\hat{\boldsymbol{\epsilon}}
	\label{eq:const}
\end{equation}
The strain-displacement relation is given by:
\begin{equation}
	\hat{\boldsymbol{\epsilon}}=\mathbf{L}\mathbf{u}
	\label{eq:kine}
\end{equation}
\par
The P-SV problem describes in-plane motion due to in-plane wave propagation. The out-of-plane displacement component $\hat{u}_y$ is equal to zero and the displacement vector reduces to ${\hat{\mathbf{u}}=\{\hat{u}_x,\hat{u}_z\}^\mathrm{T}}$. Due to the two-dimensional nature of the problem, the dependence of all variables on the coordinate~$y$ vanishes, resulting in zero values for the strains $\hat{\epsilon}_{yy}$, $\hat{\gamma}_{xy}$, and $\hat{\gamma}_{yz}$. The remaining strain components are collected in the strain vector $\hat{\boldsymbol{\epsilon}}=\{\hat{\epsilon}_{xx},\hat{\epsilon}_{zz},\hat{\gamma}_{zx}\}^\mathrm{T}$, where the operator $\mathbf{L}$ is defined as:
\begin{equation}
	\mathbf{L}=\mathbf{L}_x\dpd{}{x}+\mathbf{L}_z\dpd{}{z}=\begin{bmatrix}
	1 & 0 \\
	0 & 0 \\
	0 & 1
	\end{bmatrix}\dpd{}{x}+\begin{bmatrix}
	0 & 0 \\
	0 & 1 \\
	1 & 0
	\end{bmatrix}\dpd{}{z}
\end{equation}
The constitutive matrix $\mathbf{C}$, which relates the stress vector $\hat{\boldsymbol{\sigma}}=\{\hat{\sigma}_{xx},\hat{\sigma}_{zz},\hat{\sigma}_{zx}\}^\mathrm{T}$ to the strain vector~$\hat{\boldsymbol{\epsilon}}$, equals: 
\begin{equation}
	\mathbf{C}=\begin{bmatrix}
	\lambda+2\mu & \lambda & 0 \\
	\lambda & \lambda+2\mu & 0 \\
	0 & 0 & \mu
	\end{bmatrix}
\end{equation}
where $\lambda$ denotes the first Lam\'e constant and $\mu$ is the shear modulus.
\par
In the SH problem, the motion is purely out-of-plane and the in-plane displacement components $\hat{u}_x$ and $\hat{u}_z$ are equal to zero. The displacement vector reduces to $\hat{\mathbf{u}}=\{{\hat{u}_y\}^\mathrm{T}}$. In this case, only the shear strains remain non-zero, reducing the strain vector to $\hat{\boldsymbol{\epsilon}}=\{\hat{\gamma}_{xy},\hat{\gamma}_{zy}\}^\mathrm{T}$. The corresponding stress vector becomes $\hat{\boldsymbol{\sigma}}=\{\hat{\sigma}_{xy},\hat{\sigma}_{zy}\}^\mathrm{T}$. The differential operator $\mathbf{L}$ for SH case is defined as:
\begin{equation}
	\mathbf{L}=\mathbf{L}_x\dpd{}{x}+\mathbf{L}_z\dpd{}{z}=\begin{bmatrix}
	1 \\
	0
	\end{bmatrix}\dpd{}{x}+\begin{bmatrix}
	0 \\
	1
	\end{bmatrix}\dpd{}{z}
\end{equation}
The corresponding constitutive matrix $\mathbf{C}$ equals: 
\begin{equation}
	\mathbf{C}=\begin{bmatrix}
	\mu & 0 \\
	0 & \mu 
	\end{bmatrix}
\end{equation}

\subsection{Boundary conditions}
\begin{figure}[H]
  \centering
  \includegraphics[width=0.35\linewidth]{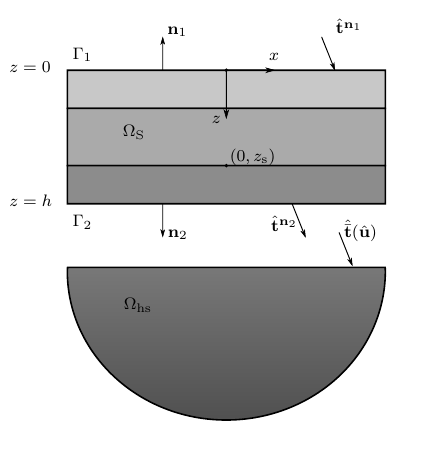}
  \caption{Boundary conditions for a horizontally layered soil medium over a homogeneous halfspace. }
  \label{fig:halfspace2}
\end{figure}

The response of a layered halfspace can be represented using different numerical formulations, depending on how the soil layers and the underlying medium are modeled, thereby defining the corresponding boundary conditions. A hybrid formulation is adopted, coupling thin layer elements for the soil layers with a halfspace element to model a layered halfspace efficiently. The corresponding formulation is detailed in subsection~\ref{subsec:stiffness}. The layered soil medium in figure~\ref{fig:halfspace} is infinite and invariant in the $x$-direction. Wave propagation towards infinity in the $x$-direction is accounted for by the transformation of the coordinate $x$ to the wavenumber domain $k_x$ (subsection~\ref{subsec:transform}). The free surface of the soil domain is denoted by the boundary $\Gamma_1$ (figure~\ref{fig:halfspace2}) on which the following boundary condition is imposed: 
\begin{equation}
	\hat{\mathbf{t}}^{\mathbf{n}_1}=\hat{\boldsymbol{\sigma}}\cdot\mathbf{n}_1=\mathbf{0}\quad\text{on}\quad\Gamma_1
	\label{eq:BC1}
\end{equation}
with $\hat{\mathbf{t}}^{\mathbf{n}_1}$ the traction vector on $\Gamma_1$ with outward unit normal vector $\mathbf{n}_1$. On the horizontal boundary $\Gamma_2$ at $z=h$ (figure~\ref{fig:halfspace}a), the following Dirichlet boundary condition is imposed if the soil domain rests on bedrock: 
\begin{equation}
	\hat{\mathbf{u}}=\mathbf{0}\quad\text{on}\quad\Gamma_2
	\label{eq:BC2}
\end{equation}

Alternatively, if the soil domain overlays a homogeneous halfspace (figure~\ref{fig:halfspace2}), the following boundary condition is imposed on $\Gamma_2$: 
\begin{equation}
	\hat{\mathbf{t}}^{\mathbf{n}_2}=\hat{\boldsymbol{\sigma}}\cdot\mathbf{n}_2=-\hat{\bar{\mathbf{t}}}(\hat{\mathbf{u}})
	\label{eq:BC3}
\end{equation}
\noindent
where $\hat{\mathbf{t}}^{\mathbf{n}_2}$ is the traction vector at $z=h$ with outward unit normal vector $\mathbf{n}_2$. This boundary condition expresses equilibrium between the traction vector $\hat{\mathbf{t}}^{\mathbf{n}_2}$ on the boundary $\Gamma_2$ and the traction vector $\hat{\bar{\mathbf{t}}}(\hat{\mathbf{u}})$ on the surface of the homogeneous halfspace. 

\subsection{Variational formulation}
The governing equilibrium equation, together with the natural boundary conditions, is reformulated in the frequency–spatial domain using the weighted residual method. A kinematically admissible virtual displacement field $\hat{\mathbf{v}}$, satisfying the Dirichlet boundary conditions, is employed as the test function. This leads to the strong form of equations~\eqref{eq:equi}, \eqref{eq:BC1}, and \eqref{eq:BC3} in the spatial-frequency domain:

\begin{equation}
	-\int\limits_{\Omega_\mathrm{S}} \hat{\mathbf{v}}^\mathrm{T}\left(\mathbf{L}^\mathrm{T}\hat{\boldsymbol{\sigma}}+\rho\hat{\mathbf{b}}+\rho\omega^2\hat{\mathbf{u}}\right)\mathrm{d}\Omega
	+\int\limits_{\Gamma_1}\hat{\mathbf{v}}^\mathrm{T}\hat{\mathbf{t}}^{\mathbf{n}_1}\mathrm{d}\Gamma
	+\int\limits_{\Gamma_2}\hat{\mathbf{v}}^\mathrm{T}\left(\hat{\mathbf{t}}^{\mathbf{n}_2}+\hat{\bar{\mathbf{t}}}(\hat{\mathbf{u}}))\right)\mathrm{d}\Gamma=0.
	\label{eq:strong}
\end{equation}
The third term in equation~\eqref{eq:strong} vanishes when the Dirichlet boundary condition in equation~\eqref{eq:BC2} is imposed instead of the boundary condition specified in equation~\eqref{eq:BC3}. Integration by parts of equation~\eqref{eq:strong} and application of the divergence theorem yields the weak form: 
\begin{equation}
	\int\limits_{\Omega_\mathrm{S}}(\mathbf{L}\hat{\mathbf{v}})^\mathrm{T}\mathbf{C}\mathbf{L}\hat{\mathbf{u}}\mathrm{d}\Omega
	-\omega^2\int\limits_{\Omega_\mathrm{S}}\rho\hat{\mathbf{v}}^\mathrm{T}\hat{\mathbf{u}}\mathrm{d}\Omega 
	+\int\limits_{\Gamma_2}\hat{\mathbf{v}}^\mathrm{T} \, \hat{\bar{\mathbf{t}}}(\hat{\mathbf{u}})\mathrm{d}\Gamma
	=\int\limits_{\Omega_\mathrm{S}}\hat{\mathbf{v}}^\mathrm{T}\rho\hat{\mathbf{b}}\mathrm{d}\Omega.
	\label{eq:weak}
\end{equation}
Substituting $\mathbf{L}=\mathbf{L}_x\dpd{}{x}+\mathbf{L}_z\dpd{}{z}$ in equation~\eqref{eq:weak} results in: 
\begin{multline}
	\int\limits_{\Omega_\mathrm{S}}\left(\dpd{\hat{\mathbf{v}}}{x}\right)^\mathrm{T}\mathbf{L}_x^\mathrm{T}\mathbf{C}\mathbf{L}_x\dpd{\hat{\mathbf{u}}}{x}\mathrm{d}\Omega
	+\int\limits_{\Omega_\mathrm{S}}\left(\dpd{\hat{\mathbf{v}}}{x}\right)^\mathrm{T}\mathbf{L}_x^\mathrm{T}\mathbf{C}\mathbf{L}_z\dpd{\hat{\mathbf{u}}}{z}\mathrm{d}\Omega 
	+\int\limits_{\Omega_\mathrm{S}}\left(\dpd{\hat{\mathbf{v}}}{z}\right)^\mathrm{T}\mathbf{L}_z^\mathrm{T}\mathbf{C}\mathbf{L}_x\dpd{\hat{\mathbf{u}}}{x}\mathrm{d}\Omega \\
	+\int\limits_{\Omega_\mathrm{S}}\left(\dpd{\hat{\mathbf{v}}}{z}\right)^\mathrm{T}\mathbf{L}_z^\mathrm{T}\mathbf{C}\mathbf{L}_z\dpd{\hat{\mathbf{u}}}{z}\mathrm{d}\Omega
	-\omega^2\int\limits_{\Omega_\mathrm{S}}\rho\hat{\mathbf{v}}^\mathrm{T}\hat{\mathbf{u}}\mathrm{d}\Omega
	+\int\limits_{\Gamma_2}\hat{\mathbf{v}}^\mathrm{T} \,\hat{\bar{\mathbf{t}}}(\hat{\mathbf{u}})\mathrm{d}\Gamma
	=\int\limits_{\Omega_\mathrm{S}}\hat{\mathbf{v}}^\mathrm{T}\rho\hat{\mathbf{b}}\mathrm{d}\Omega
	\label{eq:weak2}
\end{multline}
The first two terms in equation~\eqref{eq:weak2} containing $\dpd{\hat{\mathbf{v}}}{x}$ are integrated by parts: 
\begin{multline}
	-\int\limits_{\Omega_\mathrm{S}}\hat{\mathbf{v}}^\mathrm{T}\mathbf{L}_x^\mathrm{T}\mathbf{C}\mathbf{L}_x\dpd[2]{\hat{\mathbf{u}}}{x}\mathrm{d}\Omega
	-\int\limits_{\Omega_\mathrm{S}}\hat{\mathbf{v}}^\mathrm{T}\mathbf{L}_x^\mathrm{T}\mathbf{C}\mathbf{L}_z\dmd{\hat{\mathbf{u}}}{2}{x}{}{z}{}\mathrm{d}\Omega
	+\int\limits_{\Omega_\mathrm{S}}\left(\dpd{\hat{\mathbf{v}}}{z}\right)^\mathrm{T}\mathbf{L}_z^\mathrm{T}\mathbf{C}\mathbf{L}_x\dpd{\hat{\mathbf{u}}}{x}\mathrm{d}\Omega \\
	+\int\limits_{\Omega_\mathrm{S}}\left(\dpd{\hat{\mathbf{v}}}{z}\right)^\mathrm{T}\mathbf{L}_z^\mathrm{T}\mathbf{C}\mathbf{L}_z\dpd{\hat{\mathbf{u}}}{z}\mathrm{d}\Omega
	-\omega^2\int\limits_{\Omega_\mathrm{S}}\rho\hat{\mathbf{v}}^\mathrm{T}\hat{\mathbf{u}}\mathrm{d}\Omega
	+\int\limits_{\Gamma_2}\hat{\mathbf{v}}^\mathrm{T} \,\hat{\bar{\mathbf{t}}}(\hat{\mathbf{u}})\mathrm{d}\Gamma
	=\int\limits_{\Omega_\mathrm{S}}\hat{\mathbf{v}}^\mathrm{T}\rho\hat{\mathbf{b}}\mathrm{d}\Omega
	\label{eq:weak3}
\end{multline}
Equation~\eqref{eq:weak3} represents the weak form in the spatial-frequency domain, which is subsequently transformed to the wavenumber-frequency domain.

\subsection{Coordinate transformation}
\label{subsec:transform}
The integrals over $\Omega_\mathrm{S}$ in equation~\eqref{eq:weak3} are reduced to integrals over the coordinate $z$ by transforming the coordinate $x$ to the wavenumber $k_x$ using the forward Fourier transform defined as: 
\begin{equation}
	\tilde{\mathbf{u}}(k_x,z,\omega)=\int\limits_{-\infty}^{\infty} \hat{\mathbf{u}}(x,z,\omega)\exp(+\mathrm{i}k_xx)\mathrm{d}x.
	\label{eq:ff}
\end{equation}
A tilde above a variable denotes its representation in the wavenumber-frequency domain. The inverse Fourier transform to the frequency-spatial domain is defined as: 
\begin{equation}
	\hat{\mathbf{u}}(x,z,\omega)=\frac{1}{2\pi}\int\limits_{-\infty}^{\infty} \tilde{\mathbf{u}}(k_x,z,\omega)\exp(-\mathrm{i}k_xx)\mathrm{d}k_x.
	\label{eq:iff}
\end{equation}
The weak form in the frequency-spatial domain is obtained by applying the forward Fourier transform defined in equation~\eqref{eq:ff} to equation~\eqref{eq:weak3}: 
\begin{multline}
		k_x^2\int\limits_{0}^{h}\tilde{\mathbf{v}}^\mathrm{T}\mathbf{L}_x^\mathrm{T}\mathbf{C}\mathbf{L}_x\tilde{\mathbf{u}}\mathrm{d}z
		+\mathrm{i}k_x\int\limits_{0}^{h}\left(\tilde{\mathbf{v}}^\mathrm{T}\mathbf{L}_x^\mathrm{T}\mathbf{C}\mathbf{L}_z\dpd{\tilde{\mathbf{u}}}{z}-\left(\dpd{\tilde{\mathbf{v}}}{z}\right)^\mathrm{T}\mathbf{L}_z^\mathrm{T}\mathbf{C}\mathbf{L}_x\tilde{\mathbf{u}}\right)\mathrm{d}z\\
		+\int\limits_{0}^{h}\left(\dpd{\tilde{\mathbf{v}}}{z}\right)^\mathrm{T}\mathbf{L}_z^\mathrm{T}\mathbf{C}\mathbf{L}_z\dpd{\tilde{\mathbf{u}}}{z}\mathrm{d}z
		-\omega^2\int\limits_{0}^{h}\rho\tilde{\mathbf{v}}^\mathrm{T}\tilde{\mathbf{u}}\mathrm{d}z
		+\tilde{\mathbf{v}}^\mathrm{T}(z=h)\tilde{\bar{\mathbf{t}}}(\tilde{\mathbf{u}})(z=h)
		=\int\limits_{0}^{h}\tilde{\mathbf{v}}^\mathrm{T}\rho\tilde{\mathbf{b}}\mathrm{d}z
	\label{eq:weak4}
\end{multline}
The wavenumber $k_x$ is alternatively expressed in terms of the slowness $p_x = k_x / \omega$, allowing for frequency-independent sampling in the slowness domain. The inverse transformation from the frequency–wavenumber domain to the frequency–spatial domain remains valid and can still be performed using equation~\eqref{eq:iff} for a given excitation frequency~$\omega$.

\subsection{Dynamic impedance of layered soil}
\label{subsec:stiffness}\
The soil stiffness can be determined using either the direct stiffness method or the thin layer method proposed by Kausel and Ro\"esset~\cite{kaus81a}. The thin layer method is based on the use of polynomial shape functions to represent the vertical variation of displacements and tractions. Compared to the direct stiffness method, the thin layer method leads to mathematically more tractable stiffness matrices involving only polynomial functions instead of transcendental functions. Due to its approximative nature, the method requires a small thickness of the layer elements compared to the smallest wavelength. The element stiffness matrix can be expressed as:

\begin{equation}
\tilde{\mathbf{K}}^{e} 
= k_x^{2} \tilde{\mathbf{A}}^{e} 
+ k_x \tilde{\mathbf{B}}^{e} 
+ \tilde{\mathbf{G}}^{e} 
- \omega^{2} \tilde{\mathbf{M}}^{e}
\end{equation}
where the system matrices $\tilde{\mathbf{A}}^{e}$, $\tilde{\mathbf{B}}^{e}$, $\tilde{\mathbf{G}}^{e}$ and $\tilde{\mathbf{M}}^{e}$ are independent of the wavenumber $k_x$ and the frequency $\omega$. The analytical expressions, following the formulation of Kausel and Ro\"esset~\cite{kaus81a}, are presented in ~\ref{sec:appendix1}.

The traction vector $\tilde{\bar{\mathbf{t}}}(\tilde{\mathbf{u}})$ from equation~\eqref{eq:weak4} is related to the soil displacements~$\tilde{\mathbf{u}}$ on the boundary $\Gamma_2$ by the stiffness matrix~$\tilde{\mathbf{K}}_\mathrm{s}(p_x,\omega)=\omega\tilde{\mathbf{K}}_\mathrm{p}(p_x)$ of the homogeneous halfspace:

\begin{equation}
	\tilde{\bar{\mathbf{t}}}(\tilde{\mathbf{u}})=\tilde{\mathbf{K}}_\mathrm{s}(p_x,\omega)\tilde{\mathbf{u}}=\omega\tilde{\mathbf{K}}_\mathrm{p}(p_x).
	\label{eq:Ks}
\end{equation}
This separated form is required in the reduced order formulation that is elaborated in section~\ref{sec:mor}. The stiffness matrix $\tilde{\mathbf{K}}_\mathrm{s}(p_x,\omega)$ relates tractions to displacements at the surface of a homogeneous halfspace.  The analytical expression ~\cite{kaus81a} for $\tilde{\mathbf{K}}_\mathrm{s}(p_x,\omega)$ for the SH problem is:
\begin{equation}
	\tilde{\mathbf{K}}_\mathrm{s}(p_x,\omega)=\mu\omega\sqrt{p_x^2-\frac{1}{C_\mathrm{s}^2}} \quad \text{with} \quad \omega>0 \quad \text{and} \quad p_x\geq 0
	\label{eq:Ks_SH}
\end{equation}
and for the P-SV problem is: 
\begin{equation}
	\tilde{\mathbf{K}}_\mathrm{s}(p_x,\omega)=\frac{2\mu\omega}{\xi}		
	\begin{bmatrix}
	\sqrt{p_x^2-\frac{1}{C_\mathrm{p}^2}} 
	& p_x\left(1-\xi\right) \\
	p_x\left(1-\xi\right)
	& \sqrt{p_x^2-\frac{1}{C_\mathrm{s}^2}} 
	\end{bmatrix} \quad \text{with} \quad \omega>0 \quad \text{and} \quad p_x\geq 0.
	\label{eq:Ks_PSV}
\end{equation}
The coefficient $\xi$ is defined as:
\begin{equation}
	\xi=2C_\mathrm{s}^2\left(p_x^2-\sqrt{p_x^2-\frac{1}{C_\mathrm{p}^2}}\sqrt{p_x^2-\frac{1}{C_\mathrm{s}^2}}\right)
\end{equation}
where $C_\mathrm{s}$ and $C_\mathrm{p}$ denote the shear and dilatational wave velocities of the homogeneous halfspace, respectively. Note that in~\cite{kaus81a} the $z$-axis is pointing upwards instead of downwards, and that the displacement and stress components in the $z$-direction are made imaginary to obtain symmetric matrices. Hence, the stiffness matrices given in this paper deviate slightly from those in~\cite{kaus81a}. Introducing \eqref{eq:Ks} into \eqref{eq:weak4} yields the final weak form:
\begin{multline}
		p_x^2\omega^2\int\limits_{0}^{h}\tilde{\mathbf{v}}^\mathrm{T}\mathbf{L}_x^\mathrm{T}\mathbf{C}\mathbf{L}_x\tilde{\mathbf{u}}\mathrm{d}z
		+\mathrm{i}p_x\omega\int\limits_{0}^{h}\left(\tilde{\mathbf{v}}^\mathrm{T}\mathbf{L}_x^\mathrm{T}\mathbf{C}\mathbf{L}_z\dpd{\tilde{\mathbf{u}}}{z}-\left(\dpd{\tilde{\mathbf{v}}}{z}\right)^\mathrm{T}\mathbf{L}_z^\mathrm{T}\mathbf{C}\mathbf{L}_x\tilde{\mathbf{u}}\right)\mathrm{d}z\\
		+\int\limits_{0}^{h}\left(\dpd{\tilde{\mathbf{v}}}{z}\right)^\mathrm{T}\mathbf{L}_z^\mathrm{T}\mathbf{C}\mathbf{L}_z\dpd{\tilde{\mathbf{u}}}{z}\mathrm{d}z
		-\omega^2\int\limits_{0}^{h}\rho\tilde{\mathbf{v}}^\mathrm{T}\tilde{\mathbf{u}}\mathrm{d}z
		+\omega\tilde{\mathbf{v}}^\mathrm{T}(z=h)\tilde{\mathbf{K}}_\mathrm{p}(p_x)\tilde{\mathbf{u}}(z=h)
		=\int\limits_{0}^{h}\tilde{\mathbf{v}}^\mathrm{T}\rho\tilde{\mathbf{b}}\mathrm{d}z
	\label{eq:weak5}
\end{multline}

\section{Low-rank approximation of Green's functions}
\label{sec:mor}
To mitigate the computational complexity of high-dimensional problems, tensor decomposition techniques are employed to represent the solution in a compact low-rank format. As mentioned earlier, the Greedy Rank-one (PGD) \cite{chin13a, nouy10a} and GTA \cite{geor19a} algorithms are adopted as algebraic solvers and a priori model reduction strategies. 

%
\subsection{Model reduction approach}
\label{subsec:pgd_gta}
The PGD provides a flexible framework by expressing the solution as a sum of separable rank-one terms.
From a mathematical perspective, PGD can be interpreted as a Canonical Polyadic Decomposition (CPD) \cite{kold09a} of a high-dimensional tensor, in which the solution is approximated as a sum of outer products of mode-wise functions. Consider a problem defined in a space of dimension $\mathcal{D} \in I_{\theta_1} \times \ldots \times I_{\theta_\mathcal{D}}$.  The coordinates $\theta_d \, (d = 1, \ldots, \mathcal{D})$ denote any system variables related to space, time, wavenumber, frequency but can also include model parameters such as material and geometrical parameters. PGD yields an approximation of the unknown displacement field $\mathbf{u}(\theta_1, \ldots, \theta_{\mathcal{D}})$ in a rank-$R_u$ separated form:
\begin{equation}
        \mathbf{u}(\theta_1, \ldots, \theta_{\mathcal{D}}) \approx \sum_{r_u=1}^{R_u} \mathbf{u}_{r_u} \quad \text{with} \quad \mathbf{u}_{r_u} = \bigotimes_{d=1}^{\mathcal{D}} \mathbf{u}_{r_u}^{(\theta_d)}
\end{equation} 
where $\otimes$ denotes the Kronecker product. Both the number of terms $R_u$ and the functions $\mathbf{u}_{r_u}^{(\theta_d)}$ are unknown a priori. The solution is constructed by sequentially adding enrichments in a greedy manner. 
  \begin{equation}
        \mathbf{u}(\theta_1, \ldots, \theta_{\mathcal{D}}) \approx \sum_{r_u=1}^{R_u} \bigotimes_{d=1}^{\mathcal{D}} \mathbf{u}_{r_u}^{(\theta_d)} + \bigotimes_{d=1}^{\mathcal{D}} \mathbf{u}_{R_u+1}^{(\theta_d)}
    \end{equation}
The number of enrichments $R_u$ is the key factor which indicates the efficiency of the reduction.

The GTA algorithm builds directly on this rank-one enrichment strategy. In each greedy step, a new separable mode-wise vector is computed using an alternating least-squares (ALS) procedure identical in spirit to PGD and appended to the basis in each tensor mode. After orthonormalizing these mode-wise bases, the governing linear system is projected onto the resulting low-dimensional subspace and the Tucker core tensor $\boldsymbol{\mathcal{G}}$ is computed. Iterating this process yields a Tucker-form approximation of the Green’s displacement tensor $\tilde{\mathbf{u}}^{\mathrm{G}}$ a priori, without computing the full tensor explicitly.

\par
The weak form \eqref{eq:weak5} assumes a fixed value for the source depth $z_\mathrm{s}$, the slowness $p_x$, and the frequency $\omega$. By treating these quantities as additional coordinates in the low-rank approximation, the Green’s displacement is represented as the tensor $\tilde{\mathbf{u}}^{\mathrm{G}} (x_\mathrm{s}, z_\mathrm{s}, p_x, z, \omega)$. The weighted residual formulation is correspondingly extended to these coordinates, leading to the following high-dimensional weak form:
\begin{multline}
		\int\limits_{\Omega}p_x^2\omega^2\tilde{\mathbf{v}}^\mathrm{T}\mathbf{L}_x^\mathrm{T}\mathbf{C}\mathbf{L}_x\tilde{\mathbf{u}}\mathrm{d}\Omega
		+\int\limits_{\Omega}\mathrm{i}p_x\omega\left(\tilde{\mathbf{v}}^\mathrm{T}\mathbf{L}_x^\mathrm{T}\mathbf{C}\mathbf{L}_z\dpd{\tilde{\mathbf{u}}}{z}-\left(\dpd{\tilde{\mathbf{v}}}{z}\right)^\mathrm{T}\mathbf{L}_z^\mathrm{T}\mathbf{C}\mathbf{L}_x\tilde{\mathbf{u}}\right)\mathrm{d}\Omega
		+\int\limits_{\Omega}\left(\dpd{\tilde{\mathbf{v}}}{z}\right)^\mathrm{T}\mathbf{L}_z^\mathrm{T}\mathbf{C}\mathbf{L}_z\dpd{\tilde{\mathbf{u}}}{z}\mathrm{d}\Omega \\
		-\int\limits_{\Omega}\rho\omega^2\tilde{\mathbf{v}}^\mathrm{T}\tilde{\mathbf{u}}\mathrm{d}\Omega
		+\int\limits_{\mathcal{I}_{x_\mathrm{s}}} \int\limits_{\mathcal{I}_{z_\mathrm{s}}} \int\limits_{\mathcal{I}_{p_x}}\int\limits_{\mathcal{I}_\omega}\omega\tilde{\mathbf{v}}^\mathrm{T}(z=h)\mathbf{K}_\mathrm{p}(p_x)\tilde{\mathbf{u}}(z=h) \; \mathrm{d}\omega \, \mathrm{d}p_x \,\mathrm{d}z_\mathrm{s} \, \mathrm{d}x_\mathrm{s} 
		=\int\limits_{\Omega}\tilde{\mathbf{v}}^\mathrm{T}\rho\tilde{\mathbf{b}}\mathrm{d}\Omega
	\label{eq:weak6}
\end{multline}
where $\Omega  
= \mathcal{I}_{x_\mathrm{s}} \times 
\mathcal{I}_{z_\mathrm{s}} \times 
\mathcal{I}_{p_x} \times 
\mathcal{I}_z \times 
\mathcal{I}_\omega \subset \mathbb{R}^5$
is the five-dimensional domain of the problem. The soil impedance term is evaluated at $z=h$ and is therefore integrated only over the parametric dimensions.

\subsection{Separated representation of the Green’s function}
\label{subsec:gta}
The Green's displacement tensor is approximated by the following Tucker approximation:
\begin{multline}\label{eq:uPGD}
\tilde{\mathbf{u}}^{\mathrm{G}} (x_\mathrm{s}, z_\mathrm{s}, p_x, z, \omega) \approx \boldsymbol{\mathcal{U}}= \\
\sum_{r_{x_{\mathrm{s}}}=1}^{R_{x_{\mathrm{s}}}} 
\sum_{r_{z_{\mathrm{s}}}=1}^{R_{z_{\mathrm{s}}}} 
\sum_{r_{p_x}=1}^{R_{p_x}}
\sum_{r_z=1}^{R_z} 
\sum_{r_\omega=1}^{R_\omega}
\displaystyle g_{r_{x_{\mathrm{s}}} r_{z_{\mathrm{s}}} r_{p_x} r_z r_{\omega}} 
\mathbf{u}^{\mathrm{G}(x_{\mathrm{s}})}_{r_{x_{\mathrm{s}}}}(x_{\mathrm{s}}) \otimes 
\mathbf{u}^{\mathrm{G}(z_{\mathrm{s}})}_{r_{z_{\mathrm{s}}}}(z_{\mathrm{s}}) \otimes 
{\mathbf{u}}^{\mathrm{G}(p_x)}_{r_{p_x}}(p_x) \otimes
\mathbf{u}^{\mathrm{G}(z)}_{r_z}(z) \otimes
{\mathbf{u}}^{\mathrm{G}(\omega)}_{r_{\omega}}(\omega)
\end{multline}
where the problem coordinates include source positions along the $x$- and $z$-directions $(x_\mathrm{s}, z_\mathrm{s})$, slowness $p_x$, receiver positions  $z$ and frequency $\omega$. 
${R_{x_{\mathrm{s}}}}, {R_{z_{\mathrm{s}}}}, {R_{p_x}}, {R_z}$ and ${R_\omega}$ specify the Tucker ranks associated with each dimension of the approximation. 
$\displaystyle g_{r_{x_{\mathrm{s}}} r_{z_{\mathrm{s}}} r_{p_x} r_z r_{\omega}}$ is the core tensor and  
$\mathbf{u}^{\mathrm{G}(x_{\mathrm{s}})}_{r_{x_{\mathrm{s}}}}(x_{\mathrm{s}})$, 
$\mathbf{u}^{\mathrm{G}(z_{\mathrm{s}})}_{r_{z_{\mathrm{s}}}}(z_{\mathrm{s}})$,
${\mathbf{u}}^{\mathrm{G}(p_x)}_{r_{p_x}}(p_x)$, 
$\mathbf{u}^{\mathrm{G}(z)}_{r_z}(z)$ and 
${\mathbf{u}}^{\mathrm{G}(\omega)}_{r_{\omega}}(\omega)$
are the factor matrices defined on their respective one-dimensional domains, 
$x_\mathrm{s}\in\mathcal{I}_{x_\mathrm{s}}$, 
$z_\mathrm{s}\in\mathcal{I}_{z_\mathrm{s}}$,  
${p_x}\in\mathcal{I}_{p_x}$,
$z\in\mathcal{I}_z$, and
${\omega}\in\mathcal{I}_{\omega}$. 
The spatial representation $\mathbf{u}^{\mathrm{G}(x)}_{r_{x}}(x)$ is obtained by performing an inverse Fourier transform on ${\mathbf{u}}^{\mathrm{G}(p_x)}_{r_{p_x}}(p_x)$ per frequency as defined in equation~\eqref{eq:iff}. Each factor matrix corresponds to either a scalar or a vector function, depending on whether the SH or P–SV problem is analyzed. As the rank increases, the accuracy of the approximation improves. The core tensor specifies the interaction among the one-dimensional modes, and its multilinear contraction with the factor matrices determines how these modes combine to reconstruct the full Green’s tensor.

The virtual displacement field $\tilde{\mathbf{v}} (x_\mathrm{s}, z_\mathrm{s}, p_x, z, \omega)$ is separated similarly: 
\begin{equation}
	\tilde{\mathbf{v}}(x_\mathrm{s}, z_\mathrm{s}, p_x, z, \omega) \approx \boldsymbol{\mathcal{V}}=
\mathbf{v}^{(x_{\mathrm{s}})}(x_{\mathrm{s}}) \otimes 
\mathbf{v}^{(z_{\mathrm{s}})}(z_{\mathrm{s}}) \otimes
{\mathbf{v}}^{(p_x)}(p_x) \otimes
\mathbf{v}^{(z)}(z) \otimes
{\mathbf{v}}^{(\omega)}(\omega)
	\label{eq:vPGD}
\end{equation}

In order to solve equation~\eqref{eq:weak6} numerically, the domain is discretized using finite elements, where a one-dimensional discretization is performed along each dimension. The Green's displacement tensor $\tilde{\mathbf{u}}^{\mathrm{G}} (x_\mathrm{s}, z_\mathrm{s}, p_x, z, \omega)$ in equation~\eqref{eq:uPGD} and virtual displacement field $\tilde{\mathbf{v}}(x_\mathrm{s}, z_\mathrm{s}, p_x, z, \omega)$ in equation~\eqref{eq:vPGD} are approximated using shape functions 
$\mathbf{N}_{x_{\mathrm{s}}}$, 
$\mathbf{N}_{z_{\mathrm{s}}}$, 
$\mathbf{N}_{p_x}$,
$\mathbf{N}_{z}$ and
$\mathbf{N}_{\omega}$. For brevity, the explicit dependence of the shape functions on their respective coordinates is omitted in the notation below. This leads to the discretized representation $\boldsymbol{\mathcal{U}}$ of the Green's displacement tensor:


\begin{align}
\boldsymbol{\mathcal{U}} 
&\approx \sum_{r_{x_{\mathrm{s}}}=1}^{R_{x_{\mathrm{s}}}} 
\sum_{r_{z_{\mathrm{s}}}=1}^{R_{z_{\mathrm{s}}}} 
\sum_{r_{p_x}=1}^{R_{p_x}} 
\sum_{r_z=1}^{R_z} 
\sum_{r_{\omega}=1}^{R_{\omega}} 
g_{r_{x_{\mathrm{s}}} r_{z_{\mathrm{s}}} r_{p_x} r_z r_{\omega}} 
\mathbf{N}_{x_{\mathrm{s}}} \, \underline{\mathbf{u}}^{\mathrm{G}(x_{\mathrm{s}})}_{r_{x_{\mathrm{s}}}} \otimes 
\mathbf{N}_{z_{\mathrm{s}}} \, \underline{\mathbf{u}}^{\mathrm{G}(z_{\mathrm{s}})}_{r_{z_{\mathrm{s}}}} \otimes 
\mathbf{N}_{p_x} \, \underline{{\mathbf{u}}}^{\mathrm{G}(p_x)}_{r_{p_x}} \otimes
\mathbf{N}_{z} \, \underline{\mathbf{u}}^{\mathrm{G}(z)}_{r_z} \otimes 
\mathbf{N}_{\omega} \, \underline{{\mathbf{u}}}^{\mathrm{G}(\omega)}_{r_{\omega}} \nonumber \\
&= \left(
\mathbf{N}_{x_{\mathrm{s}}} \otimes
\mathbf{N}_{z_{\mathrm{s}}} \otimes
\mathbf{N}_{p_x} \otimes
\mathbf{N}_{z} \otimes
\mathbf{N}_{\omega}
\right)  \sum_{r_{x_{\mathrm{s}}}=1}^{R_{x_{\mathrm{s}}}} 
\sum_{r_{z_{\mathrm{s}}}=1}^{R_{z_{\mathrm{s}}}} 
\sum_{r_{p_x}=1}^{R_{p_x}} 
\sum_{r_z=1}^{R_z} 
\sum_{r_{\omega}=1}^{R_{\omega}} 
g_{r_{x_{\mathrm{s}}} r_{z_{\mathrm{s}}} r_{p_x} r_z r_{\omega}} \,
\underline{\mathbf{u}}^{\mathrm{G}(x_{\mathrm{s}})}_{r_{x_{\mathrm{s}}}}\otimes 
\underline{\mathbf{u}}^{\mathrm{G}(z_{\mathrm{s}})}_{r_{z_{\mathrm{s}}}} \otimes 
\underline{{\mathbf{u}}}^{\mathrm{G}(p_x)}_{r_{p_x}} \nonumber \\ &\quad \otimes 
\underline{\mathbf{u}}^{\mathrm{G}(z)}_{r_z} \otimes \underline{{\mathbf{u}}}^{\mathrm{G}(\omega)}_{r_{\omega}} 
\label{eq:uPGDdiscretized}
\end{align}

The virtual displacement field is restricted to a rank-one separable form. Its corresponding tensor representation, denoted by $\boldsymbol{\mathcal{V}}$, is expressed as follows:
\begin{align}
\boldsymbol{\mathcal{V}} 
&\approx 
\mathbf{N}_{x_{\mathrm{s}}} \, \underline{\mathbf{v}}^{(x_{\mathrm{s}})}
\otimes \mathbf{N}_{z_{\mathrm{s}}} \, \underline{\mathbf{v}}^{(z_{\mathrm{s}})}
\otimes \mathbf{N}_{p_x} \, \underline{{\mathbf{v}}}^{(p_x)} 
\otimes \mathbf{N}_{z} \, \underline{\mathbf{v}}^{(z)} 
\otimes \mathbf{N}_{\omega} \, \underline{{\mathbf{v}}}^{(\omega)} \nonumber \\
&= \left(
\mathbf{N}_{x_{\mathrm{s}}} \otimes
\mathbf{N}_{z_{\mathrm{s}}} \otimes
\mathbf{N}_{p_x} \otimes
\mathbf{N}_{z} \otimes
\mathbf{N}_{\omega}
\right) 
\underline{\mathbf{v}}^{(x_{\mathrm{s}})} 
\otimes \underline{\mathbf{v}}^{(z_{\mathrm{s}})} 
\otimes \underline{{\mathbf{v}}}^{(p_x)}
\otimes \underline{\mathbf{v}}^{(z)} 
\otimes \underline{{\mathbf{v}}}^{(\omega)} 
\label{eq:vPGD2}
\end{align} 
The test functions may be identical to or different from the trial functions, depending on whether a Galerkin ($\boldsymbol{\mathcal{U}}=\boldsymbol{\mathcal{V}}$) or a Petrov–Galerkin ($\boldsymbol{\mathcal{U}} \neq \boldsymbol{\mathcal{V}}$) \cite{nouy10a} formulation is adopted. Further details are provided in subsection~\ref{subsec:rankone}.

The body forces $\rho\tilde{\mathbf{b}}(x_\mathrm{s}, z_\mathrm{s}, p_x, z, \omega)$  arising from the concentrated pulse load and appearing on the right-hand side of equation~\eqref{eq:weak6}, are expressed using a low-rank separated approximation:
\begin{equation}
\rho \tilde{\mathbf{b}}(x_\mathrm{s}, z_\mathrm{s}, p_x, z, \omega) \approx \boldsymbol{\mathcal{F}} =
\mathbf{b}^{(x_{\mathrm{s}})}_{r_{x_{\mathrm{s}}}}(x_{\mathrm{s}}) \otimes 
\mathbf{b}^{(z_{\mathrm{s}})}_{r_{z_{\mathrm{s}}}}(z_{\mathrm{s}}) \otimes
{\mathbf{b}}^{(p_x)}_{r_{p_x}}(p_x) \otimes
\mathbf{b}^{(z)}_{r_z}(z) \otimes
{\mathbf{b}}^{(\omega)}_{r_{\omega}}(\omega)
\label{eq:bPGD}
\end{equation}
This separated representation $\boldsymbol{\mathcal{F}}$ is subsequently discretized using identical finite element shape functions as those used for the Green's displacement tensor and virtual displacement field:
{\small
\begin{align}
\boldsymbol{\mathcal{F}} &\approx
\sum_{r_b=1}^{R_b}  
\mathbf{N}_{x_{\mathrm{s}}} \, 
\underline{\mathbf{b}}^{(x_{\mathrm{s}})}_{r_{x_{\mathrm{s}}}} (x_{\mathrm{s}}) \otimes 
\mathbf{N}_{z_{\mathrm{s}}} \, 
\underline{\mathbf{b}}^{(z_{\mathrm{s}})}_{r_{z_{\mathrm{s}}}} (z_{\mathrm{s}}) \otimes 
\mathbf{N}_{p_x} \, 
\underline{{\mathbf{b}}}^{(p_x)}_{r_{p_x}}(p_x) \otimes 
\mathbf{N}_{z} \, 
\underline{\mathbf{b}}^{(z)}_{r_z}(z) \otimes 
\mathbf{N}_{\omega} \, 
\underline{{\mathbf{b}}}^{(\omega)}_{r_{\omega}}(\omega) \nonumber \\
&= \left(
\mathbf{N}_{x_{\mathrm{s}}} \otimes
\mathbf{N}_{z_{\mathrm{s}}} \otimes
\mathbf{N}_{p_x} \otimes
\mathbf{N}_{z} \otimes
\mathbf{N}_{\omega}
\right)
\sum_{r_b=1}^{R_b} 
\underline{\mathbf{b}}^{(x_{\mathrm{s}})}_{r_{x_{\mathrm{s}}}} (x_{\mathrm{s}}) \otimes 
\underline{\mathbf{b}}^{(z_{\mathrm{s}})}_{r_{z_{\mathrm{s}}}} (z_{\mathrm{s}}) \otimes 
\underline{{\mathbf{b}}}^{(p_x)}_{r_{p_x}}(p_x) \otimes 
\underline{\mathbf{b}}^{(z)}_{r_z}(z) \otimes 
\underline{{\mathbf{b}}}^{(\omega)}_{r_{\omega}}(\omega)
\label{eq:bPGD2}
\end{align}
}

\noindent 
Introducing the low rank approximations \eqref{eq:uPGDdiscretized}, \eqref{eq:vPGD2}, and \eqref{eq:bPGD2} in the weak form \eqref{eq:weak6}, results in:
\begin{equation}
	\boldsymbol{\mathcal{V}}^\mathrm{T}\boldsymbol{\mathcal{A}}\boldsymbol{\mathcal{U}}=\boldsymbol{\mathcal{V}}^\mathrm{T}\boldsymbol{\mathcal{F}}
	\label{eq:weakPGD}
\end{equation}
where $\boldsymbol{\mathcal{A}}$ is the dynamic stiffness tensor provided in \ref{sec:appendix2} and $\boldsymbol{\mathcal{F}}$ is the external load tensor expressed using equation~\eqref{eq:bPGD}.

\subsection{Computation of rank-one enrichments}
\label{subsec:rankone}
The rank-one updates in the algorithm are computed using an ALS procedure \cite{geor19a, kold09a}. At enrichment step ${R_u+1}$, all unknown terms associated with the new rank-one mode are grouped on the left-hand side, and the problem is restated as

\begin{equation}
\boldsymbol{\mathcal{V}}^\mathrm{T} \boldsymbol{\mathcal{A}} \boldsymbol{\mathcal{U}}_{R_u+1} = \boldsymbol{\mathcal{V}}^\mathrm{T} \boldsymbol{\mathcal{F}} - \boldsymbol{\mathcal{V}}^\mathrm{T} \boldsymbol{\mathcal{A}} \boldsymbol{\mathcal{U}}_{R_u}
\end{equation}
for all test function factors $\mathbf{\underline{v}}^{(x_{\mathrm{s}})} \in \mathbb{C}^{x_{\mathrm{s}}}, 
\mathbf{\underline{v}}^{(z_{\mathrm{s}})} \in \mathbb{C}^{z_{\mathrm{s}}}, 
\mathbf{\underline{v}}^{(p_x)} \in \mathbb{C}^{\mathbb{N}_{p_x}}, 
\mathbf{\underline{v}}^{(z)} \in \mathbb{C}^{\mathbb{N}_{z}},
\mathbf{\underline{v}}^{(\omega)} \in \mathbb{C}^{\mathbb{N}_\omega}$.

The tensor of the displacement field $\boldsymbol{\mathcal{U}}_{R_u}$ contains the computed nodal values from all previously converged modes:
\begin{equation}
\boldsymbol{\mathcal{U}}_{R_u} = \sum_{r_u=1}^{R_u} {\mathbf{\underline{u}}}_{r_u}^{(x_\mathrm{s})} \otimes {\mathbf{\underline{u}}}_{r_u}^{(z_\mathrm{s})} \otimes {\mathbf{\underline{u}}}_{r_u}^{(p_x)} \otimes {\mathbf{\underline{u}}}_{r_u}^{(z)} \otimes {\mathbf{\underline{u}}}_{r_u}^{(\omega)}
\end{equation} 
while the unknown mode $\boldsymbol{\mathcal{U}}_{R_u+1}$ is sought in the separable format:
\begin{equation}
\boldsymbol{\mathcal{U}}_{R_u+1} = {\mathbf{\underline{u}}}_{R_u+1}^{(x_\mathrm{s})} \otimes {\mathbf{\underline{u}}}_{R_u+1}^{(z_\mathrm{s})} \otimes {\mathbf{\underline{u}}}_{R_u+1}^{(p_x)} \otimes {\mathbf{\underline{u}}}_{R_u+1}^{(z)} \otimes {\mathbf{\underline{u}}}_{R_u+1}^{(\omega)}
\end{equation} 
Similarly, the test function $\boldsymbol{\mathcal{V}}$ is expressed as:
\begin{equation*}
\boldsymbol{\mathcal{V}} = {\mathbf{\underline{v}}}^{(x_\mathrm{s})} \otimes {\mathbf{\underline{v}}}^{(z_\mathrm{s})} \otimes {\mathbf{\underline{v}}}^{(p_x)} \otimes {\mathbf{\underline{v}}}^{(z)} \otimes {\mathbf{\underline{v}}}^{(\omega)}
\end{equation*}
For the present problem, the Petrov–Galerkin approach generally yields better convergence. In this strategy, the enrichment at iteration $R_u+1$ is computed by fixing all but one mode and solving for the remaining unknown. For instance, fix $\underline{\mathbf{u}}_{R_u+1}^{(z_{\mathrm{s}})}$,  
        $\underline{\mathbf{u}}_{R_u+1}^{(p_{x})}$, 
        $\underline{\mathbf{u}}_{R_u+1}^{(z)}$, 
        $\underline{\mathbf{u}}_{R_u+1}^{(\omega)}$, 
        $\underline{\mathbf{v}}^{(z_{\mathrm{s}})}$, 
        $\underline{\mathbf{v}}^{(p_{x})}$,  
        $\underline{\mathbf{v}}^{(z)}$ and
        $\underline{\mathbf{v}}^{(\omega)}$ to solve for 
        $\underline{\mathbf{u}}_{R_u+1}^{(x_{\mathrm{s}})}$:
\vspace{-1em}
\begin{multline}
        \left( I_{x_{\mathrm{s}}} \otimes 
        \underline{\mathbf{v}}^{(z_{\mathrm{s}})} \otimes 
        \underline{\mathbf{v}}^{(p_{x})} 
        \otimes \underline{\mathbf{v}}^{(z)} 
        \otimes \underline{\mathbf{v}}^{(\omega)} \right)^{\mathrm{H}}
        \boldsymbol{\mathcal{A}} 
        \left( \underline{\mathbf{u}}_{R_u+1}^{(x_{\mathrm{s}})} \otimes 
        \underline{\mathbf{u}}_{R_u+1}^{(z_{\mathrm{s}})} \otimes 
        \underline{\mathbf{u}}_{R_u+1}^{(p_{x})} \otimes 
        \underline{\mathbf{u}}_{R_u+1}^{(z)} \otimes 
        \underline{\mathbf{u}}_{R_u+1}^{(\omega)} \right) \\ =
		\left( I_{x_{\mathrm{s}}} \otimes 
		\underline{\mathbf{v}}^{(z_{\mathrm{s}})} \otimes 
		\underline{\mathbf{v}}^{(p_{x})} \otimes 
		\underline{\mathbf{v}}^{(z)} \otimes 
		\underline{\mathbf{v}}^{(\omega)} \right)^{\mathrm{H}} \boldsymbol{\mathcal{F}} - 
		\left( I_{x_{\mathrm{s}}} \otimes 
		\underline{\mathbf{v}}^{(z_{\mathrm{s}})} \otimes 
		\underline{\mathbf{v}}^{(p_{x})} \otimes 
		\underline{\mathbf{v}}^{(z)} \otimes 
		\underline{\mathbf{v}}^{(\omega)} \right)^{\mathrm{H}} \boldsymbol{\mathcal{A}} \boldsymbol{\mathcal{U}_{R_u}}
\end{multline}
The test factors are updated analogously within the Petrov–Galerkin alternating scheme to satisfy the orthogonality condition.
The same alternating procedure is applied cyclically to
$\underline{\mathbf{u}}_{R_u+1}^{(z_{\mathrm{s}})}$,  
$\underline{\mathbf{u}}_{R_u+1}^{(p_{x})}$, 
$\underline{\mathbf{u}}_{R_u+1}^{(z)}$,  
$\underline{\mathbf{u}}_{R_u+1}^{(\omega)}$
and the corresponding test functions 
$\underline{\mathbf{v}}^{(z_{\mathrm{s}})}$,  
$\underline{\mathbf{v}}^{(p_{x})}$,  
$\underline{\mathbf{v}}^{(z)}$ and
$\underline{\mathbf{v}}^{(\omega)}$
until convergence of the enrichment is achieved. The alternating iterations are terminated once the change in the rank-one factors between successive iterations falls below a prescribed convergence tolerance. After orthonormalization, a newly computed mode is retained only if its norm exceeds a prescribed truncation tolerance, thereby preventing the inclusion of nearly linearly dependent directions. The enrichment proceeds until the prescribed maximum number of modes is reached.


\begin{algorithm}[!b]
\footnotesize  
\caption{Greedy Tucker Approximation for Low-Rank Green’s Functions (GTA-GF)\\
{\footnotesize \emph{(adopted from Georgieva and Hofreither \cite{geor19a})}}}
\label{alg:GTA_GF}
\begin{algorithmic}
\State \textbf{function} GTA-GF($\mathcal{A}: \mathbb{R}^I \to \mathbb{R}^I$, $\mathcal{F} \in \mathbb{R}^I$, Tol, nIterations)
\State \textbf{let} $U_j = \{\}$ \text{ for } $j = 1, \dots, d$
\State \text{Initialize the core tensor} $\boldsymbol{\mathcal{G}} = 0$
\For{$l = 1, 2, \dots, R$}
    \State $T \gets (U_1, \dots, U_d) \cdot \boldsymbol{\mathcal{G}}$ 
    \State \text{ALS for rank 1 approximation of the factor matrices:}
    \State Choose nonzero starting vectors $x_j \in \mathbb{R}^{N_j}$, $j = 1, \dots, d$
    \State Initialize error = Inf, itr = 1
    \While{error $>$ Tol \text{ and } itr $<$ nIterations}
        \For{$k = 1, 2, \dots, d$}
            \State Solve the linear least squares problem
            \State $x_k \gets \underset{y \in \mathbb{R}^{N_k}}{\arg \min} \| \boldsymbol{\mathcal{F}} - \boldsymbol{\mathcal{A}}  (x_1 \otimes \dots \otimes x_{k-1} \otimes y \otimes x_{k+1} \otimes \dots \otimes x_d) \|^2$
        \EndFor
        \State \text{Error computation:}
        \State $\text{error} = \max\limits_{1 \leq i \leq d} \| {x_i}^\mathrm{(itr)} - {x_i}^\mathrm{(itr-1)} \|$ 
        \State itr = itr + 1 
    \EndWhile
        \For{$j = 1, \dots, d$}
            \State $U_j \gets \text{orthogonalize}(U_j, x_j)$ 
        \EndFor
        \State $\boldsymbol{\mathcal{G}} \gets (\boldsymbol{\mathcal{V}}^T \boldsymbol{\mathcal{A}} \boldsymbol{\mathcal{U}})^{-1} \boldsymbol{\mathcal{V}}^T \boldsymbol{\mathcal{F}}$ 
\EndFor
\State \textbf{return} $(U_1, \dots, U_d), \boldsymbol{\mathcal{G}}$ 
\State \textbf{end function}
\end{algorithmic}
\end{algorithm}

\subsection{Greedy Tucker Approximation}
\label{subsec:GTA_alg}
Following each rank-one update, the resulting vectors are appended to the corresponding factor matrices of the Tucker decomposition, thereby progressively enriching the low-rank representation of the Green’s tensor. The best rank-one correction to the current residual, obtained through the alternating least-squares strategy described in the previous subsection, is incorporated into these mode matrices and orthonormalized using a modified Gram–Schmidt procedure \cite{golu13a} to ensure linear independence and numerical stability. 

Once the factor matrices are updated, the weak form~\eqref{eq:weakPGD} is projected onto the reduced trial and test spaces spanned by the current modes. This projection yields a compact reduced system that defines the core tensor, allowing matrix–vector operations to be carried out efficiently prior to inversion. The reduced system, evaluated over all frequencies and wavenumbers, is substantially smaller than the original full order formulation. Accordingly, the core tensor provides a compact representation of the reduced system, with dimensions governed by the number of modes retained in each direction of the decomposition:
\begin{equation}
\boldsymbol{\mathcal{G}} \;=\; (\boldsymbol{\mathcal{V}}^{\mathrm{T}} \boldsymbol{\mathcal{A}} \boldsymbol{\mathcal{U}})^{-1}\,\boldsymbol{\mathcal{V}}^{\mathrm{T}}\boldsymbol{\mathcal{F}}
\end{equation}
This projection step preserves the contributions of all previously computed enrichments within a consistent low-rank structure. By alternating enrichment, orthogonalization, and core updates in this manner, the GTA builds a progressively refined representation of the Green’s function tensor. The present formulation closely follows that of \cite{geor19a}. The main differences are the use of a Petrov–Galerkin projection instead of a Galerkin formulation, and the application of a modified Gram–Schmidt orthonormalization procedure. The complete procedure is summarized in Algorithm~\ref{alg:GTA_GF} \cite{geor19a}. Tensor computations are performed using the Tensor Toolbox for MATLAB \cite{bade25a}.

\section{Numerical examples}
\label{sec:examples}
\subsection{Layered soil on a bedrock} \label{subsec:ex1}
As a first example, we consider a soil stratum overlying rigid bedrock. The soil layer has a thickness $h=20\,\mathrm{m}$, shear wave velocity $C_\mathrm{s}=100\,\mathrm{m/s}$, dilatational wave velocity $C_\mathrm{p}=200\,\mathrm{m/s}$ and density $\rho=1800\,\textrm{kg/m}^3$.  A damping ratio $\beta_\mathrm{s}=\beta_\mathrm{p}=0.02$ is assumed for both shear and volumetric deformation (table \ref{tab:soilpropterties}). Material damping is assumed to be hysteretic and rate-independent. According to the correspondence principle, the complex Lame constants are defined as $\mu^{\star} = \mu(1 + 2\beta_\mathrm{s}\mathrm{i})$ and $(\lambda^{\star} + 2\mu^{\star}) = (\lambda + 2\mu)(1 + 2\beta_\mathrm{p}\mathrm{i})$ \cite{mole80a}.
The model coordinates comprise receivers in the $z$-direction, slowness $p_x$, and frequency $\omega$. The analysis is carried out for frequencies ranging from $0$ to $15 \,\mathrm{Hz}$ for both SH and P-SV wave propagation. The Green’s displacements in the wavenumber–frequency domain are approximated using the GTA algorithm (Algorithm~\ref{alg:GTA_GF}) and expressed in Tucker format as:

\begin{equation} \label{eq:gta_ex1a}
\tilde{\mathbf{u}}^{\mathrm{G}} (p_x, z, \omega) \approx
\sum_{r_{p_x}=1}^{R_{p_x}}
\sum_{r_z=1}^{R_z} 
\sum_{r_\omega=1}^{R_\omega}
\displaystyle g_{r_{p_x} r_z r_{\omega}} 
{\mathbf{u}}^{\mathrm{G}(p_x)}_{r_{p_x}}(p_x) \otimes
\mathbf{u}^{\mathrm{G}(z)}_{r_z}(z) \otimes
{\mathbf{u}}^{\mathrm{G}(\omega)}_{r_{\omega}}(\omega)
\end{equation}
\par
The approximation is performed using structured sampling over slowness $p_x$, receiver locations in the $z$-direction, and frequency $\omega$. To improve efficiency, the wavenumber domain is sampled in terms of slowness, defined as $p_x={k_x}/{\omega}$. Multiplying the slowness $p_x$ with the shear wave velocity $C_\mathrm{s}$ yields the dimensionless wavenumber $\bar{k}_x = k_x C_\mathrm{s} / \omega = p_x C_\mathrm{s}$. The frequency range is sampled uniformly from $0$ to $15~\mathrm{Hz}$ with a bin $\Delta f = 0.025~\mathrm{Hz}$, while the wavenumber domain is sampled logarithmically using $1500$ points over the range $\bar{k}_x \in [10^{-2},\,10^{1}]$. The receiver depth is sampled uniformly over the soil layer as $z \in [0,\,20]~\mathrm{m}$ with a spacing of $0.2~\mathrm{m}$. The source location is fixed at $x_{\mathrm{s}} = 0$ and $z_{\mathrm{s}} = 0$, while the horizontal receiver coordinate is sampled over the interval $x \in [0,\,50]~\mathrm{m}$ using $250$ uniformly distributed points.
The GTA algorithm constructs a reduced subspace greedily, expanding iteratively along each tensor dimension. The approximation process is controlled by a stopping criterion based on a prescribed tolerance of $10^{-5}$ and a maximum number of 40 modes, which collectively influence the convergence behavior and accuracy of the solution.

\begin{table}[!htb]
	\centering
        \caption{Soil properties for numerical examples.}
        \label{tab:soilpropterties}
	\begin{tabular}{l c c c c c c c c}
		\toprule
		Soil profile & $h$ & $C_\mathrm{s}$ & $C_\mathrm{p}$ & $\beta_\mathrm{s}$ & $\beta_\mathrm{p}$ & $\rho$\\
		& [m] & [m/s] & [m/s] & [-] & [-] & [kg/m$^3$] \\
		\midrule
		Layer on bedrock & 20.0 & 100 & 200 & 0.020 & 0.020 & 1800\\
		\midrule
		& 3.7 & 50 & 1761 & 0.025 & 0.025 & 1107 \\
		& 7.0 & 75 & 1719 & 0.025 & 0.025 & 1500 \\
		Groene Hart & 8.3 & 180 & 1686 & 0.025 & 0.025 & 1970\\
		& 9.3 & 240 & 1715 & 0.025 & 0.025 & 1970 \\
		& $\infty$ & 260 & 1726 & 0.025 & 0.025 & 1970 \\
		\bottomrule
	\end{tabular}
\end{table}
 
\begin{figure}[!htb]
   \begin{center}
	 (a) \includegraphics[width=0.27\textwidth]{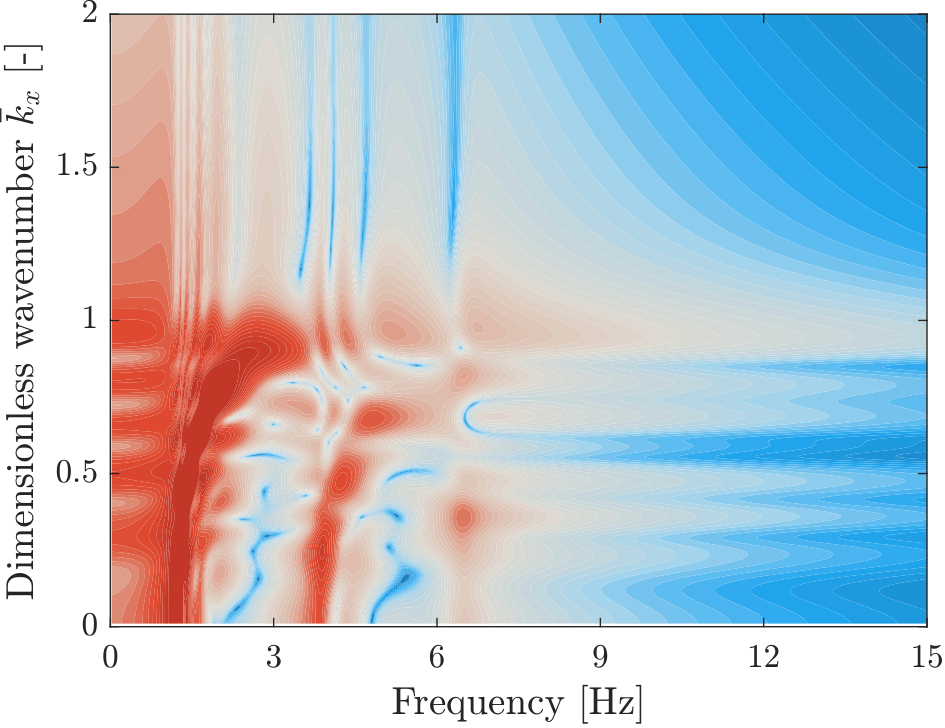}
	 (b) \includegraphics[width=0.27\textwidth]{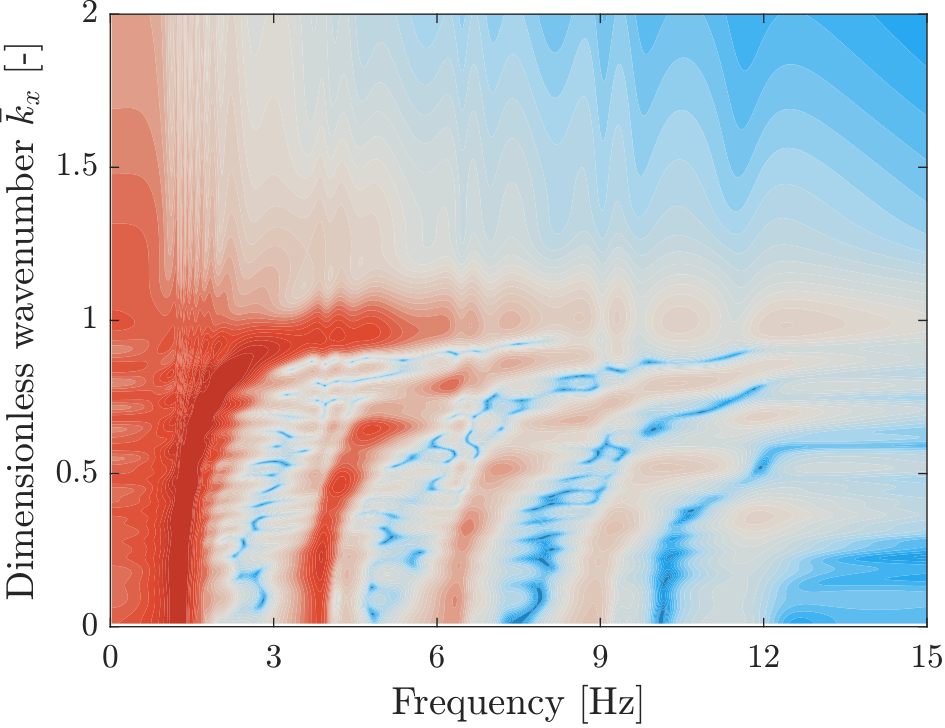} \\
	 (c) \includegraphics[width=0.27\textwidth]{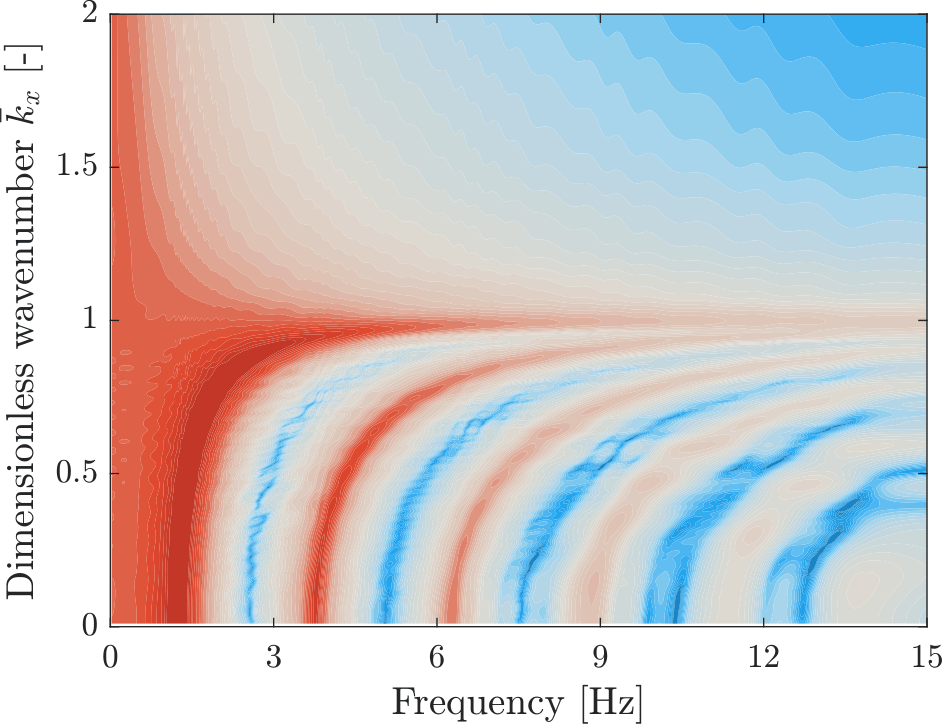} 
	 (d) \includegraphics[width=0.27\textwidth]{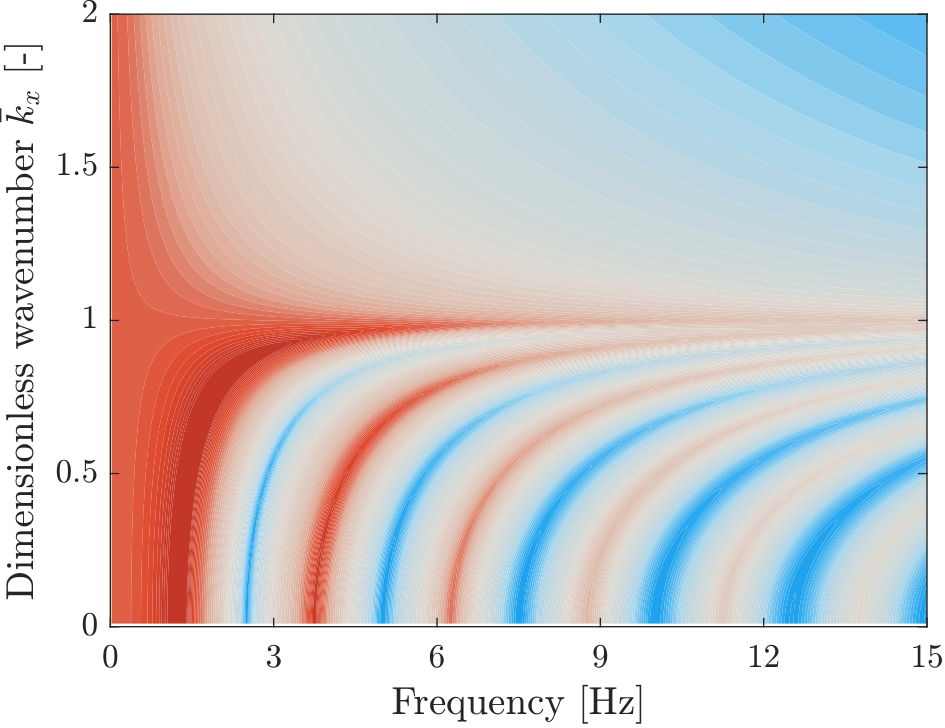} \\
	 \vspace{0.2cm} \hspace{0.6cm}
	     \includegraphics[width=0.27\textwidth]{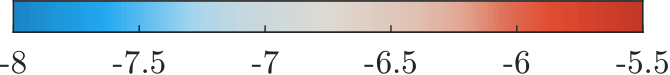} 
   \end{center}
	 \caption{SH Green's displacement $\tilde{u}^{\mathrm{G}}_{yy} (p_x, z, \omega)$ [$\mathrm{m}^2/\mathrm{Hz}$] in the wavenumber-frequency domain evaluated at $z = 0 \, \mathrm{m}$ for an out-of-plane concentrated load applied at $x_\mathrm{s}=0 \, \mathrm{m}$ and $z_\mathrm{s}=0 \, \mathrm{m}$ using ROM with (a) 10 modes, (b) 20 modes, (c) 40 modes and (d) FOM.}
    \label{fig:ex1a_sh_f_kxd}
\end{figure}

Figure~\ref{fig:ex1a_sh_f_kxd} illustrates the SH Green’s displacement in the wavenumber-frequency domain computed at $z = 0 \, \mathrm{m}$.
Figure~\ref{fig:ex1a_sh_f_kxd}d presents the solution computed with the direct stiffness method, referred to as the full order model (FOM). 
The results describe the response of a soil stratum overlying rigid bedrock, where wave propagation is governed primarily by Love waves. The fundamental Love mode appears at $1.25 \, \mathrm{Hz}$ and corresponds to the first eigenfrequency of the clamped soil layer, estimated as $C_\mathrm{s}/(4h)$. At this frequency, the phase velocity tends to infinity. The next mode appears at $3C_\mathrm{s}/(4h)$, corresponding to the second eigenfrequency. Higher-order modes emerge at higher frequencies when the layer thickness equals one quarter of the shear wavelength, reflecting the dispersive behavior of the soil layer constrained by the rigid bedrock.
At higher frequencies, where the shear wavelength becomes small relative to the layer thickness, the response tends toward that of a homogeneous halfspace governed by the same shear wave velocity.
Figure~\ref{fig:ex1a_sh_f_kxd} shows the ROM approximations based on the thin layer formulation for 10, 20, and 40 modes, respectively.   For a low number of modes, corresponding to 10 modes (figure~\ref{fig:ex1a_sh_f_kxd}a), the ROM captures the overall trend of the fundamental Love mode but remains less accurate in representing the higher-order branches (figure~\ref{fig:ex1a_sh_f_kxd}a). Increasing the number of modes to 20 (figure~\ref{fig:ex1a_sh_f_kxd}b)  and 40 (figure~\ref{fig:ex1a_sh_f_kxd}c) progressively enhances the accuracy, enabling the ROM to reproduce multiple higher-order dispersion branches and to closely match the FOM across the entire frequency range.

\begin{figure}[!htb]
   \begin{center}
	 (a) \includegraphics[width=0.27\textwidth]{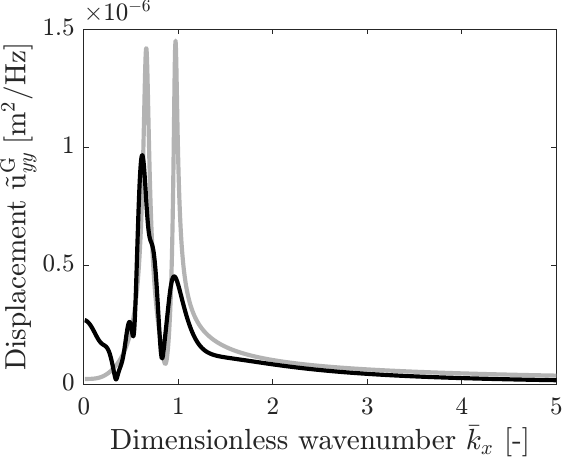}
	 (b) \includegraphics[width=0.27\textwidth]{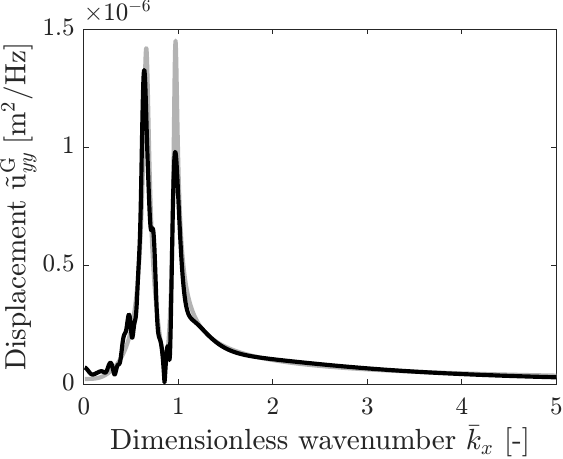}
	 (c) \includegraphics[width=0.27\textwidth]{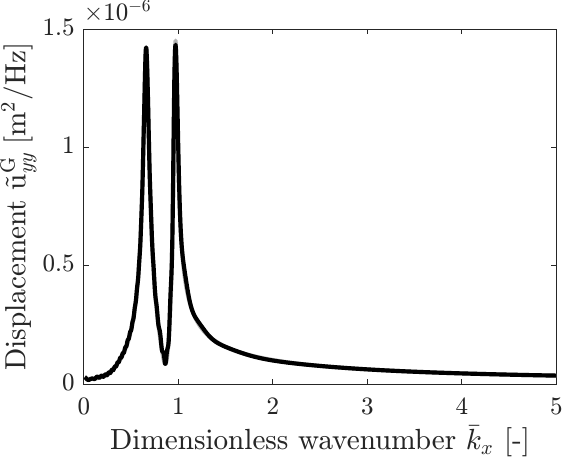} 
	\end{center}
	 \caption{SH Green's displacement $\tilde{u}^{\mathrm{G}}_{yy} (p_x, z, \omega)$ as a function of dimensionless wavenumber evaluated at $z = 0 \, \mathrm{m}$ and $f=5 \, \mathrm{Hz}$ for an out-of-plane concentrated load applied at $x_\mathrm{s}=0 \, \mathrm{m}$ and $z_\mathrm{s}=0 \, \mathrm{m}$ using ROM (black) with (a) 10 modes, (b) 20 modes and (c) 40 modes, compared with FOM (gray).}\label{fig:ex1a_sh_kxd}
\end{figure}

\begin{figure}[!htb]
   \begin{center}
	 (a) \includegraphics[width=0.27\textwidth]{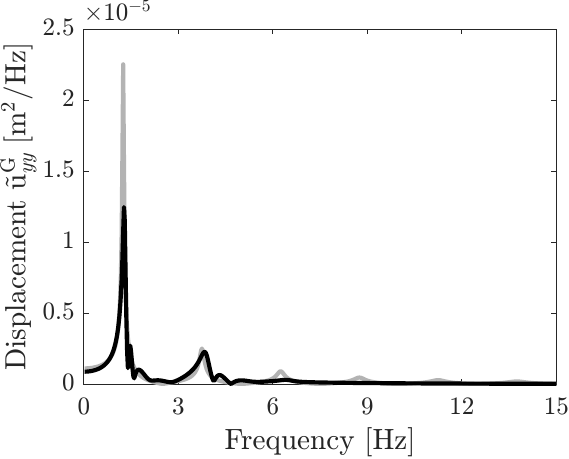}
	 (b) \includegraphics[width=0.27\textwidth]{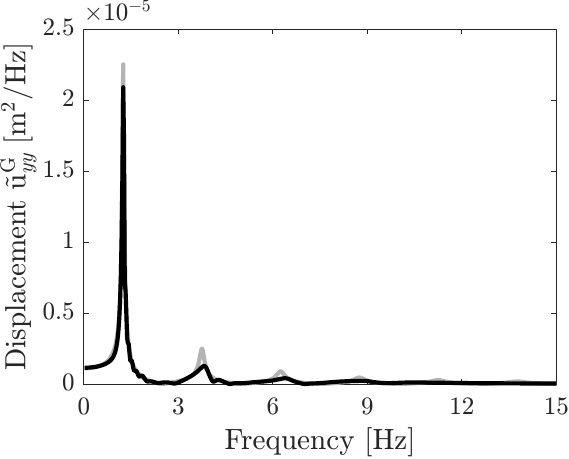}
	 (c) \includegraphics[width=0.27\textwidth]{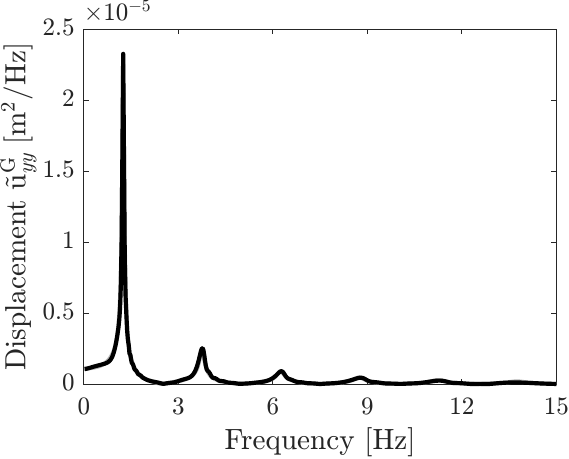} 
   \end{center}
	 \caption{SH Green's displacement $\tilde{u}^{\mathrm{G}}_{yy} (p_x, z, \omega)$ as a function of frequency evaluated at  $\bar{k}_{x} = 0 \, \mathrm{rad/m}$ and $z = 0 \, \mathrm{m}$ for an out-of-plane concentrated load applied at $x_\mathrm{s}=0 \, \mathrm{m}$ and $z_\mathrm{s}=0 \, \mathrm{m}$ using ROM (black) with (a) 10 modes, (b) 20 modes and (c) 40 modes, compared with FOM (gray).}
    \label{fig:ex1a_sh_ftilde}
\end{figure}

To provide a clearer comparison of the convergence between the ROM and FOM, slices of the SH Green's displacements for a single frequency and dimensionless wavenumber $\bar{k}_x$ are examined. In figure~\ref{fig:ex1a_sh_kxd}, the distinct peaks corresponding to the fundamental and higher-order Love modes in the soil layer are reproduced by the ROM with 40 modes in close agreement with the FOM. These peaks represent the wavenumbers at which constructive interference occurs between multiple reflections within the layer, resulting in guided wave propagation. In Figure~\ref{fig:ex1a_sh_ftilde}, the pronounced peaks associated with the eigenfrequencies of the clamped soil layer, estimated as $C_\mathrm{s}/4h$ and $3C_\mathrm{s} / 4h$ are also accurately approximated by the ROM.

\begin{figure}[!htb]
   \begin{center}
	 (a) \includegraphics[width=0.25\textwidth]{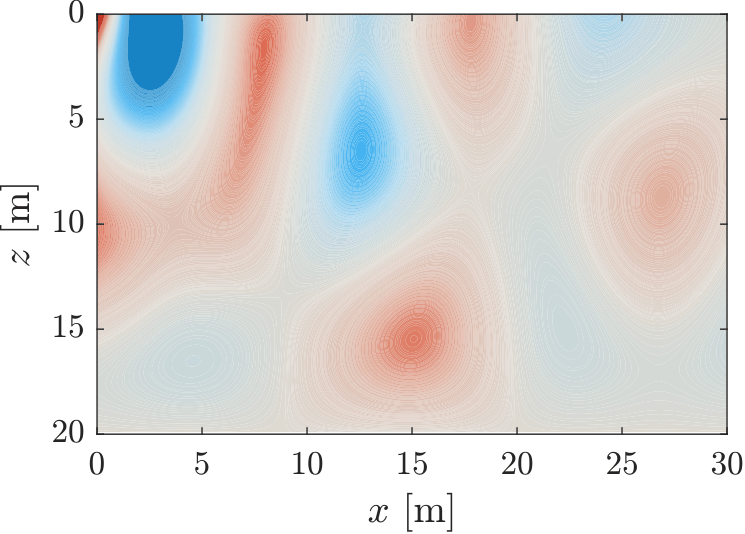}
	 (b) \includegraphics[width=0.25\textwidth]{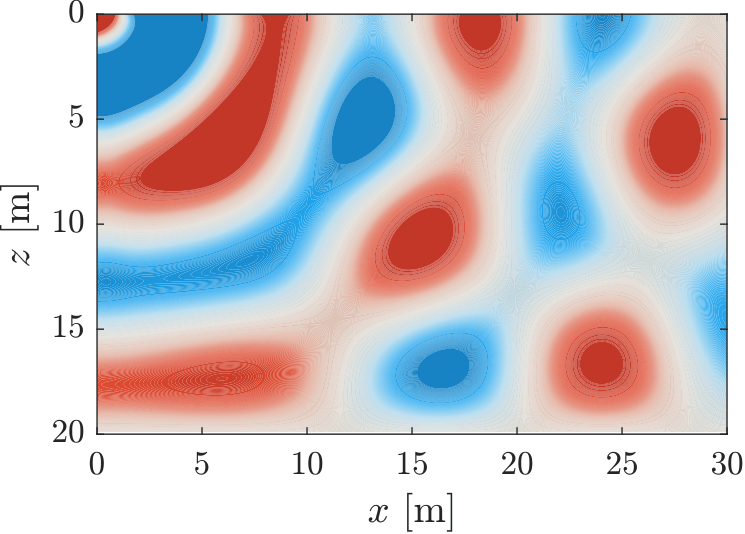} \\
	 (c) \includegraphics[width=0.25\textwidth]{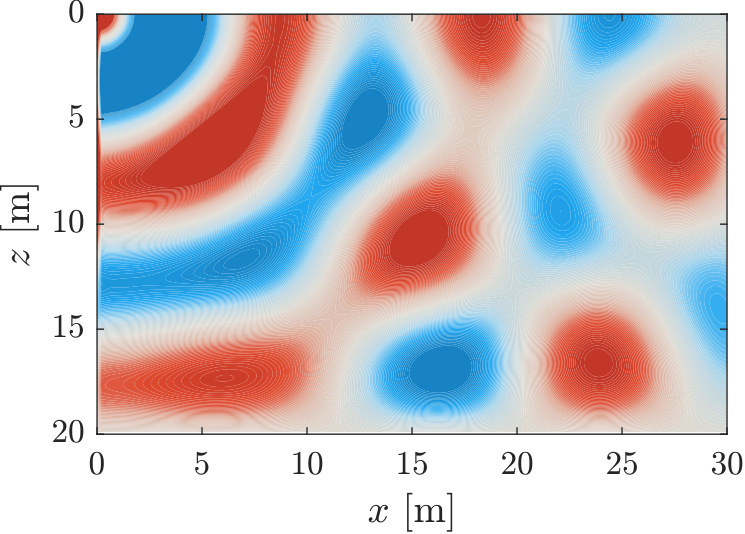} 
	 (d) \includegraphics[width=0.25\textwidth]{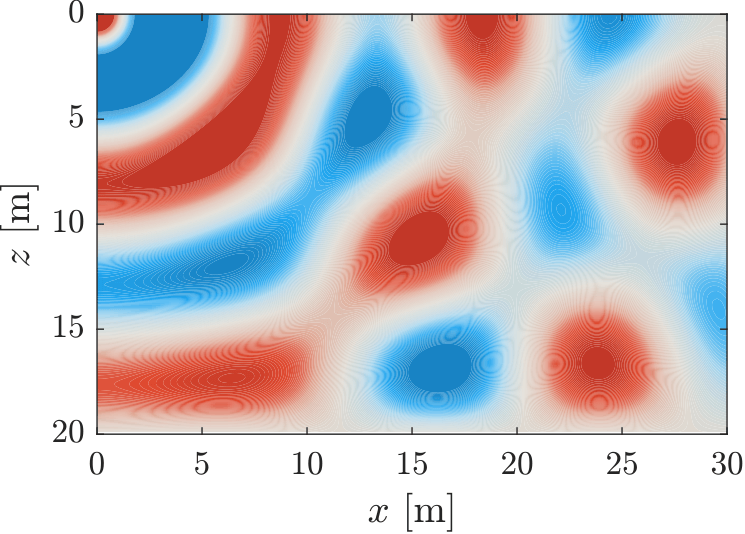} \\
	 \vspace{0.2cm} \hspace{0.6cm}
	     \includegraphics[width=0.27\textwidth]{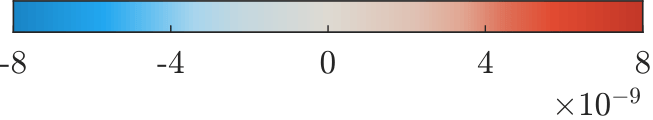} 
   \end{center}
	 \caption{SH Green's displacement $\hat{u}^{\mathrm{G}}_{yy} (x, z, \omega)$ [$\mathrm{m}/\mathrm{Hz}$] in the spatial-frequency domain at $f=10 \, \mathrm{Hz}$ for an out-of-plane concentrated load applied at $x_\mathrm{s}=0 \, \mathrm{m}$ and $z_\mathrm{s}=0 \, \mathrm{m}$ using ROM with (a) 10 modes, (b) 20 modes, (c) 40 modes and (d) FOM.}\label{fig:ex1_sh_contour_x}
\end{figure}

Green’s displacements in the spatial–frequency domain are obtained by applying an inverse logarithmic Fourier transform \cite{talm78a, sche09a} to equation (\ref{eq:gta_ex1a}). The low-rank structure of the decomposition allows this transformation to be carried out directly on the relevant factor matrix, i.e. ${\mathbf{u}}^{\mathrm{G}(p_x)}_{r_{p_x}}(p_x)$, rather than to the full Green’s tensor. Unlike the conventional approach, which requires inverse transforms over the entire discretized grid of size $n(p_x) \times n(z) \times n(\omega)$, where $n(p_x)$ is the number of slowness samples, $n(z)$ the number of vertical discretization points, and $n(\omega)$ the number of frequency samples, the proposed method requires the transformation only over the factor matrix ${\mathbf{u}}^{\mathrm{G}(p_x)}_{r_{p_x}}(p_x)$. The reduced order model therefore involves only $R_{p_x}$ inverse Fourier transforms per frequency, which are considerably fewer than those required in the full order model, resulting in faster computation. 
Figure~\ref{fig:ex1_sh_contour_x} illustrates the SH Green’s displacement $\hat{u}^{\mathrm{G}}_{yy} (x, z, \omega)$ evaluated at $f = 10 \, \mathrm{Hz}$ for an out-of-plane concentrated load applied at $x_\mathrm{s}=0 \, \mathrm{m}$ and $z_\mathrm{s}=0 \, \mathrm{m}$. The contours exhibit the depth-dependent pattern of horizontally polarized shear motion within the soil layer, with amplitudes concentrated near the surface and gradually decaying away from the source position. The horizontal oscillation pattern corresponds to the shear wavelength of approximately $10 \, \mathrm{m}$ at this frequency, while the gradual amplitude decay along the surface is primarily governed by material damping.
The reduced order solution with 40 modes (figure~\ref{fig:ex1_sh_contour_x}c) accurately reproduces the spatial characteristics of the reference full order solution shown in figure~\ref{fig:ex1_sh_contour_x}d.

\begin{figure}[!htb]
   \begin{center}
	 (a) \includegraphics[width=0.27\textwidth]{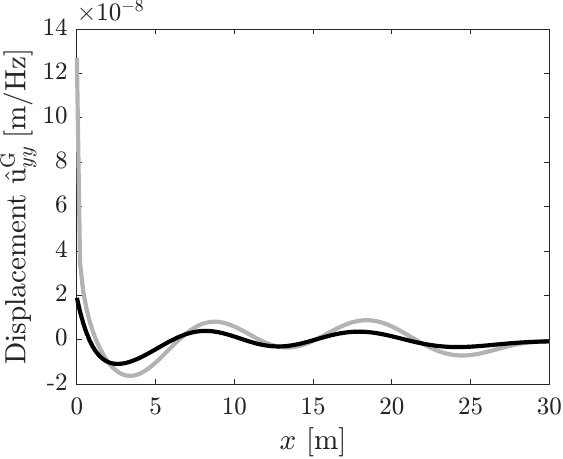}
	 (b) \includegraphics[width=0.27\textwidth]{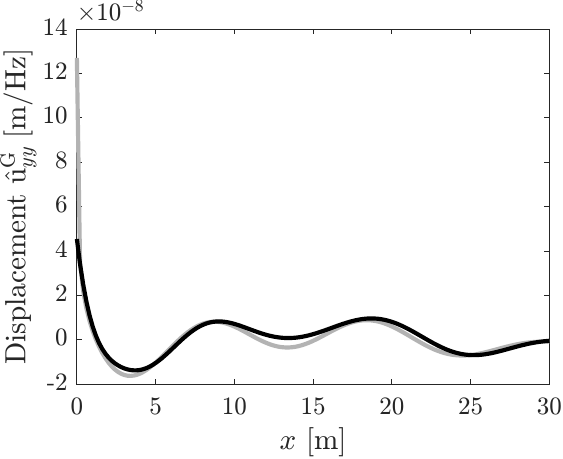}
	 (c) \includegraphics[width=0.27\textwidth]{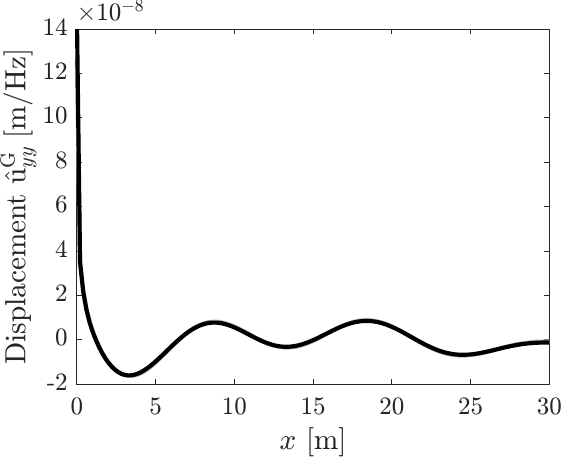} 
	\end{center}
	 \caption{SH Green's displacement $\hat{u}^{\mathrm{G}}_{yy} (x, z, \omega)$ as a function of $x$ evaluated at $z = 0 \, \mathrm{m}$ and $f=10 \, \mathrm{Hz}$ for an out-of-plane concentrated load applied at $x_\mathrm{s}=0 \, \mathrm{m}$ and $z_\mathrm{s}=0 \, \mathrm{m}$ using ROM (black) with (a) 10 modes, (b) 20 modes and (c) 40 modes, compared with FOM (gray).}
	 \label{fig:ex1_sh_x}
\end{figure}

\begin{figure}[!htb]
   \begin{center}
	 (a) \includegraphics[width=0.27\textwidth]{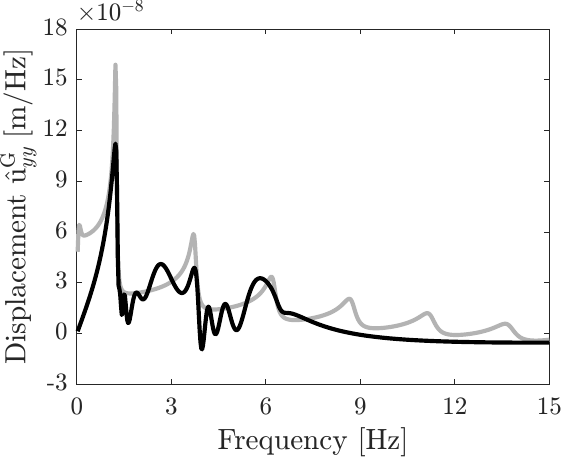}
	 (b) \includegraphics[width=0.27\textwidth]{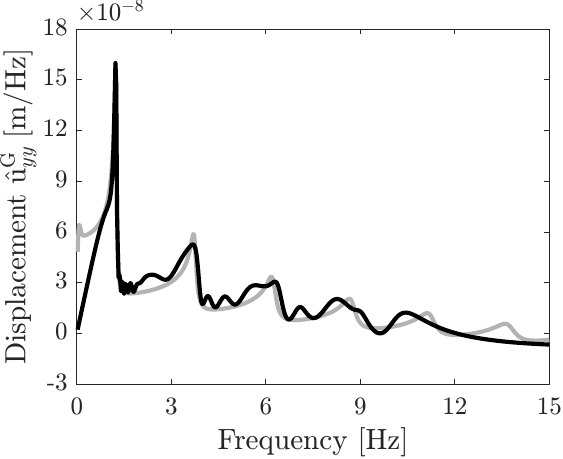}
	 (c) \includegraphics[width=0.27\textwidth]{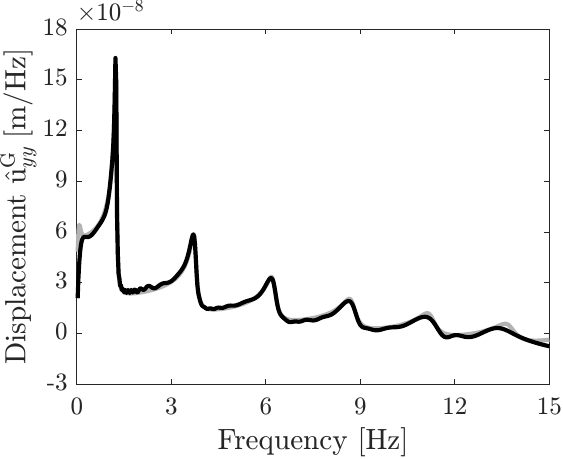} 
   \end{center}
	 \caption{SH Green's displacement $\hat{u}^{\mathrm{G}}_{yy} (x, z, \omega)$ as a function of frequency evaluated at $x = 1 \, \mathrm{m}$ and $z = 0 \, \mathrm{m}$ for an out-of-plane concentrated load applied at $x_\mathrm{s}=0 \, \mathrm{m}$ and $z_\mathrm{s}=0 \, \mathrm{m}$ using ROM (black) with (a) 10 modes, (b) 20 modes and (c) 40 modes, compared with FOM (gray).}\label{fig:ex1_sh_fx}
\end{figure}

The convergence of the ROM is assessed further. Figure~\ref{fig:ex1_sh_x} presents the SH displacement $\hat{u}^{\mathrm{G}}_{yy} (x, z, \omega)$ as a function of $x$ evaluated at $z=0 \, \mathrm{m}$ and $f = 10 \, \mathrm{Hz}$.
The ROM accurately captures the oscillatory behavior and the gradual decay in amplitude along the surface, capturing the shear wavelength ($10 \, \mathrm{m}$) at this frequency as well as the gradual attenuation of horizontally polarized shear motion within the layer.
Figure~\ref{fig:ex1_sh_fx} demonstrates the same displacement component as a function of frequency evaluated for a receiver located at $x=1 \, \mathrm{m}$ and $z=0 \, \mathrm{m}$.
The frequency-dependent behavior exhibits distinct amplitude peaks corresponding to the resonant frequencies of the soil layer. These peaks are increasingly well captured as the number of ROM modes increases. The ROM with 40 modes reproduces the FOM within the prescribed tolerance.

\begin{figure}[!htb]
   \begin{center}
	 (a) \includegraphics[width=0.27\textwidth]{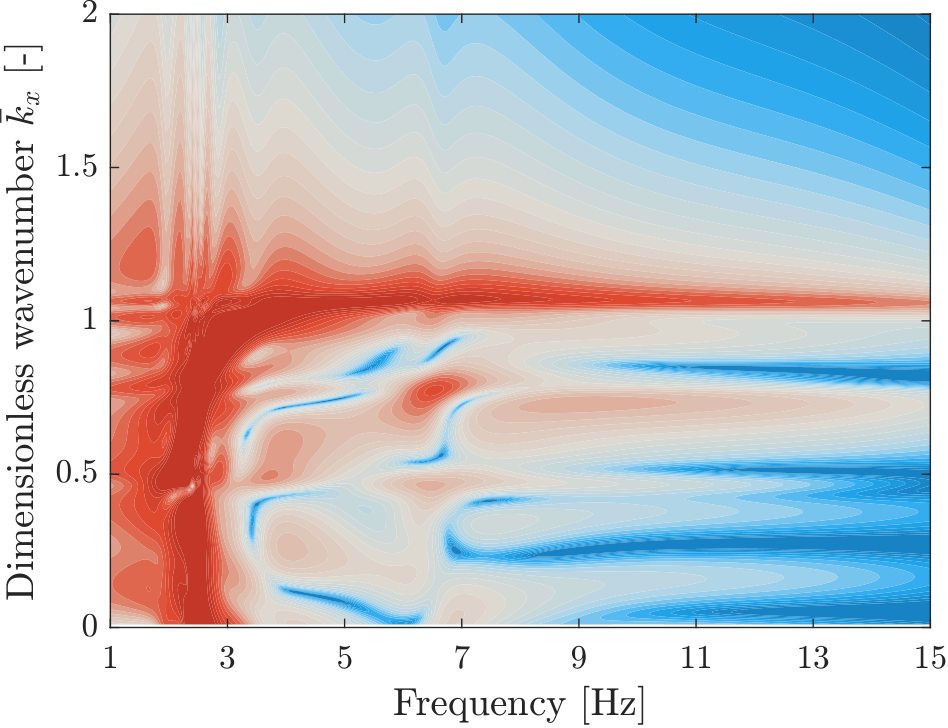}
	 (b) \includegraphics[width=0.27\textwidth]{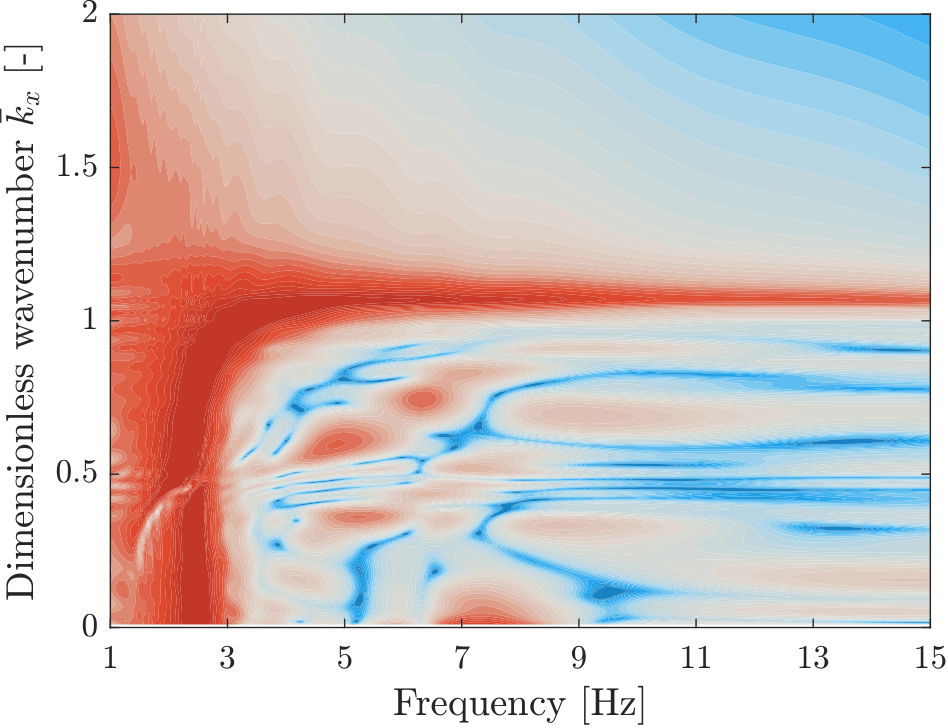} \\
	 (c) \includegraphics[width=0.27\textwidth]{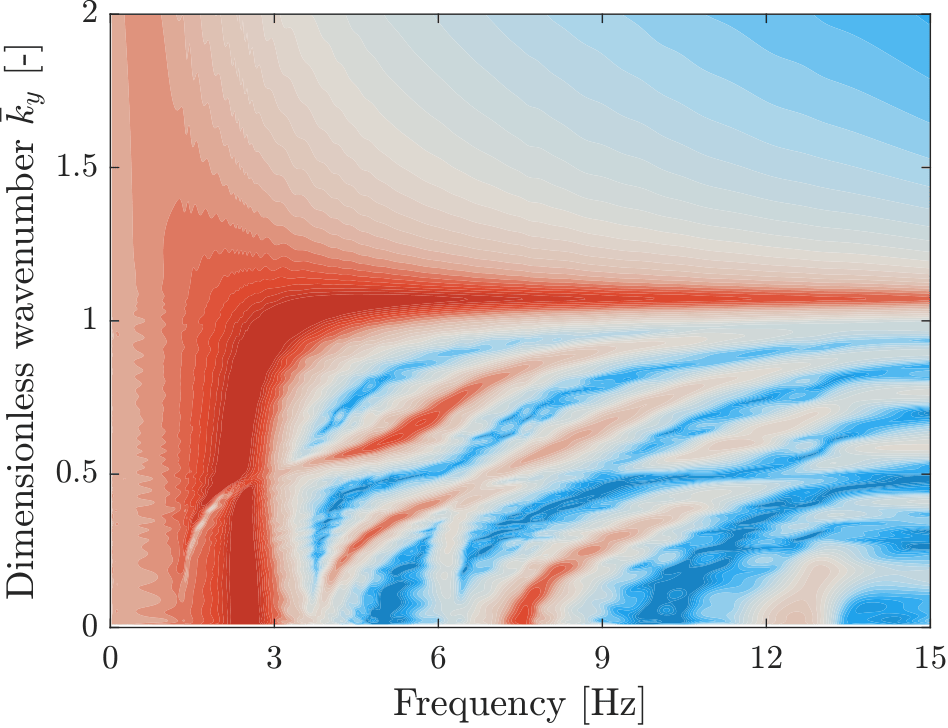} 
	 (d) \includegraphics[width=0.27\textwidth]{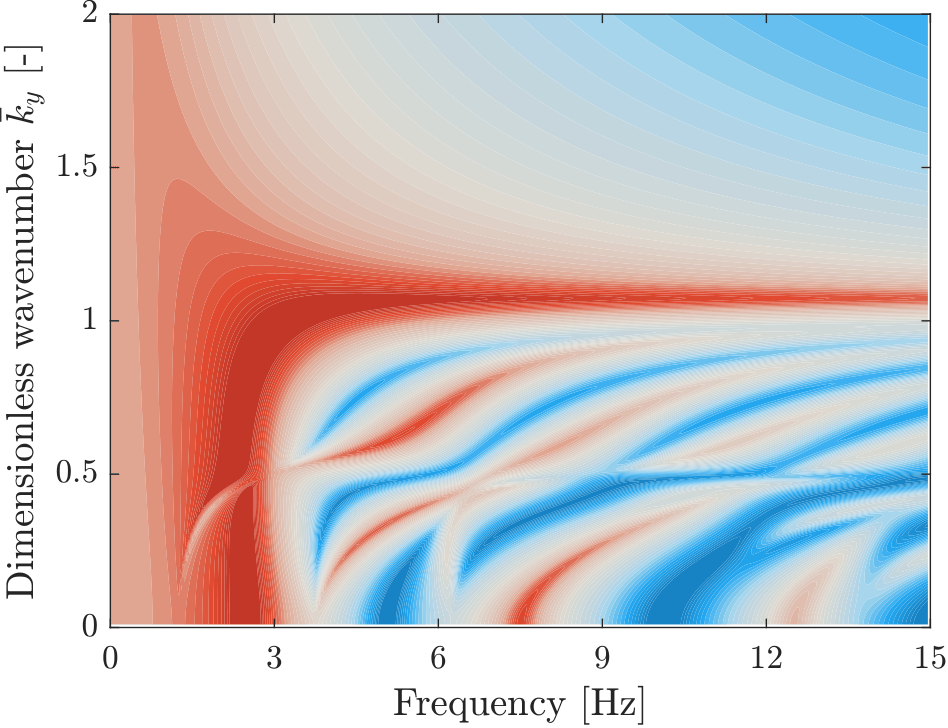}  \\
	     \vspace{0.2cm} \hspace{0.6cm}
	     \includegraphics[width=0.27\textwidth]{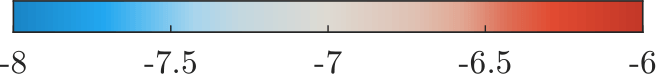} 
   \end{center}
	 \caption{P-SV Green's displacement $\tilde{u}^{\mathrm{G}}_{zz} (p_x, z, \omega)$ [$\mathrm{m}^2/\mathrm{Hz}$] in the wavenumber-frequency domain evaluated at $z=0 \, \mathrm{m}$ for a vertical concentrated load applied at $x_\mathrm{s}=0 \, \mathrm{m}$ and $z_\mathrm{s}=0 \, \mathrm{m}$ using ROM with (a) 10 modes, (b) 20 modes, (c) 40 modes and (d) FOM.}
    \label{fig:ex1a_psv_f_kxd}
\end{figure}

\begin{figure}[!htb]
   \begin{center}
	 (a) \includegraphics[width=0.27\textwidth]{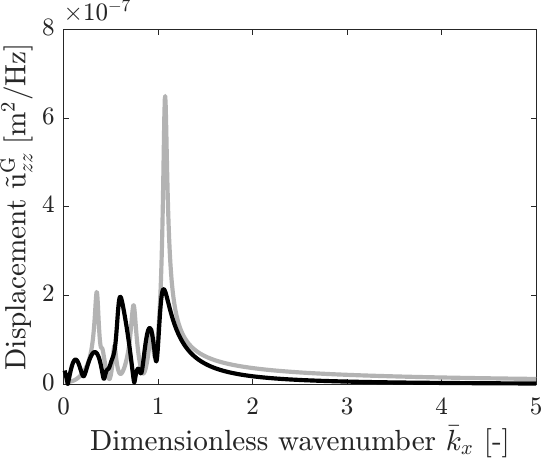}
	 (b) \includegraphics[width=0.27\textwidth]{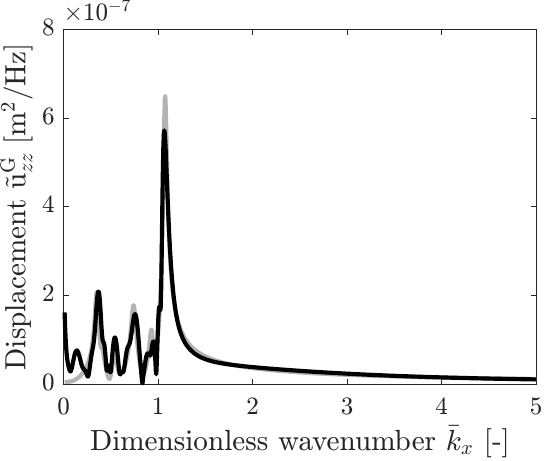}
	 (c) \includegraphics[width=0.27\textwidth]{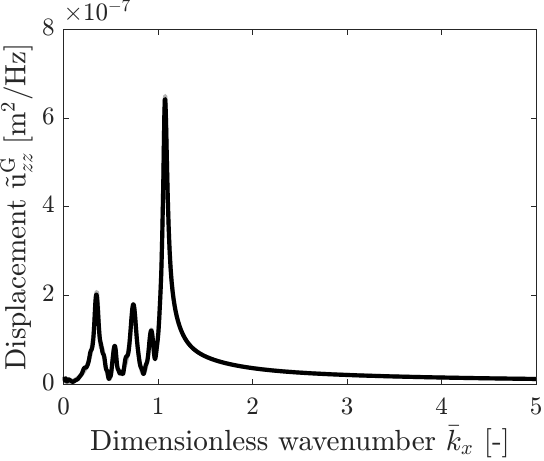} 
	\end{center}
	 \caption{P-SV Green's displacement $\tilde{u}^{\mathrm{G}}_{zz} (p_x, z, \omega)$ as a function of dimensionless wavenumber evaluated at $z = 0 \, \mathrm{m}$ and $f=10 \, \mathrm{Hz}$ for a vertical concentrated load applied at $x_\mathrm{s}=0 \, \mathrm{m}$ and $z_\mathrm{s}=0 \, \mathrm{m}$ using ROM (black) with (a) 10 modes, (b) 20 modes and (c) 40 modes, compared with FOM (gray).}\label{fig:ex1a_psv_kxd}
\end{figure}

\begin{figure}[!htb]
   \begin{center}
	 (a) \includegraphics[width=0.27\textwidth]{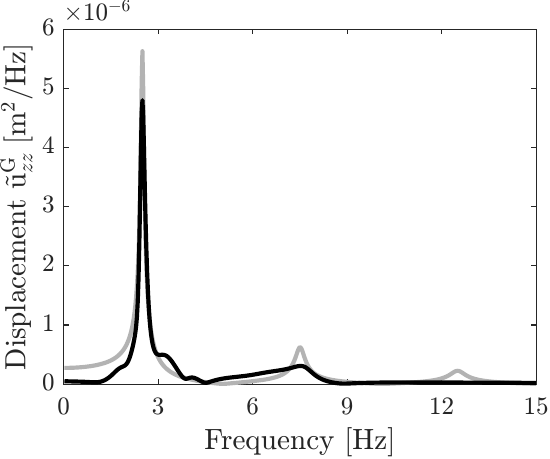}
	 (b) \includegraphics[width=0.27\textwidth]{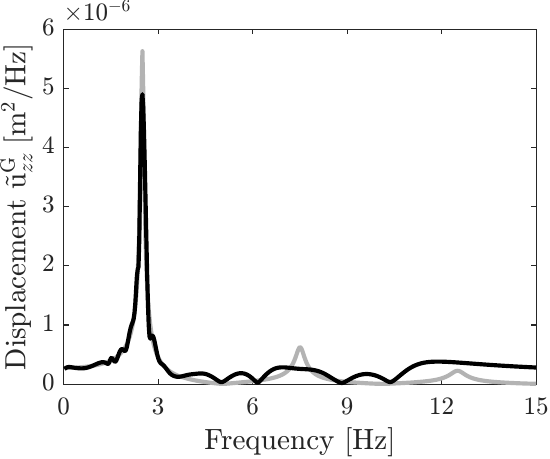}
	 (c) \includegraphics[width=0.27\textwidth]{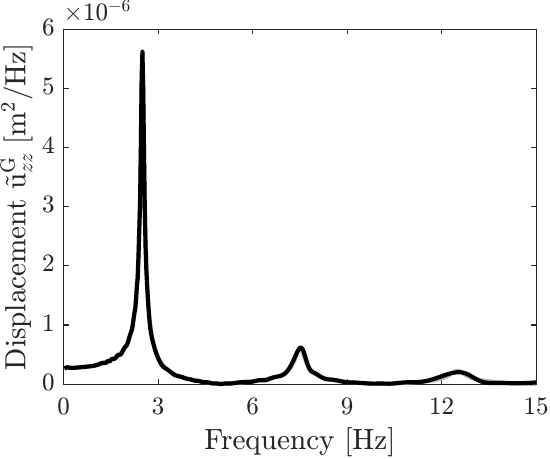} 
   \end{center}
	 \caption{P-SV Green's displacement $\tilde{u}^{\mathrm{G}}_{zz} (p_x, z, \omega)$ versus frequency evaluated at $\bar{k}_{x} = 0 \, \mathrm{rad/m}$ and $z = 0 \, \mathrm{m}$ for a vertical concentrated load applied at $x_\mathrm{s}=0 \, \mathrm{m}$ and $z_\mathrm{s}=0 \, \mathrm{m}$ using ROM (black) with (a) 10 modes, (b) 20 modes and (c) 40 modes, compared with FOM (gray).}
    \label{fig:ex1a_psv_ftilde}
\end{figure}

The analysis is now extended to the P-SV problem to illustrate the performance of the reduced order model for mixed compressional and shear motions. Figure~\ref{fig:ex1a_psv_f_kxd} shows the P–SV Green’s displacement $\tilde{u}^{\mathrm{G}}_{zz} (p_x, z, \omega)$ in the wavenumber–frequency domain evaluated at the surface of the soil layer ($z = 0 \, \mathrm{m}$) due to a vertical concentrated load applied at $x_\mathrm{s}=0 \, \mathrm{m}$ and $z_\mathrm{s}=0 \, \mathrm{m}$. Figure~\ref{fig:ex1a_psv_f_kxd}d shows the FOM reference solution. In this coupled problem, the response arises from the interaction between compressional and vertically polarized shear motions, producing both surface and body wave components. The wavenumber–frequency spectrum is dominated by the fundamental Rayleigh branch at low wavenumbers, while higher-order P–SV modes appear at larger wavenumbers and frequencies. Distinct peaks near $2.5 \, \mathrm{Hz}$ and $7.5 \, \mathrm{Hz}$ correspond to the first and second vertical eigenfrequencies of the clamped soil layer, equal to $C_\mathrm{p}/(4h)$ and $3C_\mathrm{p}/(4h)$, respectively. With increasing frequency, the influence of layering diminishes as the wavelengths become small compared to the layer thickness, and the response approaches that of a homogeneous halfspace. 
The ROM approximation obtained with 10 modes (figure~\ref{fig:ex1a_psv_f_kxd}a) yields a relatively coarse approximation. As the number of modes is increased to 20 (figure~\ref{fig:ex1a_psv_f_kxd}b) and 40 (figure~\ref{fig:ex1a_psv_f_kxd}c), the approximation quality improves progressively.
The ROM with 40 modes accurately reproduces the coupled dispersion characteristics, including both the Rayleigh mode and higher-order body wave contributions in close agreement with the FOM across the entire frequency range.
To provide further insight, figure~\ref{fig:ex1a_psv_kxd} compares the P–SV Green’s displacement $\tilde{u}^{\mathrm{G}}_{zz}(p_x, z, \omega)$ as a function of the dimensionless wavenumber $\bar{k}_x$ at $f=10 \, \mathrm{Hz}$ obtained using both ROM and FOM. The wavenumber-dependent behavior demonstrates multiple coupled P–SV branches, with the fundamental Rayleigh mode dominant at low wavenumbers and higher-order modes emerging at larger wavenumbers and higher frequencies.
Figure~\ref{fig:ex1a_psv_ftilde} shows the same displacement component as a function of frequency evaluated for $\bar{k}_{x} = 0 \, \mathrm{rad/m}$. The frequency-dependent displacement field is well captured by the ROM, which closely reproduces the modal peaks associated with the vertical eigenfrequencies of the clamped soil layer.

\begin{figure}[!htb]
   \begin{center}
	 (a) \includegraphics[width=0.25\textwidth]{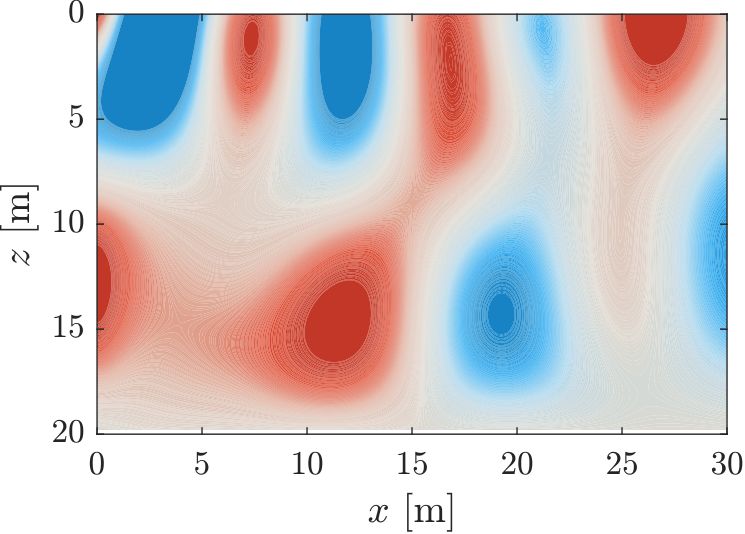}
	 (b) \includegraphics[width=0.25\textwidth]{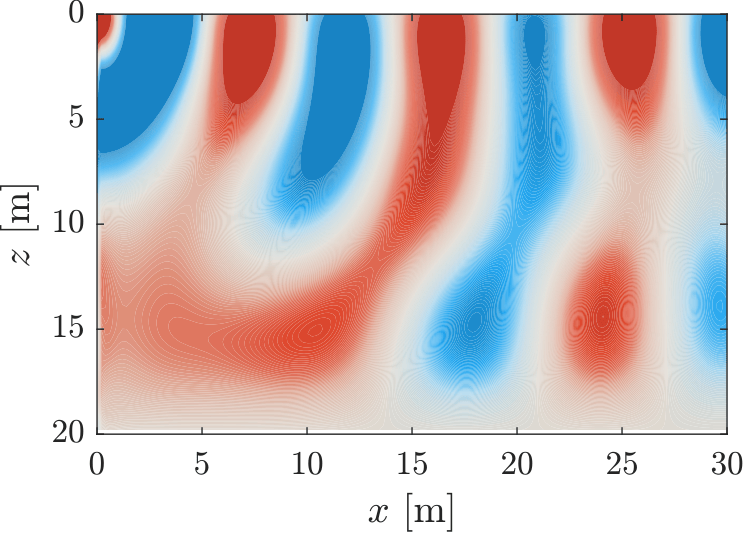} \\
	 (c) \includegraphics[width=0.25\textwidth]{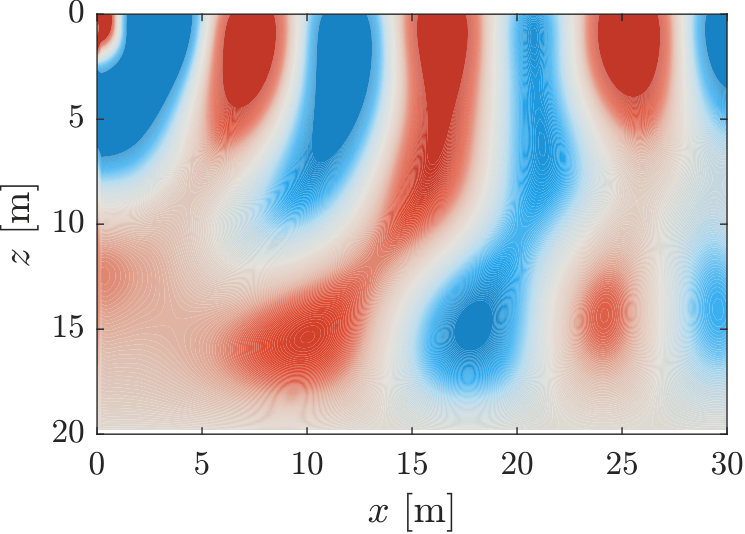} 
	 (d) \includegraphics[width=0.25\textwidth]{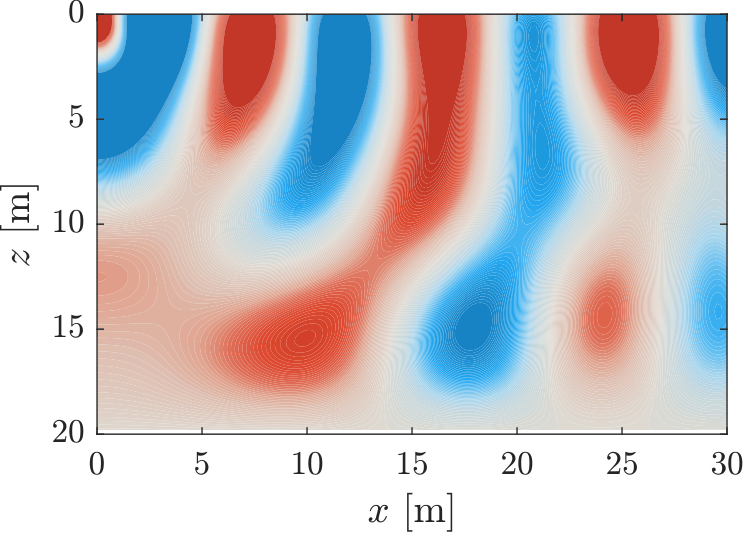} \\
	 \vspace{0.2cm} \hspace{0.6cm}
	     \includegraphics[width=0.27\textwidth]{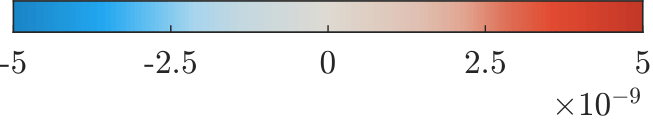} 
   \end{center}
	 \caption{P-SV Green's displacement $\hat{u}^{\mathrm{G}}_{zz} (x, z, \omega)$ [$\mathrm{m}/\mathrm{Hz}$] in the spatial-frequency domain at $f=10 \,\mathrm{Hz}$ for a vertical concentrated load applied at $x_\mathrm{s}=0 \, \mathrm{m}$ and $z_\mathrm{s}=0 \, \mathrm{m}$ using ROM with (a) 10 modes, (b) 20 modes, (c) 40 modes and (d) FOM.}
	 \label{fig:ex1_psv_contour_x}
\end{figure}

\begin{figure}[!htb]
   \begin{center}
	 (a) \includegraphics[width=0.27\textwidth]{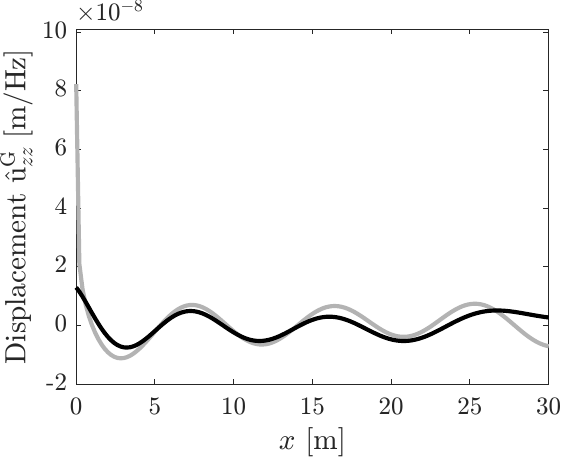}
	 (b) \includegraphics[width=0.27\textwidth]{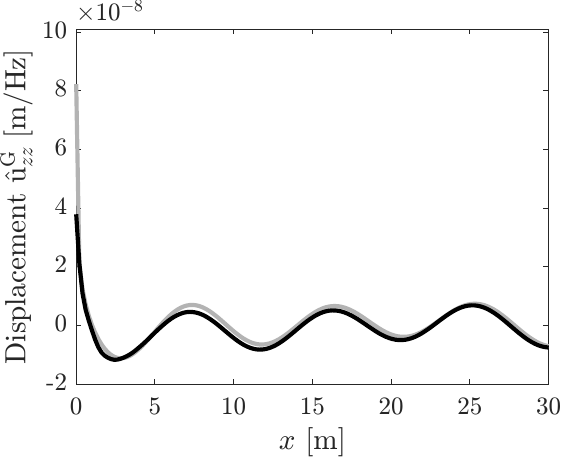}
	 (c) \includegraphics[width=0.27\textwidth]{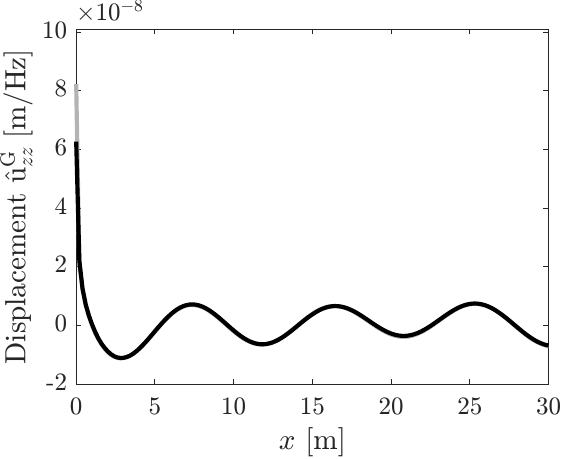} 
	\end{center}
	 \caption{P-SV Green's displacement $\hat{u}^{\mathrm{G}}_{zz} (x, z, \omega)$ as a function of $x$ evaluated at $z = 0 \, \mathrm{m}$ and $f=10 \,\mathrm{Hz}$ for a vertical concentrated load applied at $x_\mathrm{s}=0 \, \mathrm{m}$ and $z_\mathrm{s}=0 \, \mathrm{m}$ using ROM (black) with (a) 10 modes, (b) 20 modes and (c) 40 modes, compared with FOM (gray).}
	 \label{fig:ex1_psv_x}
\end{figure}

\begin{figure}[!htb]
   \begin{center}
	 (a) \includegraphics[width=0.27\textwidth]{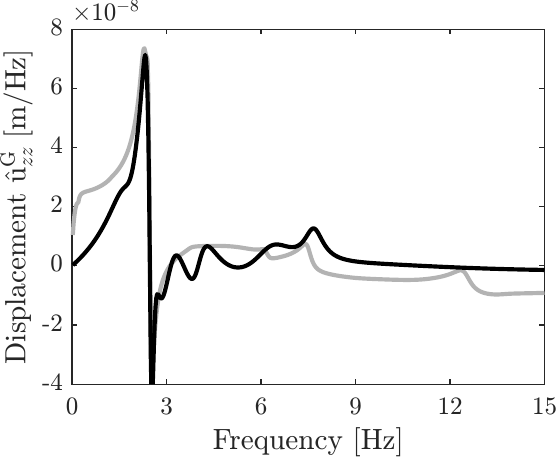}
	 (b) \includegraphics[width=0.27\textwidth]{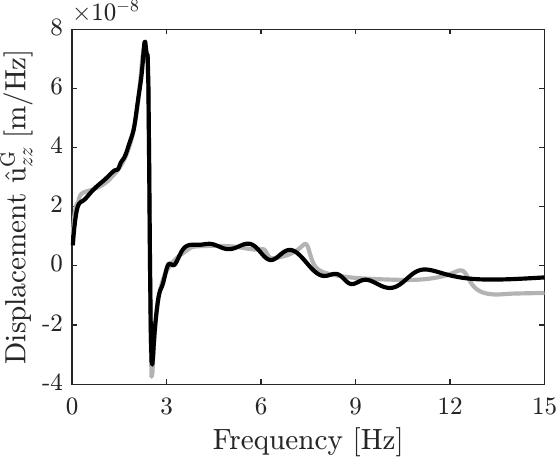}
	 (c) \includegraphics[width=0.27\textwidth]{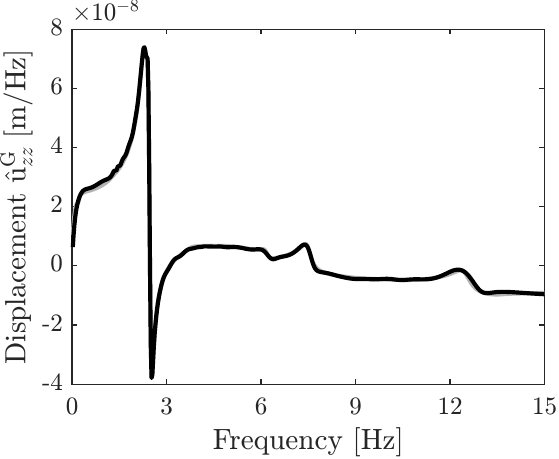} 
   \end{center}
	 \caption{P-SV Green's displacement $\hat{u}^{\mathrm{G}}_{zz} (x, z, \omega)$ as a function of frequency evaluated at $x = 1 \, \mathrm{m}$ and $z = 0 \, \mathrm{m}$ for a vertical concentrated load applied at $x_\mathrm{s}=0 \, \mathrm{m}$ and $z_\mathrm{s}=0 \, \mathrm{m}$ using ROM (black) with (a) 10 modes, (b) 20 modes and (c) 40 modes, compared with FOM (gray).}\label{fig:ex1_psv_fx}
\end{figure}

As for the SH problem, Green’s displacements in the spatial–frequency domain are obtained through an inverse logarithmic Fourier transform applied to the low-rank representation. 
Figure~\ref{fig:ex1_psv_contour_x} shows the P–SV Green’s displacement at $f = 10 \, \mathrm{Hz}$. The response exhibits oscillatory patterns characteristic of coupled compressional-shear motion, with amplitudes decaying away from the source. 
At this frequency, the Rayleigh wavelength is expressed as $\lambda_\mathrm{r} = C_\mathrm{r}/f$. The Rayleigh wave velocity $C_\mathrm{r}$ can be approximated in terms of the Poisson's ratio $\nu$ as \cite{vikt67a}:
\begin{equation}
C_\mathrm{r} \approx C_\mathrm{s}\frac{0.87 + 1.12 \nu}{1 + \nu}
\label{eq:rayleigh_wavelength}
\end{equation}
For the present configuration, this yields $C_\mathrm{r} = 93.2 \, \mathrm{m/s}$ and a Rayleigh wavelength of $\lambda_\mathrm{r} = 9.31 \, \mathrm{m}$ at  $10 \, \mathrm{Hz}$. Since the Rayleigh wavelength is small compared to the layer thickness, the associated penetration depth remains confined to the upper part of the layer and the rigid bedrock has only a limited effect on the surface wave motion. This behavior is consistent with the oscillation pattern observed in the displacement field. Figure~\ref{fig:ex1_psv_x} further illustrates this behavior by showing the vertical displacement $\hat{u}^{\mathrm{G}}_{zz}(x, z, \omega)$ along $z=0 \, \mathrm{m}$ at $f = 10 \,\mathrm{Hz}$ as a function of distance $x$.
 The ROM results progressively approach the FOM, and with 40 modes, both amplitude and phase match well. Figure~\ref{fig:ex1_psv_fx} shows the same displacement component as a function of frequency for a receiver positioned at $x = 1 \,\mathrm{m}$ and $z = 0 \,\mathrm{m}$. The ROM closely follows the FOM across the frequency range, demonstrating consistent agreement between the two models. Although the P–SV problem involves roughly twice the system size of the SH problem due to coupling between displacement components, the ROM reproduces the reference results with comparable accuracy.

The ROM solution is stored in memory as a compact combination of a core tensor and factor matrices corresponding to the separated directions, allowing efficient reconstruction for any parameter configuration \cite{bade25a}. The corresponding storage requirements are compared in table~\ref{tab:storage_comparison}. This comparison highlights the strong compression capability of the ROM, achieving storage reductions of nearly three orders of magnitude while preserving the accuracy and physical fidelity of the full order Green’s function within the prescribed tolerance limits.

The computation time of the FOM and ROM are summarized in table~\ref{tab:time_comparison_domains}. For both the SH and P–SV problems, the offline construction of the ROM requires higher computation time than the corresponding FOM computation in the wavenumber-frequency domain. Once constructed, however, the ROM exhibits markedly lower reassembly times, as the solution is recovered through low-rank tensor contractions rather than through the assembly and solution of large linear systems. In the spatial–frequency domain, the ROM is faster in both the offline and online phases, owing to the reduced number of inverse Fourier transforms required in comparison with the FOM. These results quantify the trade-off between the higher offline construction cost of the ROM and its reduced cost for subsequent evaluations. Moreover, the offline ROM representation is highly compact in storage and can be reassembled rapidly for subsequent use, without recomputing the full order solution.

\begin{table}[!h]
	\centering
	\caption{Storage requirement for ROM and FOM in the layer on bedrock configuration.}
	\label{tab:storage_comparison}
	\begin{tabular}{l c c c c}
		\toprule
		& \multicolumn{2}{c}{Wavenumber-frequency domain}
		& \multicolumn{2}{c}{Spatial-frequency domain} \\
		\cmidrule(lr){2-3} \cmidrule(lr){4-5}
		Case & ROM & FOM & ROM & FOM \\
		& [MB] & [MB] & [MB] & [MB] \\
		\midrule
		SH & 1.7 & 1373.0 & 1.9 & 200.0 \\
		P-SV & 2.5 & 2700.0 & 48.0 & 458.0 \\
		\bottomrule
	\end{tabular}
\end{table}

\begin{table}[!h]
	\centering
	\caption{Computation time for ROM and FOM in the layer on bedrock configuration.}
	\label{tab:time_comparison_domains}
	\begin{tabular}{l c c c c c c}
		\toprule
		& \multicolumn{3}{c}{Wavenumber-frequency domain} 
		& \multicolumn{3}{c}{Spatial-frequency domain} \\
		\cmidrule(lr){2-4} \cmidrule(lr){5-7}
		Case 
		& ROM offline & ROM online & FOM 
		& ROM offline & ROM online & FOM \\
		& [s] & [s] & [s] & [s] & [s] & [s] \\
		\midrule
		SH   
		& 947.0  & 0.55 & 15.0 
		& 230.0  & 0.4  & 1636.0 \\
		P-SV 
		& 12072.0 & 1.50 & 23.7 
		& 298.0  & 0.7  & 7235.0 \\
		\bottomrule
	\end{tabular}
\end{table}

\subsection{Layered soil on bedrock with separation of a material parameter}
\label{subsec:ex1d}
This example considers a layered soil resting on rigid bedrock, where the shear modulus \(\mu\) is introduced as a separable dimension in the reduced order formulation. The objective is to extend the low-rank representation of the Green’s function to account for variations in material stiffness directly within the decomposition, thereby avoiding the need to reassemble the global matrices and solve the full order system for each material configuration. The thin layer matrices \(\mathbf{A}\), \(\mathbf{B}\), and \(\mathbf{G}\) are expressed as linear functions of the Lamé constants \(\mu\) and \(\lambda\), while the mass matrix \(\mathbf{M}\) remains independent of this decomposition:

\begin{equation}
\begin{aligned}
\mathbf{A}(\mu,\lambda) &= \mu\,\mathbf{A}_{\mu} + \lambda\,\mathbf{A}_{\lambda}\\
\mathbf{B}(\mu,\lambda) &= \mu\,\mathbf{B}_{\mu} + \lambda\,\mathbf{B}_{\lambda}\\
\mathbf{G}(\mu,\lambda) &= \mu\,\mathbf{G}_{\mu} + \lambda\,\mathbf{G}_{\lambda}
\end{aligned}
\label{eq:lame}
\end{equation}
This reformulation enables the shear modulus \(\mu\) to be treated as an explicit separable coordinate in the reduced order algorithm. In the present study, the Lam\'e constant $\lambda$ is fixed, although the formulation readily permits $\lambda$ to be introduced as an additional separable parameter, allowing multiple material properties to be considered within a single reduced order framework.
The shear modulus of the soil layer is sampled over the range $\mu \in [18.0,\,26.0]\,\mathrm{MPa}$ using 21 uniformly spaced values.
The linear dependence of the thin layer matrices on the Lamé parameters allows the contributions associated with \(\mu\) and \(\lambda\) to be separated, as shown in equation~\eqref{eq:lame}. To this end, the matrices corresponding to unit values of the Lamé parameters, namely  \((\mu, \lambda) = (1, 0)\) and \((\mu, \lambda) = (0, 1)\) are computed once. These matrices represent the individual contributions of the shear and volumetric stiffness, respectively. For any given material configuration, the thin layer matrices are then obtained by a linear combination of these precomputed unit matrices, scaled by the desired values of $\mu$ and $\lambda$. This linear parameterization is directly exploited within the GTA framework by treating the shear modulus $\mu$ as an explicit separable coordinate in the reduced order algorithm. Consequently, variations in material stiffness are handled within the same low-rank approximation, without reassembling the global matrices for each set of material parameters, and the precomputed unit matrices can be reused consistently across the entire wavenumber–frequency grid.
From a computational perspective, this approach decouples the material dependence from the frequency–wavenumber solution space. Once the reduced bases are constructed, the system response can be evaluated at any value of $\mu$ within the sampled range, without reassembling or resolving the full system.

\begin{figure}[!htb]
   \begin{center}
     (a) \includegraphics[width=0.27\textwidth]{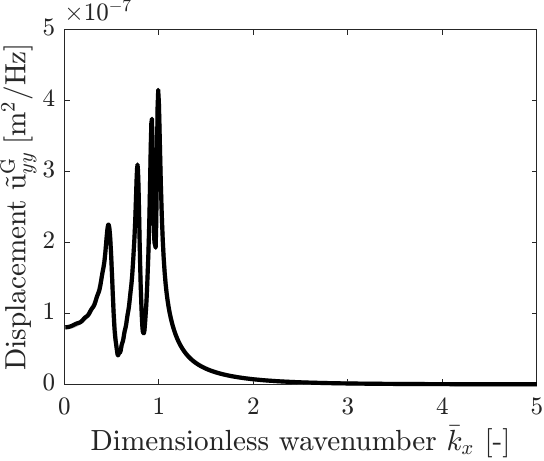}	
	 (b) \includegraphics[width=0.27\textwidth]{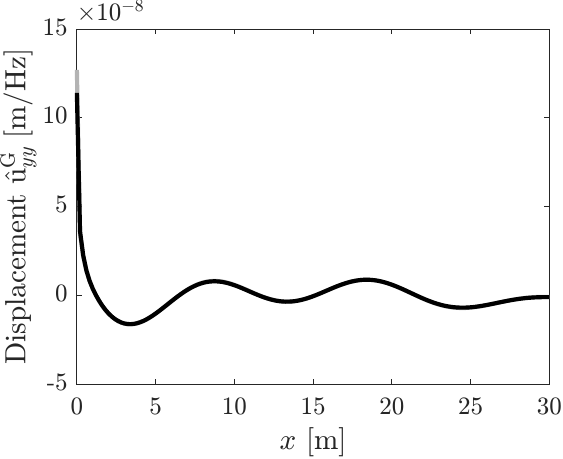} 
   \end{center}
	 \caption{SH Green’s displacement (a) $\tilde{u}^{\mathrm{G}}_{yy}(p_x,z,\mu)$ as a function of dimensionless wavenumber and (b) $\hat{u}^{\mathrm{G}}_{yy}(x,z,\mu)$ as a function of $x$, evaluated at $z=0 \, \mathrm{m}$, $\mu = 18~\mathrm{MPa}$ and $f=10 \; \mathrm{Hz}$ for an out-of-plane concentrated load applied at $x_\mathrm{s}=0 \, \mathrm{m}$ and $z_\mathrm{s}=0 \, \mathrm{m}$ using ROM using 30 modes (black) and FOM (gray).}\label{fig:ex1d_sh_a}
\end{figure}

\begin{figure}[!htb]
   \begin{center}
	 (a) \includegraphics[width=0.27\textwidth]{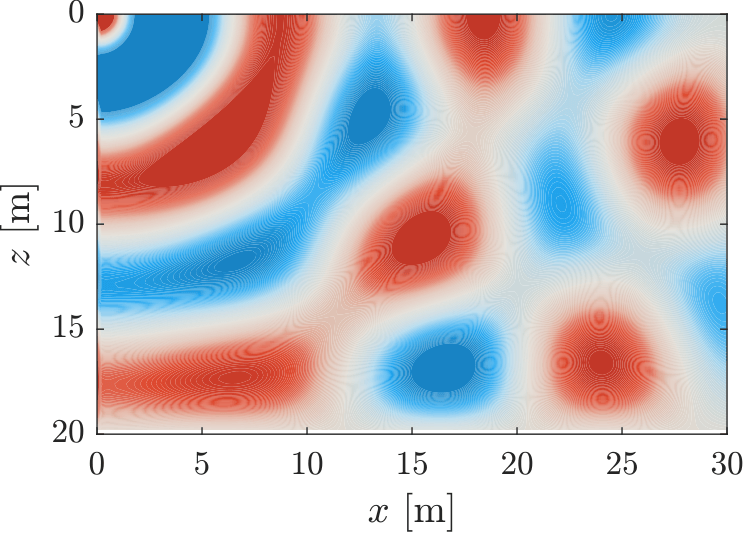}	
     (b) \includegraphics[width=0.27\textwidth]{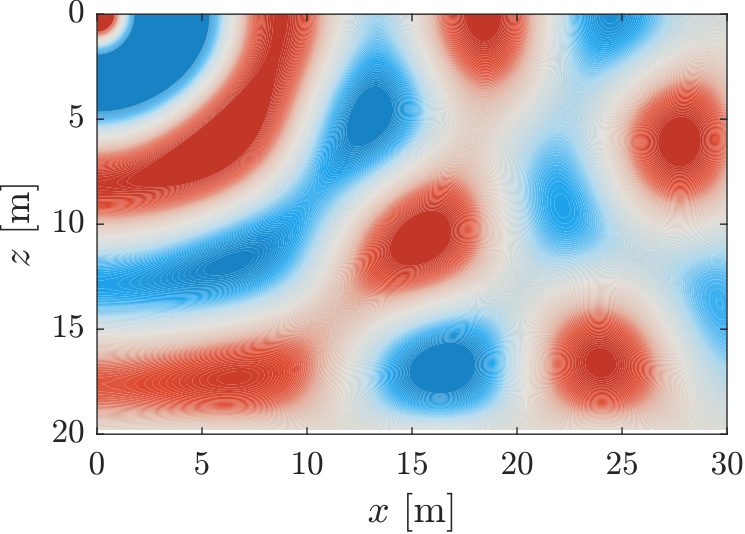} \\
     \vspace{0.2cm} \hspace{0.6cm}
	     \includegraphics[width=0.27\textwidth]{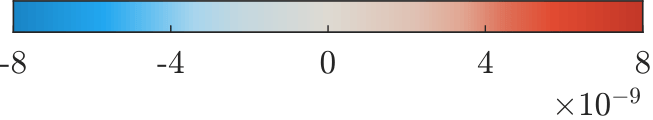} 
   \end{center}
	 \caption{SH Green’s displacement $\hat{u}^{\mathrm{G}}_{yy}(x,z,\mu)$ [$\mathrm{m}/\mathrm{Hz}$] evaluated at $\mu = 18~\mathrm{MPa}$ and $f=10 \; \mathrm{Hz}$ for an out-of-plane concentrated load applied at $x_\mathrm{s}=0 \, \mathrm{m}$ and $z_\mathrm{s}=0 \, \mathrm{m}$ using (a) ROM with 30 modes and (b) FOM.}\label{fig:ex1d_sh_b}
\end{figure}

Figure~\ref{fig:ex1d_sh_a}a shows the SH Green’s displacement $\hat{u}^{\mathrm{G}}_{yy}(x,z,\mu)$ in the wavenumber domain evaluated at $z=0 \, \mathrm{m}$, $\mu = 18~\mathrm{MPa}$ and $f=10 \; \mathrm{Hz}$ due to a horizontal load applied at the surface of the soil layer ($x_\mathrm{s} = 0 \, \mathrm{m}$, $z_\mathrm{s} = 0 \, \mathrm{m}$). Distinct peaks corresponding to the fundamental and higher-order Love modes in the soil layer are observed. The corresponding wavefield of $\hat{u}^{\mathrm{G}}_{yy}(x,z,\mu)$ in the spatial-frequency domain is shown in figure~\ref{fig:ex1d_sh_b}. The displacement exhibits an oscillatory pattern with a wavelength of approximately $10 \, \mathrm{m}$. The displacement amplitude decreases with distance from the source due to material damping. Figure~\ref{fig:ex1d_sh_a}b presents a slice of the spatial wavefield evaluated at $z=0 \, \mathrm{m}$, providing a complementary view of the spatial response.  The ROM with 30 modes reproduces these characteristics with high accuracy.

\begin{figure}[!htb]
   \begin{center}
	 (a) \includegraphics[width=0.27\textwidth]{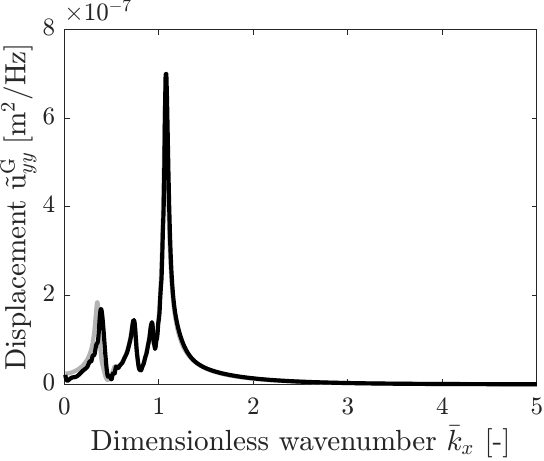}
	 (b) \includegraphics[width=0.27\textwidth]{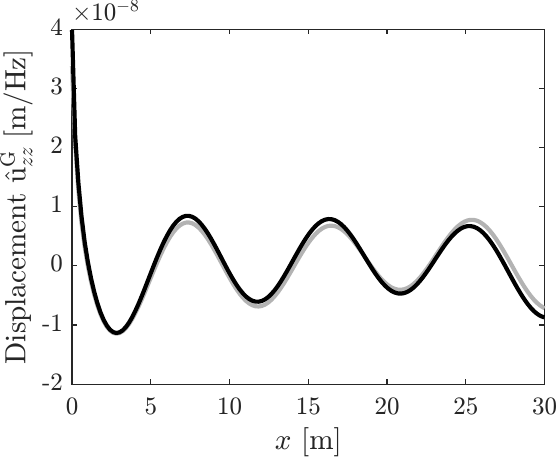} 
   \end{center}
	 \caption{P-SV Green’s displacement (a) $\tilde{u}^{\mathrm{G}}_{zz}(p_x,z,\mu)$ as a function of dimensionless wavenumber and (b) $\hat{u}^{\mathrm{G}}_{zz}(x,z,\mu)$ as a function of $x$, evaluated at $z=0 \, \mathrm{m}$, $\mu = 18~\mathrm{MPa}$  and $f=10 \; \mathrm{Hz}$ for a vertical concentrated load applied at $x_\mathrm{s}=0 \, \mathrm{m}$ and $z_\mathrm{s}=0 \, \mathrm{m}$ using ROM (black) and FOM (gray).}
\label{fig:ex1d_psv_a}
\end{figure}

\begin{figure}[!htb]
   \begin{center}
	 (a) \includegraphics[width=0.27\textwidth]{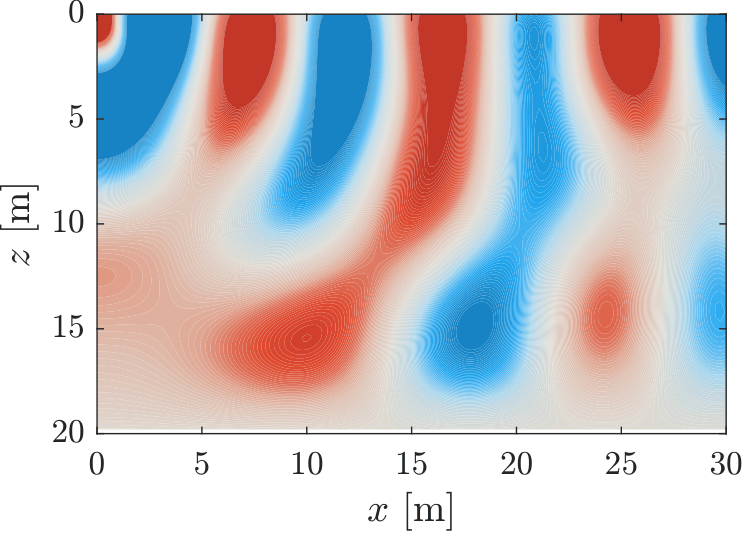}
	 (b) \includegraphics[width=0.27\textwidth]{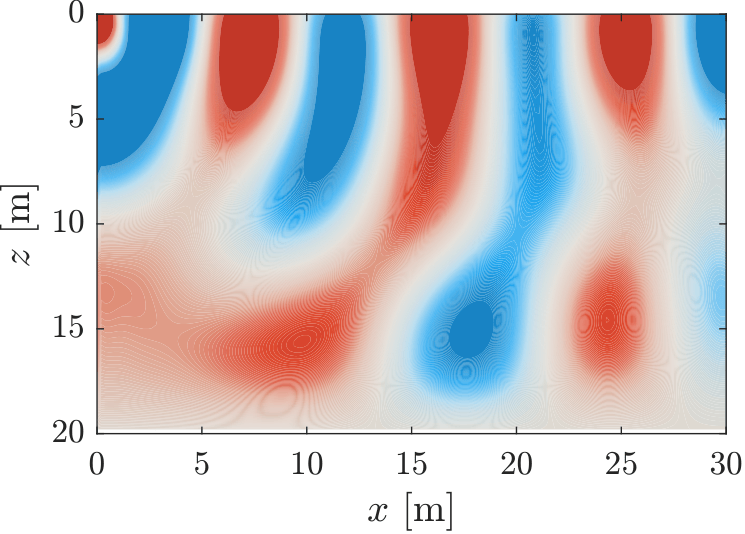} \\
	  \vspace{0.2cm} \hspace{0.6cm}
	     \includegraphics[width=0.27\textwidth]{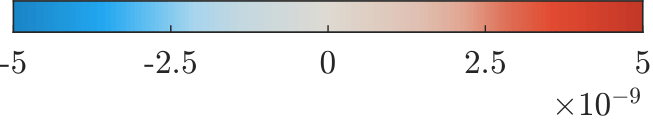}  
   \end{center}
	 \caption{P-SV Green’s displacement $\hat{u}^{\mathrm{G}}_{zz}(x, z, \mu)$ [$\mathrm{m}/\mathrm{Hz}$] evaluated at $\mu = 18~\mathrm{MPa}$ and $f = 10~\mathrm{Hz}$ for a vertical concentrated load applied at $x_\mathrm{s}=0 \, \mathrm{m}$ and $z_\mathrm{s}=0 \, \mathrm{m}$ using (c) ROM with 35 modes and(d) FOM.}
\label{fig:ex1d_psv_b}
\end{figure}

For the P–SV case, figure~\ref{fig:ex1d_psv_a}a illustrates the Green’s displacement $\hat{u}^{\mathrm{G}}_{zz}(x,z,\mu)$ in the wavenumber domain evaluated at $z=0 \, \mathrm{m}$, $\mu = 18~\mathrm{MPa}$  and $f=10 \; \mathrm{Hz}$. Figure~\ref{fig:ex1d_psv_b} shows the corresponding wavefield in the spatial domain. Figure~\ref{fig:ex1d_psv_a}b demonstrates the splice through the spatial-frequency waveplot computed at $z=0 \, \mathrm{m}$. The coupled compressional–shear response is dominated near the surface by Rayleigh wave propagation and decays with depth into the soil layer. The horizontal oscillation pattern corresponds to a Rayleigh wavelength of approximately $\lambda_\mathrm{r} \approx 9.32~\mathrm{m}$, with amplitude decay is governed by material damping. These features are accurately captured by the ROM with 35 modes. These results are in agreement with example~\ref{subsec:ex1}, demonstrating that the ROM maintains its accuracy when a material parameter is introduced as a separable coordinate and can efficiently reconstruct the Green’s function for any specified configuration without reassembling the system matrices.

\begin{table}[!htb]
	\centering
	\caption{Storage requirement for ROM and FOM in layer on bedrock case with material parameter separation.}
	\label{tab:storage_comparison_mu}
	\begin{tabular}{l c c c c}
		\toprule
		& \multicolumn{2}{c}{Wavenumber-frequency domain}
		& \multicolumn{2}{c}{Spatial-frequency domain} \\
		\cmidrule(lr){2-3} \cmidrule(lr){4-5}
		Case & ROM & FOM & ROM & FOM \\
		& [MB] & [MB] & [MB] & [MB] \\
		\midrule
		SH & 0.52 & 1.53 & 0.18 & 0.38 \\
		P-SV & 0.75 & 6.10 & 0.35 & 1.53 \\
		\bottomrule
	\end{tabular}
\end{table}

\begin{table}[!htb]
	\centering
	\caption{Computation time for ROM and FOM in the layer on bedrock case with material parameter separation.}
	\label{tab:time_comparison_mu}
	\begin{tabular}{l c c c c c c}
		\toprule
		& \multicolumn{3}{c}{Wavenumber-frequency domain}
		& \multicolumn{3}{c}{Spatial-frequency domain} \\
		\cmidrule(lr){2-4} \cmidrule(lr){5-7}
		Case 
		& ROM offline & ROM online & FOM
		& ROM offline & ROM online & FOM \\
		& [s] & [s] & [s] & [s] & [s] & [s] \\
		\midrule
		SH   
		& 140  & 0.08  & 0.07
		& 0.077 & 0.018 & 0.11 \\
		P-SV 
		& 229  & 0.03  & 0.08
		& 0.094 & 0.023 & 0.14 \\
		\bottomrule
	\end{tabular}
\end{table}

The storage requirements and computation time for the material-parametric SH and P-SV problems are summarized in tables~\ref{tab:storage_comparison_mu} and~\ref{tab:time_comparison_mu}, respectively. The ROM requires a substantially smaller memory than the corresponding FOM, owing to its low-rank representation based on factor matrices and a compact core. In terms of computational effort, the offline construction of the ROM entails a higher upfront cost than the total computation time required for FOM evaluations over all sampled values of the shear modulus. The ROM allows for rapid online reconstruction of the system response for any configuration within the sampled range. This enables efficient reuse of the parametric solution for different material configurations without recomputing the full order system. 

\begin{figure}[!htb]
   \begin{center}
	 (a) \includegraphics[width=0.27\textwidth]{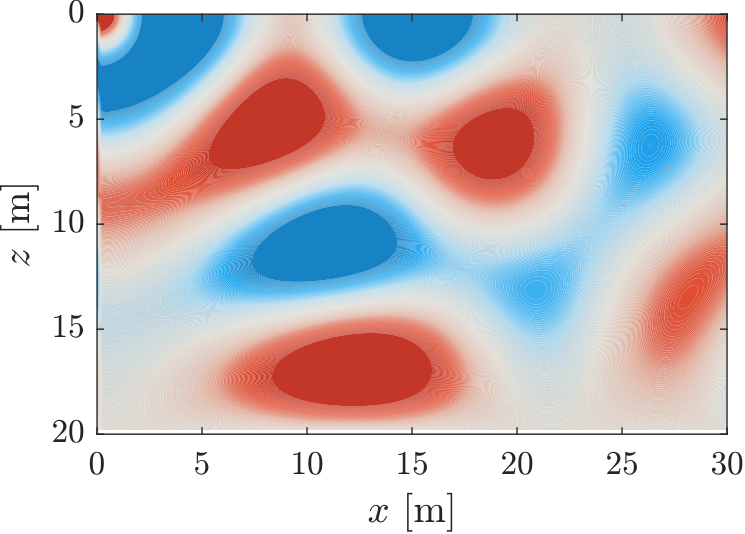}
	 (c) \includegraphics[width=0.27\textwidth]{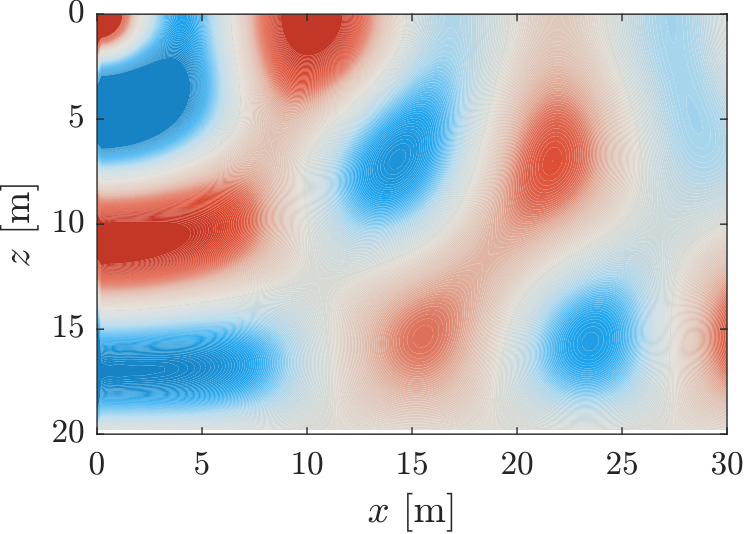}\\
	 \vspace{0.2cm} \hspace{0.6cm}
	     \includegraphics[width=0.25\textwidth]{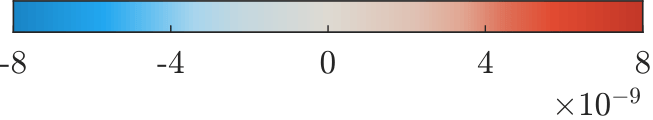}  
   \end{center}
	 \caption{SH Green’s displacement $\hat{u}^{\mathrm{G}}_{yy}(x,z,\mu)$ [$\mathrm{m}/\mathrm{Hz}$] evaluated at $f=10 \; \mathrm{Hz}$ for shear moduli $\mu$ (a) $22~\mathrm{MPa}$ and (b) $26~\mathrm{MPa}$ due to an out-of-plane concentrated load applied at $x_\mathrm{s}=0 \, \mathrm{m}$ and $z_\mathrm{s}=0 \, \mathrm{m}$, computed using ROM with 30 modes.}
\label{fig:ex1d_demo_sh}
\end{figure}

\begin{figure}[!htb]
   \begin{center}
	 (a) \includegraphics[width=0.27\textwidth]{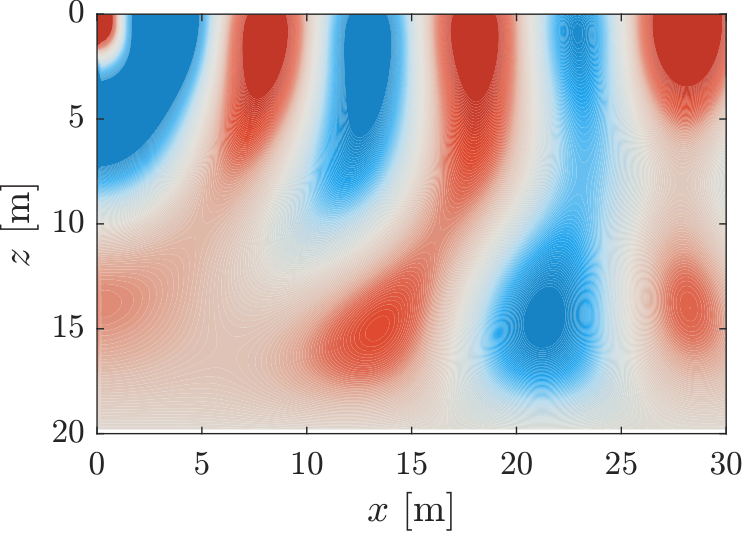}
	 (c) \includegraphics[width=0.27\textwidth]{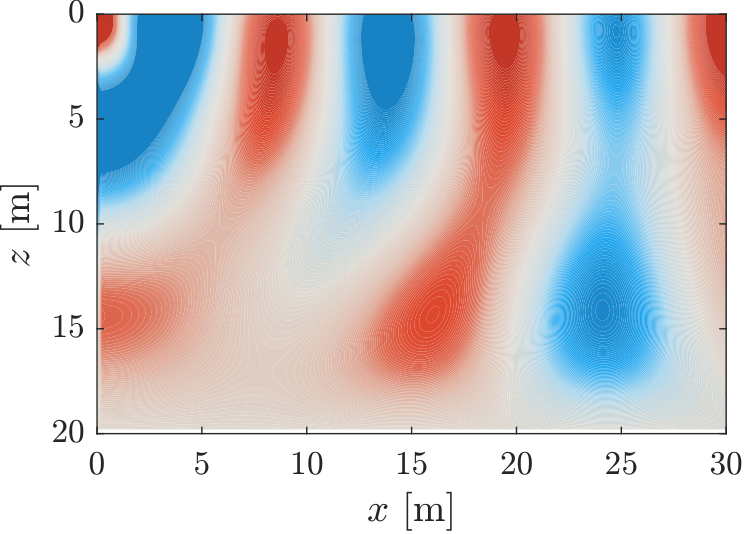}\\
	  \vspace{0.2cm} \hspace{0.6cm}
	     \includegraphics[width=0.25\textwidth]{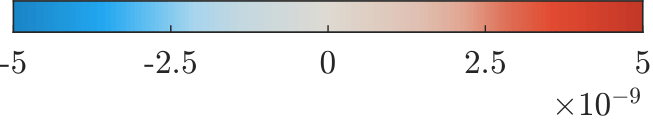}   
   \end{center}
	 \caption{P-SV Green’s displacement $\hat{u}^{\mathrm{G}}_{zz}(x,z,\mu)$ [$\mathrm{m}/\mathrm{Hz}$] evaluated at $f=10 \; \mathrm{Hz}$ for shear moduli $\mu$ (a) $22~\mathrm{MPa}$ and (b) $26~\mathrm{MPa}$ due to a vertical concentrated load applied at $x_\mathrm{s}=0 \, \mathrm{m}$ and $z_\mathrm{s}=0 \, \mathrm{m}$, computed  using ROM with 35 modes.}
\label{fig:ex1d_demo_psv}
\end{figure}

To illustrate the computational advantage of the ROM in a parametric setting, a case is considered in which the shear modulus of the soil is varied to represent different soil stiffness conditions. Figure~\ref{fig:ex1d_demo_sh} presents the SH Green’s displacement field $\hat{u}^{\mathrm{G}}_{yy}(x,z,\mu)$ for shear moduli $\mu = 22$ and $26~\mathrm{MPa}$ obtained using ROM. These results have been verified against the corresponding FOM solutions to ensure accuracy. For a fixed excitation frequency, increasing the shear modulus leads to higher shear wave velocities, resulting in longer shear wavelengths. In the present example, the considered range of shear moduli corresponds to shear wavelengths increasing from $10~\mathrm{m}$ for $\mu = 18~\mathrm{MPa}$ to about $12~\mathrm{m}$ for $\mu = 26~\mathrm{MPa}$. As the shear modulus increases, the displacement amplitudes decrease and the overall response becomes less pronounced.
The corresponding P–SV case is shown in Figure~\ref{fig:ex1d_demo_psv}, where the vertical displacement component $\hat{u}^{\mathrm{G}}_{zz}(x,z,\mu)$ is reported for the same values of shear moduli. At the free surface, the P–SV response is dominated by Rayleigh waves, whose wavelength increases with increasing shear modulus and can be approximated using equation~\eqref{eq:rayleigh_wavelength}. For the considered material properties and excitation frequency, the associated Rayleigh wavelengths increase from approximately $9.3~\mathrm{m}$ for $\mu = 18~\mathrm{MPa}$ to about $11.1~\mathrm{m}$ for $\mu = 26~\mathrm{MPa}$.
Within the ROM framework, the responses associated with different shear moduli are recovered directly from the precomputed offline reduced order representation. In contrast, the FOM requires a separate assembly and solution of the full system for each parameter value. 
By embedding material dependence directly within a low-rank representation, the reduced order formulation enables efficient and rapid evaluation of the response across the parameter space.

\subsection{Groene Hart site}
\label{subsec:ex2}
This section illustrates the application of the proposed reduced order algorithm to a practical case study, in which the Green’s functions are approximated for the Groene Hart site in the Netherlands. The corresponding soil properties are summarized in table~\ref{tab:soilpropterties}.The analysis considers SH and P-SV wave propagation with separation over slowness $p_x$, receiver depth $z$ and frequency $\omega$  using the same decomposition framework as in example~\ref{subsec:ex1}. The soil profile consists of several layers underlain by a homogeneous halfspace. While the layers are modeled through thin layer matrices, the halfspace contribution is represented analytically using an impedance relation that links displacements and tractions at the interface (eq.~\eqref{eq:Ks}). This impedance term is appended to the dynamic stiffness matrix of the layered system, as detailed in section~\ref{subsec:stiffness}.

For the SH case, the frequency range is sampled from $f = 0$ to $10\,\mathrm{Hz}$ with a bin $\Delta f = 0.025\,\mathrm{Hz}$.
The wavenumber domain is sampled in terms of slowness and discretized using a logarithmic grid with $1500$ points over the range $\bar{k}_{x} \in [10^{-2},\,10^{1}]$, with the corresponding slowness defined as $p_x = \bar{k}_{x} / C_{\mathrm{s}}^\mathrm{min}$.  The receiver depth is sampled using the thin layer mesh, extending uniformly down to $z = 28\,\mathrm{m}$ at a resolution of $0.05\,\mathrm{m}$. To enable the subsequent evaluation of the spatial–frequency response, the horizontal receiver coordinate is sampled as $x \in [0, \,50]\,\mathrm{m}$ using $250$ uniformly spaced points. A concentrated load is applied at the source position $z_{\mathrm{s}} = 0$ and $x_{\mathrm{s}} = 0$. The greedy Petrov--Galerkin algorithm is employed.
The stopping criterion for the greedy algorithm is set to a tolerance of $10^{-5}$ with a maximum of 40 modes.

\begin{figure}[!htb]
   \begin{center}
	 (a) \includegraphics[width=0.27\textwidth]{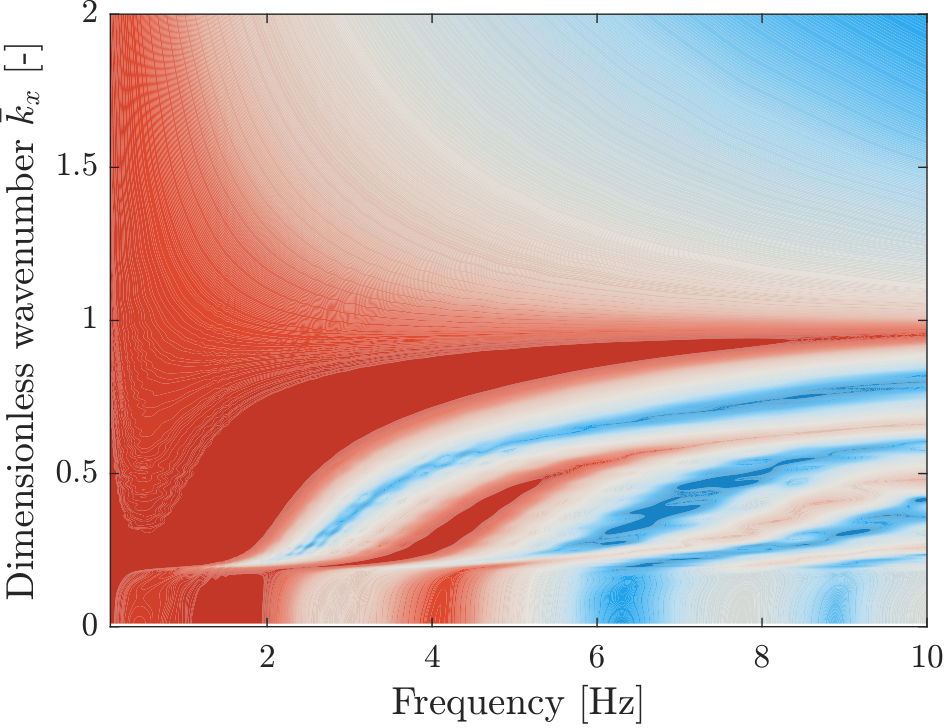}
	 (b) \includegraphics[width=0.27\textwidth]{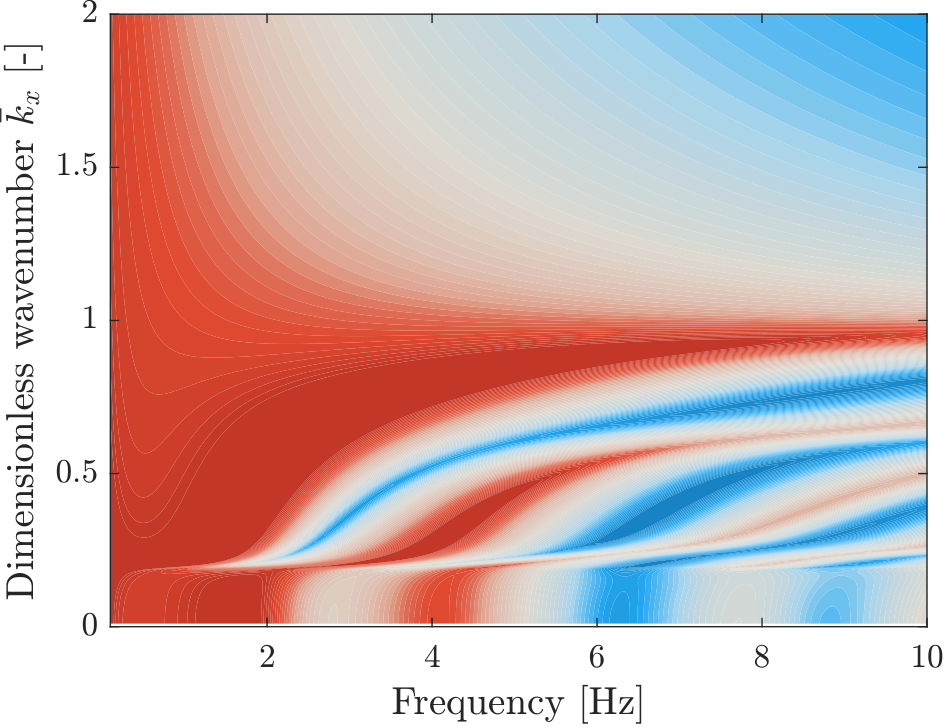} 
	 (c) \includegraphics[width=0.27\textwidth]{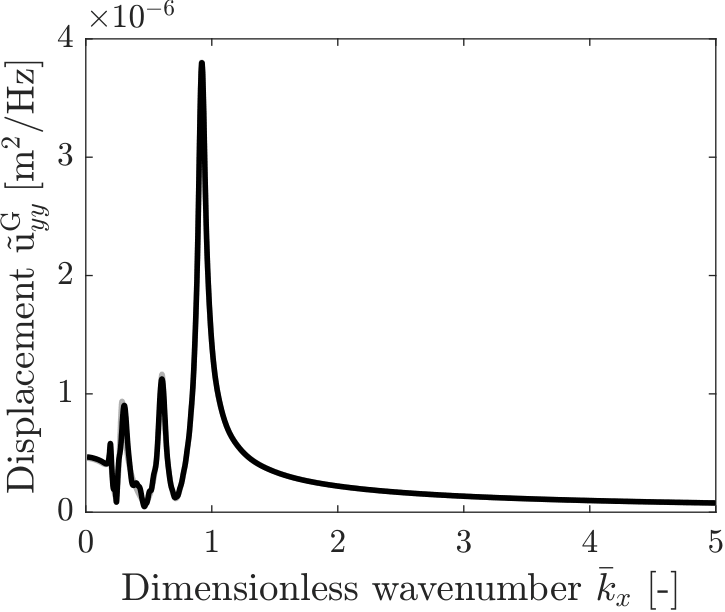}\\
     \vspace{0.2cm} \hspace{-4cm}
	     \includegraphics[width=0.27\textwidth]{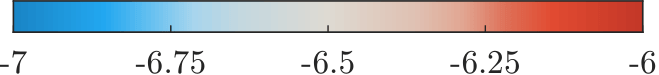}
   \end{center}
	 \caption{SH Green’s displacements $\tilde{u}^{\mathrm{G}}_{yy}(p_x, z, \omega)$ [$\mathrm{m}^2/\mathrm{Hz}$] in the wavenumber-frequency domain evaluated at $z=0 \, \mathrm{m}$ for an out-of-plane concentrated load applied at $x_\mathrm{s}=0 \, \mathrm{m}$ and $z_\mathrm{s}=0 \, \mathrm{m}$ using (a) ROM with 40 modes, (b) FOM and (c) as a function of dimensionless wavenumber for $f = 7.5\, \mathrm{Hz}$ using ROM with 40 modes (black) and FOM (gray).}
\label{fig:ex2_sh1}
\end{figure}

\begin{figure}[!htb]
   \begin{center}
	 (a) \includegraphics[width=0.27\textwidth]{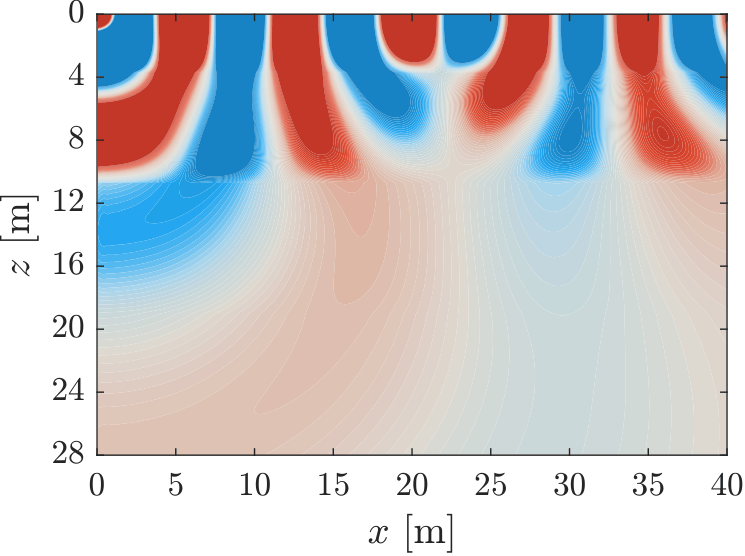}	
     (b) \includegraphics[width=0.27\textwidth]{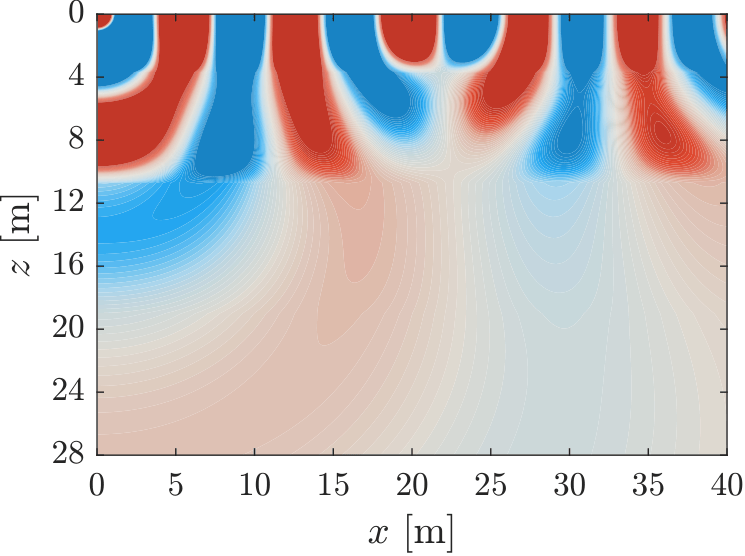}
     (c) \includegraphics[width=0.27\textwidth]{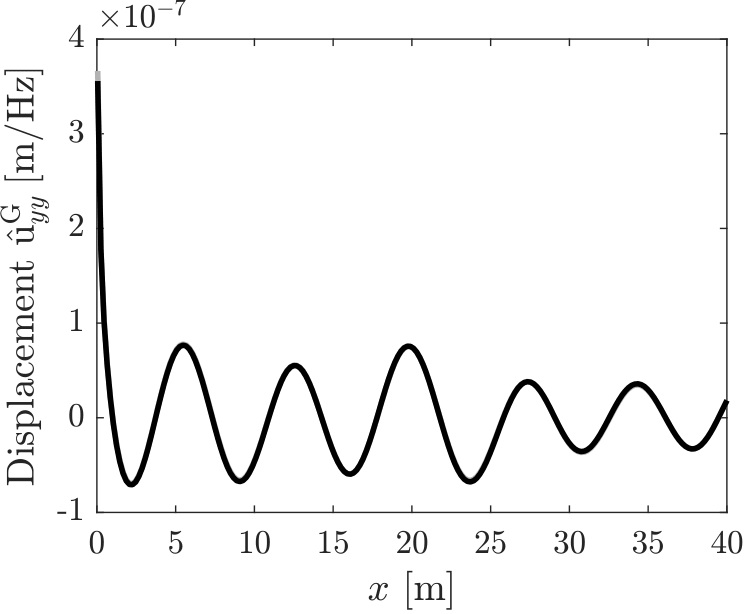}\\
     \vspace{0.2cm} \hspace{-4cm}
	     \includegraphics[width=0.27\textwidth]{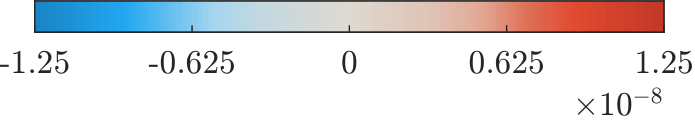}
   \end{center}
	 \caption{SH Green’s displacements $\hat{u}^{\mathrm{G}}_{yy}(x, z, \omega)$ [$\mathrm{m}/\mathrm{Hz}$] in the spatial-frequency domain at $f=7.5 \, \mathrm{Hz}$ for an out-of-plane concentrated load applied at $x_\mathrm{s}=0 \, \mathrm{m}$ and $z_\mathrm{s}=0 \, \mathrm{m}$ using (a) ROM with 40 modes, (b) FOM  and (c) as a function of $x$ evaluated for $z = 0\, \mathrm{m}$ using ROM with 40 modes (black) and FOM (gray).}
\label{fig:ex2_sh2}
\end{figure}

Figure~\ref{fig:ex2_sh1} compares the Green's displacement $\tilde{u}^{\mathrm{G}}_{yy}(p_x, z, \omega)$ computed at the surface of the soil ($z=0 \, \mathrm{m}$) in the wavenumber–frequency domain using ROM and FOM. Figure~\ref{fig:ex2_sh1}b shows the full order Green’s function. The dispersion curves describe the wave modes generated by an out-of-plane Dirac excitation at the surface. At low wavenumbers and frequencies, the response is associated with modes that penetrate more deeply into the soil, whereas at higher wavenumbers and frequencies it is increasingly dominated by surface wave modes with limited penetration depth. The ROM (figure~\ref{fig:ex2_sh1}a) reproduces the main dispersion characteristics of the site, accurately capturing both the fundamental and higher-order Love modes with close agreement to the FOM. Figure~\ref{fig:ex2_sh1}c shows a slice for \(f = 7.5\,\mathrm{Hz}\), demonstrating close convergence between the ROM and FOM solutions and confirming that the reduced model reliably represents the full site response over the considered frequency range.
Figure~\ref{fig:ex2_sh2} presents the SH Green’s displacements $\hat{u}^{\mathrm{G}}_{yy}(x, z, \omega)$ in the spatial–frequency domain.  
Figures~\ref{fig:ex2_sh2}a and~\ref{fig:ex2_sh2}b show the displacement field computed using the ROM and FOM, respectively, while figure~\ref{fig:ex2_sh2}c plots the displacement variation along the surface at \(f = 10~\mathrm{Hz}\).  
The spatial response is characterized by oscillatory behavior predominantly concentrated in the upper part of the soil , with amplitudes decreasing with both depth and horizontal distance due to material damping.
The ROM with 40 modes accurately captures these spatial features, reproducing the full order response within prescribed accuracy.

\begin{figure}[!htb]
  \centering
  \includegraphics[width=0.35\linewidth]{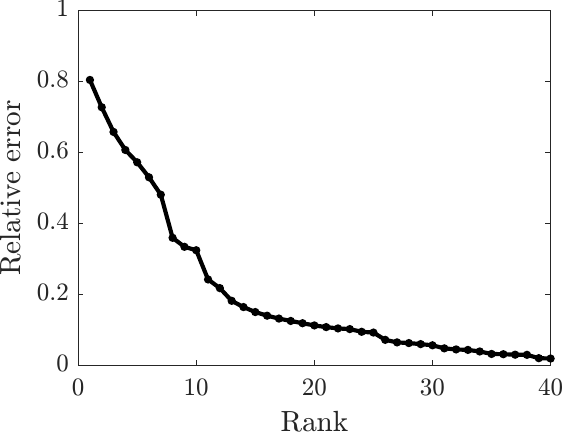}
  \caption{Relative Frobenius error of the ROM for the SH case as a function of the rank.}
  \label{fig:ex2_sh_error}
\end{figure}

To assess the convergence of the reduced order approximation, the relative Frobenius error between the full Green’s tensor and a sequence of Tucker reconstructions of increasing multilinear rank is evaluated. For each rank level, a truncated Tucker model is formed by retaining only the precribed number of modes in each parametric direction of the factor matrices and core tensor, yielding a hierarchy of ROMs with progressively richer separated representations. The error at each rank is computed as the Frobenius norm of the difference between the truncated ROM and the full order solution, normalized by the Frobenius norm of the latter:

\begin{equation}
\varepsilon_{\mathrm{F}}=
\frac{\left\lVert \tilde{\mathbf{u}}^{\mathrm{G}}_{\mathrm{FOM}}-
\tilde{\mathbf{u}}^{\mathrm{G}}_{\mathrm{ROM}(\mathbf{r})} \right\rVert_{\mathrm{F}}}
{\left\lVert \tilde{\mathbf{u}}^{\mathrm{G}}_{\mathrm{FOM}} \right\rVert_{\mathrm{F}}}
\label{eq:relFrob_error}
\end{equation}
where $\tilde{\mathbf{u}}^{\mathrm{G}}_{\mathrm{FOM}}$ denotes the full order Green’s tensor and $\tilde{\mathbf{u}}^{\mathrm{G}}_{\mathrm{ROM}(\mathbf{r})}$ the reduced order approximation obtained by truncating the Tucker decomposition
to multilinear rank $\mathbf{r} = (r_{x_\mathrm{s}}, r_{z_\mathrm{s}}, r_{p_x}, r_{z}, r_{\omega})$,
corresponding to the number of retained modes in each parametric direction. The error $\varepsilon_{\mathrm{F}}$ is evaluated for a sequence of increasing multilinear ranks $\mathbf{r}$, ranging from one up to the maximum number of retained modes in each
direction.
As shown in Figure~\ref{fig:ex2_sh_error}, increasing the multilinear rank systematically improves the approximation, with the final truncation level yielding a relative error of approximately \(2\%\) for the SH problem. 

A comparison of storage and computation time highlights the performance differences between ROM and FOM.
In the wavenumber-frequency domain, the FOM entails substantially higher storage requirements, while the ROM achieves a very compact representation. Although the offline construction of the ROM is more expensive than the corresponding FOM computation, the online reconstruction of the Green’s function is achieved efficiently through low-rank tensor contractions.
In the spatial–frequency domain, the FOM again requires significantly larger memory resources and computation time, whereas the ROM reconstruction is obtained at a markedly reduced cost. This acceleration primarily results from the reduced number of inverse Fourier transforms required in the reduced setting. These results highlight the trade-off between the higher offline cost of the ROM and its substantially reduced cost for subsequent evaluations. The corresponding storage requirements and computation times are summarized in tables~\ref{tab:storage_comparison_ex2_sh} and~\ref{tab:time_comparison_ex2_sh}. 

\begin{table}[!htb]
	\centering
	\caption{Storage requirement comparison between ROM and FOM for the SH case at the Groene Hart site.}
	\label{tab:storage_comparison_ex2_sh}
	\begin{tabular}{l c c c c}
		\toprule
		& \multicolumn{2}{c}{Wavenumber-frequency domain}
		& \multicolumn{2}{c}{Spatial-frequency domain} \\
		\cmidrule(lr){2-3} \cmidrule(lr){4-5}
		Case & ROM & FOM & ROM & FOM \\
		& [MB] & [MB] & [MB] & [MB] \\
		\midrule
		SH & 2.2 & 7300.0 & 62.0 & 1200.0 \\
		\bottomrule
	\end{tabular}
\end{table}

\begin{table}[!htb]
	\centering
	\caption{Computation time comparison between ROM and FOM for the SH case at the Groene Hart site.}
	\label{tab:time_comparison_ex2_sh}
	\begin{tabular}{l c c c c c c}
		\toprule
		& \multicolumn{3}{c}{Wavenumber-frequency domain}
		& \multicolumn{3}{c}{Spatial-frequency domain} \\
		\cmidrule(lr){2-4} \cmidrule(lr){5-7}
		Case 
		& ROM offline & ROM online & FOM
		& ROM offline & ROM online & FOM \\
		& [s] & [s] & [s] & [s] & [s] & [s] \\
		\midrule
		SH   
		& 5715.0  & 1.7  & 76.0 & 431.0  & 1.0  & 6276.0 \\
		\bottomrule
	\end{tabular}
\end{table}

\begin{figure}[!htb]
   \begin{center}
	 (a) \includegraphics[width=0.27\textwidth]{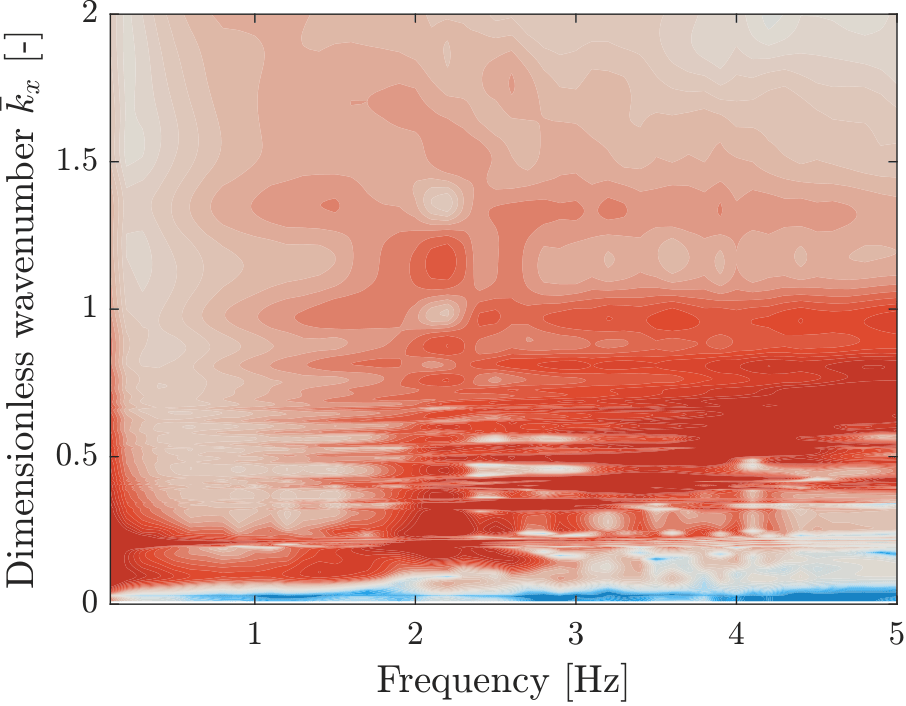}
	 (b) \includegraphics[width=0.27\textwidth]{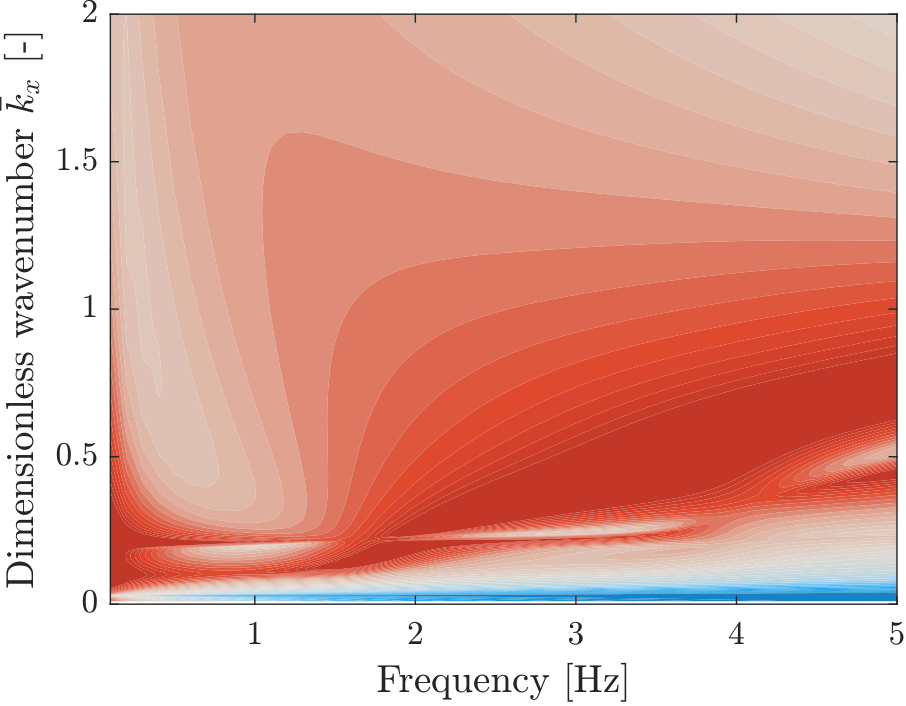} 
	 (c) \includegraphics[width=0.26\textwidth]{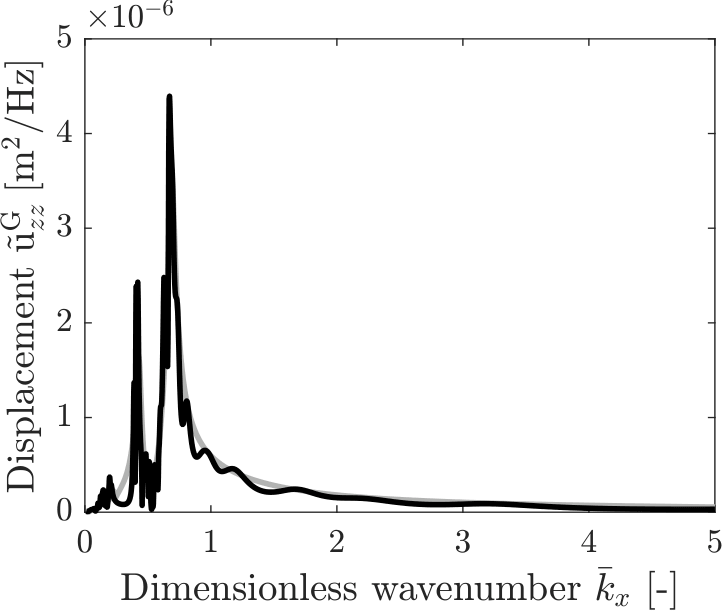}\\
     \vspace{0.2cm} \hspace{-4cm}
	     \includegraphics[width=0.27\textwidth]{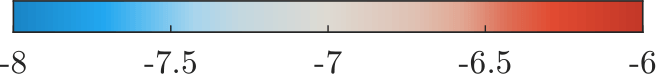}
   \end{center}
	 \caption{P-SV Green’s displacement $\tilde{u}^{\mathrm{G}}_{zz}(p_x, z, \omega)$ [$\mathrm{m}^2/\mathrm{Hz}$] in the wavenumber-frequency domain evaluated at $z=0 \, \mathrm{m}$ for a vertical concentrated load applied at $x_\mathrm{s}=0 \, \mathrm{m}$ and $z_\mathrm{s}=0 \, \mathrm{m}$ using (a) ROM with 40 modes, (b) FOM and (c) as a function of dimensionless wavenumber for $f = 5\, \mathrm{Hz}$ using ROM with 40 modes(black) and FOM (gray).}
\label{fig:ex2_psv1}
\end{figure}

For the P-SV case, the frequency range is limited to $f = 5~\mathrm{Hz}$ to account for the increased problem complexity and computational demand compared to the SH case. All remaining parameter sampling and the stopping criterion for the greedy algorithm are kept identical to those used in the SH case. Figure~\ref{fig:ex2_psv1}b presents the wavenumber–frequency response of the FOM evaluated at the soil surface, showing the characteristic dispersion branches of coupled P–SV motion dominated by the fundamental Rayleigh mode at low frequencies and higher-order modes at higher frequencies. Figure~\ref{fig:ex2_psv1}a shows the corresponding ROM response, which attempts to capture the main dispersion behavior but deviates in the higher branches due to the limited number of retained modes. The cross-sectional view in figure~\ref{fig:ex2_psv1}c, evaluated at $z=0 \, \mathrm{m}$ and  $f = 5\, \mathrm{Hz}$ further illustrates that while the ROM reproduces the overall displacement spectrum, finer-scale variations are not fully resolved.

\begin{figure}[!htb]
   \begin{center}
	 (a) \includegraphics[width=0.27\textwidth]{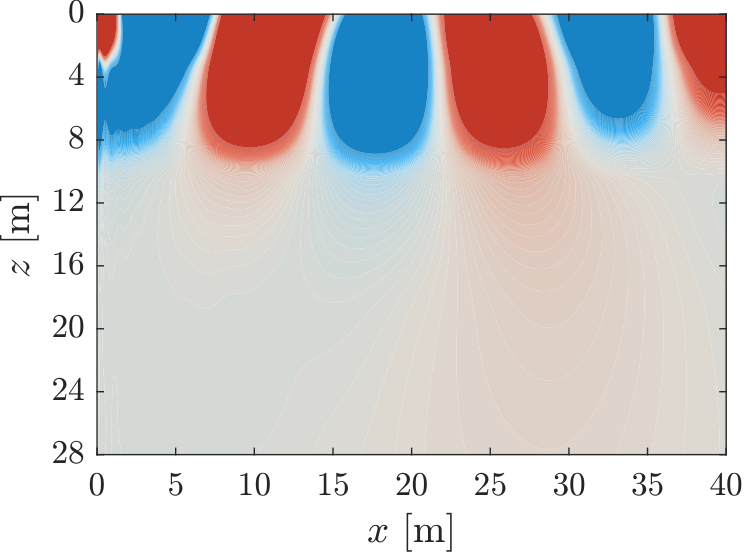}	
     (b) \includegraphics[width=0.27\textwidth]{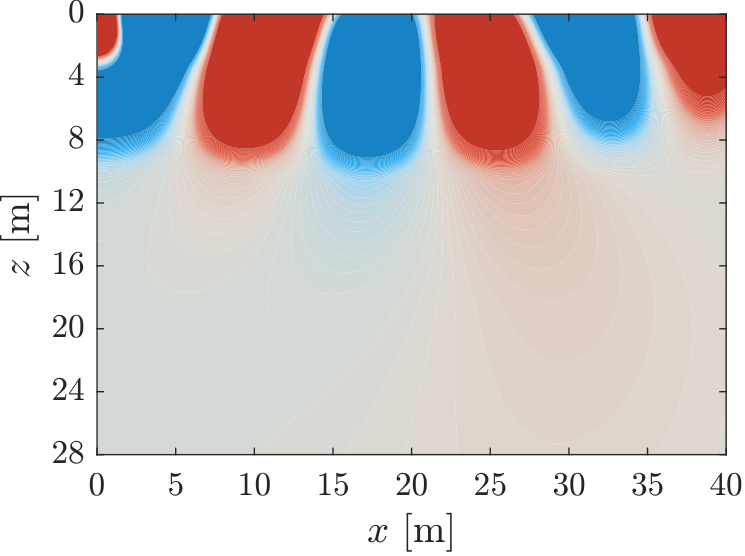}
     (c) \includegraphics[width=0.27\textwidth]{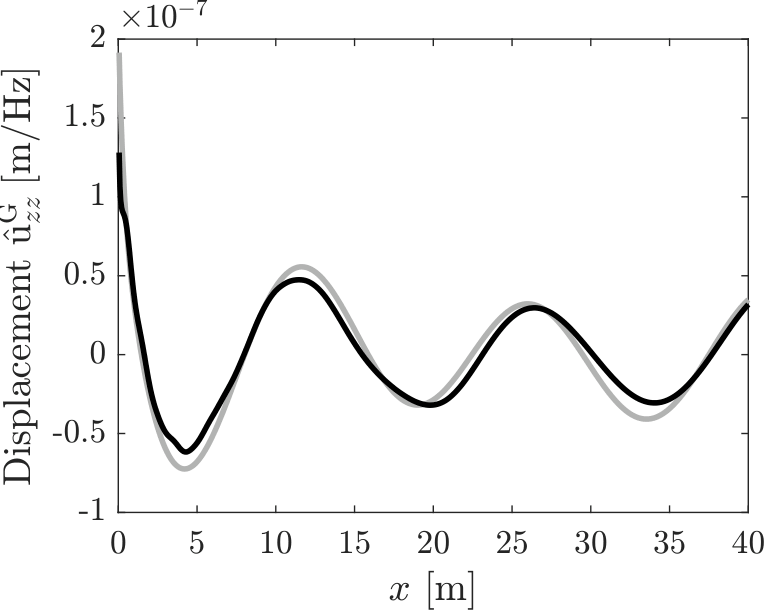}\\
     \vspace{0.2cm} \hspace{-4cm}
	     \includegraphics[width=0.26\textwidth]{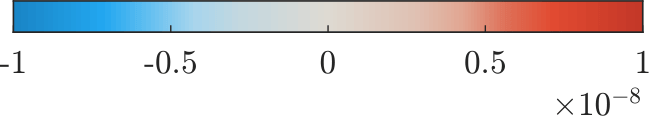}
   \end{center}
	 \caption{P-SV Green’s displacement $\hat{u}^{\mathrm{G}}_{zz}(x, z, \omega)$ [$\mathrm{m}/\mathrm{Hz}$] in the spatial-frequency domain evaluated at $f=5 \, \mathrm{Hz}$ for a vertical concentrated load applied at $x_\mathrm{s}=0 \, \mathrm{m}$ and $z_\mathrm{s}=0 \, \mathrm{m}$ using (a) ROM with 40 modes, (b) FOM and (c) as a function of $x$ evaluated at $z = 0\, \mathrm{m}$ using ROM with 40 modes (black) and FOM (gray).}
\label{fig:ex2_psv}
\end{figure}

In the spatial domain, figure~\ref{fig:ex2_psv}b shows the vertical displacement field, which is governed by Rayleigh wave propagation and decreases in amplitude with depth and distance from the source. The ROM with 40 modes (figure~\ref{fig:ex2_psv}a) captures the general interference pattern and attenuation trend, also visualized by the comparison in figure~\ref{fig:ex2_psv}c showing the displacement variation along the surface at $f = 5\, \mathrm{Hz}$. However, full convergence is not attained within 40 modes, as the coupled P–SV problem exhibits higher complexity and requires a greater number of modes for accurate representation than the simpler SH problem.

\begin{figure}[!htb]
  \centering
  \includegraphics[width=0.35\linewidth]{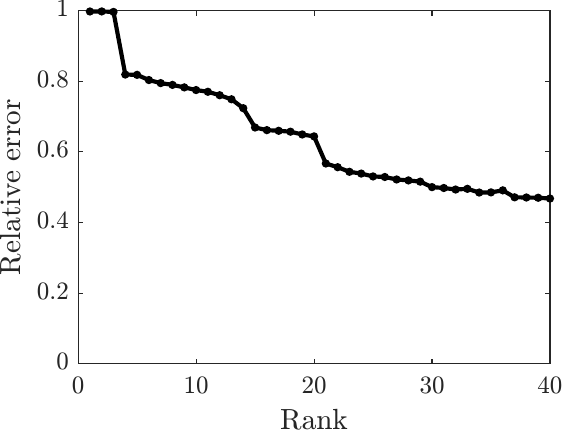}
  \caption{Relative Frobenius error of the ROM for the P-SV case as a function of the rank.}
  \label{fig:ex2_psv_error}
\end{figure}

The relative error of the ROM, shown in figure~\ref{fig:ex2_psv_error}, is approximately $45\%$. This reflects the slow convergence of the reduced basis in the P-SV case. With a rank of $40$, the ROM is able to capture the main trend of the displacement field, but the available modes are not sufficient to represent all variations across depth, wavenumber, and frequency. As a result, the ROM does not fully reproduce the complete displacement field, which leads to the observed discrepancy. Increasing the number of modes would improve the approximation accuracy and reduce the resulting error, but we are limited by computer memory.
For the current approximation of the Green's displacements, the storage requirement and computation time for ROM and FOM are summarized in tables~\ref{tab:storage_comparison_ex2} and~\ref{tab:time_comparison_ex2}, respectively. Achieving a higher level of accuracy would require additional modes, which is expected to further increase both storage demands and computational cost.

\begin{table}[!htb]
	\centering
	\caption{Storage requirement comparison between ROM and FOM for the P-SV case at the Groene Hart site.}
	\label{tab:storage_comparison_ex2}
	\begin{tabular}{l c c c c}
		\toprule
		& \multicolumn{2}{c}{Wavenumber-frequency domain}
		& \multicolumn{2}{c}{Spatial-frequency domain} \\
		\cmidrule(lr){2-3} \cmidrule(lr){4-5}
		Case & ROM & FOM & ROM & FOM \\
		& [MB] & [MB] & [MB] & [MB] \\
		\midrule
		P-SV & 2.7 & 10208.0 & 97.2 & 5444.0 \\
		\bottomrule
	\end{tabular}
\end{table}

\begin{table}[!htb]
	\centering
	\caption{Computation time comparison between ROM and FOM for the P-SV case at the Groene Hart site.}
	\label{tab:time_comparison_ex2}
	\begin{tabular}{l c c c c c c}
		\toprule
		& \multicolumn{3}{c}{Wavenumber-frequency domain}
		& \multicolumn{3}{c}{Spatial-frequency domain} \\
		\cmidrule(lr){2-4} \cmidrule(lr){5-7}
		Case 
		& ROM offline & ROM online & FOM
		& ROM offline & ROM online & FOM \\
		& [s] & [s] & [s] & [s] & [s] & [s] \\
		\midrule
        P-SV   
		& 10181.0  & 1.8  & 82.3 & 214.2  & 3.0  & 12471.0 \\
		\bottomrule
	\end{tabular}
\end{table}

\section{Conclusion}
\label{sec:conclusions}
The presented study demonstrates that low-rank tensor decompositions, in particular the Greedy Tucker Approximation, enable efficient and accurate reduced order modeling of elastodynamic Green’s functions. The proposed formulation achieves substantial reductions in computational cost and storage requirements while preserving close agreement with full order solutions across a range of frequencies and wave types. This efficiency makes the approach well suited for large-scale, high-dimensional simulations where conventional methods become prohibitive. A notable limitation, however, is the growth of the core tensor size with increasing dimensionality of the parameter space. While the separated representation significantly reduces the size of the factor matrices, the core tensor scales with the product of the modal ranks in each dimension. In the present work, separability was introduced selectively while other variables were kept fixed, resulting in moderate ranks. However, in a fully parametric setting where spatial coordinates, frequency, wavenumber, and material properties are treated simultaneously as separable dimensions, the multilinear ranks may increase. This can lead to a significant growth of the core tensor, which may limit scalability due to increased memory requirements and computational cost.  To address this, more advanced tensor formats such as the tensor train \cite{lubi13a} or the hierarchical Tucker format \cite{kres14a} are expected to provide improved scalability, and their integration into the current framework is the subject of ongoing work.
Beyond computational efficiency, the methodology presented here opens pathways for incorporating reduced order Green’s function representations into diverse simulation frameworks, including boundary element methods, structural dynamic simulations, and parametric analyses. Such extensions have the potential to significantly accelerate the analysis of wave propagation and soil–structure interaction problems, enabling more complex geometries, material heterogeneities and parametric studies to be tackled within feasible computational budgets.

\section*{Acknowledgment}
This research was performed within the frame of the project C14/20/073 “Model order reduction techniques for the prediction of railway induced vibration in the built environment” funded by the Research Council of KU Leuven. The financial support is gratefully acknowledged.

\appendix
\section{Layer stiffness matrices}
\label{sec:appendix1}
The expressions for layer stiffness matrices used in equation~\ref{eq:Ks}:
\begin{eqnarray}
	\mathbf{A}^\mathrm{e}&=&\frac{h}{6}
	\begin{bmatrix}
	2(\lambda+2\mu) & 0 & \lambda+2\mu & 0 \\
	0 & 2\mu & 0 & \mu \\
	\lambda+2\mu & 0 & 2(\lambda+2\mu) & 0 \\
	0 & \mu & 0 & 2\mu \\
	\end{bmatrix} \\
	\mathbf{B}^\mathrm{e}&=&\frac{1}{2}
	\begin{bmatrix}
	0 & \lambda-\mu & 0 & -(\lambda+\mu) \\
	\mu-\lambda & 0 & -(\lambda+\mu) & 0 \\
	0 & \lambda+\mu & 0 & \mu-\lambda \\
	\lambda+\mu & 0 & \lambda-\mu & 0 \\
	\end{bmatrix} \\
	\mathbf{G}^\mathrm{e}&=&\frac{1}{h}
	\begin{bmatrix}
	+\mu & 0 & -\mu & 0 \\
	0 & +(\lambda+2\mu) & 0 & -(\lambda+2\mu) \\
	-\mu & 0 & +\mu & 0 \\
	0 & -(\lambda+2\mu) & 0 & +(\lambda+2\mu) \\
	\end{bmatrix} \\
	\mathbf{M}^\mathrm{e}&=&\frac{\rho h}{6}
	\begin{bmatrix}
	2 & 0 & 1 & 0 \\
	0 & 2 & 0 & 1 \\
	1 & 0 & 2 & 0 \\
	0 & 1 & 0 & 2 \\
	\end{bmatrix}
\end{eqnarray}
where $\lambda$ is the Lame constant, $\mu$ is the shear modulus, $\rho$ is the mass density and $h$ is the layer thickness.

\section{Expression for $\boldsymbol{\mathcal{A}}$}
\label{sec:appendix2}
The dynamic stiffness tensor $\boldsymbol{\mathcal{A}}$ in equation~\eqref{eq:weakPGD} is equal to:
\begin{multline}
\mathcal{A} = \sum_{r_a=1}^{R_a-1} 
\mathbf{A}_{r_a}^{(x_\mathrm{s})} \otimes 
\mathbf{A}_{r_a}^{(z_\mathrm{s})} \otimes 
\mathbf{A}_{r_a}^{(p_x)} \otimes 
\mathbf{A}_{r_a}^{(z)} \otimes 
\mathbf{A}_{r_a}^{(\omega)} \\ +
\int\limits_{\mathcal{I}_{x_\mathrm{s}}}
\int\limits_{\mathcal{I}_{z_\mathrm{s}}}
\int\limits_{\mathcal{I}_{p_x}}
\int\limits_{\mathcal{I}_{\omega}}
\omega
\left(
\mathbf{N}_{x_{\mathrm{s}}} \otimes
\mathbf{N}_{z_{\mathrm{s}}} \otimes
\mathbf{N}_{p_x} \otimes
\mathbf{N}_{z = h} \otimes
\mathbf{N}_{\omega}
\right)^\mathrm{T}
\mathbf{K}_{p_x}(p_x)
\left(
\mathbf{N}_{x_{\mathrm{s}}} \otimes
\mathbf{N}_{z_{\mathrm{s}}} \otimes
\mathbf{N}_{p_x} \otimes
\mathbf{N}_{z = h} \otimes
\mathbf{N}_{\omega}
\right)
\mathrm{d}x_{\mathrm{s}}
\mathrm{d}z_{\mathrm{s}}
\mathrm{d}p_x
\mathrm{d}z
\mathrm{d}\omega
\label{A}
\end{multline}
where $R_a$ is the rank of $\boldsymbol{\mathcal{A}}$, and $\mathbf{A}_{r_a}^{(x_\mathrm{s})}$, $\mathbf{A}_{r_a}^{(z_\mathrm{s})}$, $\mathbf{A}_{r_a}^{(p_x)}$, $\mathbf{A}_{r_a}^{(z)}$ and $\mathbf{A}_{r_a}^{(\omega)}$ are determined from the separation of the integrals on the domain $\Omega$ in the weak form~\eqref{eq:weak6}, and are equal to:

{
\setlength{\jot}{1.0em}
\begin{flalign*}
&\mathbf{A}_{1}^{(x_{\mathrm{s}})} = \int\limits_{\mathcal{I}_{x_{\mathrm{s}}}}\mathbf{N}_{x_{\mathrm{s}}}^{\mathrm{T}}\mathbf{N}_{x_{\mathrm{s}}}\,\mathrm{d}x_{\mathrm{s}} &\\
&\mathbf{A}_{1}^{(z_{\mathrm{s}})} = \int\limits_{\mathcal{I}_{z_{\mathrm{s}}}}\mathbf{N}_{z_{\mathrm{s}}}^{\mathrm{T}}\mathbf{N}_{z_{\mathrm{s}}}\,\mathrm{d}z_{\mathrm{s}} &\\
&\mathbf{A}_{1}^{(p_x)} = \int\limits_{\mathcal{I}_{p_x}}p_x^{2}\,\mathbf{N}_{p_x}^{\mathrm{T}}\mathbf{N}_{p_x}\,\mathrm{d}p_x &\\
&\mathbf{A}_{1}^{(z)} = \int\limits_{\mathcal{I}_z}\mathbf{N}_{z}^{\mathrm{T}}\mathbf{L}_{x}^{\mathrm{T}}\mathbf{C}\mathbf{L}_{x}\mathbf{N}_{z}\,\mathrm{d}z &\\
&\mathbf{A}_{1}^{(\omega)} = \int\limits_{\mathcal{I}_{\omega}}\omega^{2}\,\mathbf{N}_{\omega}^{\mathrm{T}}\mathbf{N}_{\omega}\,\mathrm{d}\omega &
\end{flalign*}
}
{
\setlength{\jot}{1.0em}
\begin{flalign*}
&\mathbf{A}_{2}^{(x_{\mathrm{s}})} = \mathbf{A}_{1}^{(x_{\mathrm{s}})} &\\
&\mathbf{A}_{2}^{(z_{\mathrm{s}})} = \mathbf{A}_{1}^{(z_{\mathrm{s}})} &\\
&\mathbf{A}_{2}^{(p_x)} = \int\limits_{\mathcal{I}_{p_x}}p_x\,\mathbf{N}_{p_x}^{\mathrm{T}}\mathbf{N}_{p_x}\,\mathrm{d}p_x &\\
&\mathbf{A}_{2}^{(z)} = \int\limits_{\mathcal{I}_z}\Big(
\mathbf{N}_{z}^{\mathrm{T}}\mathbf{L}_{x}^{\mathrm{T}}\mathbf{C}\mathbf{L}_{z}\frac{\partial\mathbf{N}_{z}}{\partial z}
-\frac{\partial\mathbf{N}_{z}^{\mathrm{T}}}{\partial z}\mathbf{L}_{z}\mathbf{C}\mathbf{L}_{x}\mathbf{N}_{z}
\Big)\mathrm{d}z &\\
&\mathbf{A}_{2}^{(\omega)} = \int\limits_{\mathcal{I}_{\omega}}\omega\,\mathbf{N}_{\omega}^{\mathrm{T}}\mathbf{N}_{\omega}\,\mathrm{d}\omega &
\end{flalign*}
}
{
\setlength{\jot}{1.0em}
\begin{flalign*}
&\mathbf{A}_{3}^{(x_{\mathrm{s}})} = \mathbf{A}_{1}^{(x_{\mathrm{s}})} &\\
&\mathbf{A}_{3}^{(z_{\mathrm{s}})} = \mathbf{A}_{1}^{(z_{\mathrm{s}})} &\\
&\mathbf{A}_{3}^{(p_x)} = \int\limits_{\mathcal{I}_{p_x}}\mathbf{N}_{p_x}^{\mathrm{T}}\mathbf{N}_{p_x}\,\mathrm{d}p_x &\\
&\mathbf{A}_{3}^{(z)} = \int\limits_{\mathcal{I}_z}\frac{\partial\mathbf{N}_{z}^{\mathrm{T}}}{\partial z}\mathbf{L}_{z}^{\mathrm{T}}\mathbf{C}\mathbf{L}_{z}\frac{\partial\mathbf{N}_{z}}{\partial z}\,\mathrm{d}z &\\
&\mathbf{A}_{3}^{(\omega)} = \int\limits_{\mathcal{I}_{\omega}}\mathbf{N}_{\omega}^{\mathrm{T}}\mathbf{N}_{\omega}\,\mathrm{d}\omega &
\end{flalign*}
}
{
\setlength{\jot}{1.0em}
\begin{flalign*}
&\mathbf{A}_{4}^{(x_{\mathrm{s}})} = \mathbf{A}_{1}^{(x_{\mathrm{s}})} &\\
&\mathbf{A}_{4}^{(z_{\mathrm{s}})} = \mathbf{A}_{1}^{(z_{\mathrm{s}})} &\\
&\mathbf{A}_{4}^{(p_x)} = \mathbf{A}_{3}^{(p_x)} &\\
&\mathbf{A}_{4}^{(z)} = \int\limits_{\mathcal{I}_z}\rho\,\mathbf{N}_{z}^{\mathrm{T}}\mathbf{N}_{z}\,\mathrm{d}z &\\
&\mathbf{A}_{4}^{(\omega)} = -\mathbf{A}_{1}^{(\omega)} &
\end{flalign*}
}

\end{document}